
\catcode`\@=11


\message{Loading jyTeX fonts...}



\font\vptrm=cmr5 \font\vptmit=cmmi5 \font\vptsy=cmsy5 \font\vptbf=cmbx5

\skewchar\vptmit='177 \skewchar\vptsy='60 \fontdimen16
\vptsy=\the\fontdimen17 \vptsy

\def\vpt{\ifmmode\err@badsizechange\else
     \@mathfontinit
     \textfont0=\vptrm  \scriptfont0=\vptrm  \scriptscriptfont0=\vptrm
     \textfont1=\vptmit \scriptfont1=\vptmit \scriptscriptfont1=\vptmit
     \textfont2=\vptsy  \scriptfont2=\vptsy  \scriptscriptfont2=\vptsy
     \textfont3=\xptex  \scriptfont3=\xptex  \scriptscriptfont3=\xptex
     \textfont\bffam=\vptbf
     \scriptfont\bffam=\vptbf
     \scriptscriptfont\bffam=\vptbf
     \@fontstyleinit
     \def\rm{\vptrm\fam=\z@}%
     \def\bf{\vptbf\fam=\bffam}%
     \def\oldstyle{\vptmit\fam=\@ne}%
     \rm\fi}


\font\viptrm=cmr6 \font\viptmit=cmmi6 \font\viptsy=cmsy6
\font\viptbf=cmbx6

\skewchar\viptmit='177 \skewchar\viptsy='60 \fontdimen16
\viptsy=\the\fontdimen17 \viptsy

\def\vipt{\ifmmode\err@badsizechange\else
     \@mathfontinit
     \textfont0=\viptrm  \scriptfont0=\vptrm  \scriptscriptfont0=\vptrm
     \textfont1=\viptmit \scriptfont1=\vptmit \scriptscriptfont1=\vptmit
     \textfont2=\viptsy  \scriptfont2=\vptsy  \scriptscriptfont2=\vptsy
     \textfont3=\xptex   \scriptfont3=\xptex  \scriptscriptfont3=\xptex
     \textfont\bffam=\viptbf
     \scriptfont\bffam=\vptbf
     \scriptscriptfont\bffam=\vptbf
     \@fontstyleinit
     \def\rm{\viptrm\fam=\z@}%
     \def\bf{\viptbf\fam=\bffam}%
     \def\oldstyle{\viptmit\fam=\@ne}%
     \rm\fi}

\font\viiptrm=cmr7 \font\viiptmit=cmmi7 \font\viiptsy=cmsy7
\font\viiptit=cmti7 \font\viiptbf=cmbx7

\skewchar\viiptmit='177 \skewchar\viiptsy='60 \fontdimen16
\viiptsy=\the\fontdimen17 \viiptsy

\def\viipt{\ifmmode\err@badsizechange\else
     \@mathfontinit
     \textfont0=\viiptrm  \scriptfont0=\vptrm  \scriptscriptfont0=\vptrm
     \textfont1=\viiptmit \scriptfont1=\vptmit \scriptscriptfont1=\vptmit
     \textfont2=\viiptsy  \scriptfont2=\vptsy  \scriptscriptfont2=\vptsy
     \textfont3=\xptex    \scriptfont3=\xptex  \scriptscriptfont3=\xptex
     \textfont\itfam=\viiptit
     \scriptfont\itfam=\viiptit
     \scriptscriptfont\itfam=\viiptit
     \textfont\bffam=\viiptbf
     \scriptfont\bffam=\vptbf
     \scriptscriptfont\bffam=\vptbf
     \@fontstyleinit
     \def\rm{\viiptrm\fam=\z@}%
     \def\it{\viiptit\fam=\itfam}%
     \def\bf{\viiptbf\fam=\bffam}%
     \def\oldstyle{\viiptmit\fam=\@ne}%
     \rm\fi}


\font\viiiptrm=cmr8 \font\viiiptmit=cmmi8 \font\viiiptsy=cmsy8
\font\viiiptit=cmti8
\font\viiiptbf=cmbx8

\skewchar\viiiptmit='177 \skewchar\viiiptsy='60 \fontdimen16
\viiiptsy=\the\fontdimen17 \viiiptsy

\def\viiipt{\ifmmode\err@badsizechange\else
     \@mathfontinit
     \textfont0=\viiiptrm  \scriptfont0=\viptrm  \scriptscriptfont0=\vptrm
     \textfont1=\viiiptmit \scriptfont1=\viptmit \scriptscriptfont1=\vptmit
     \textfont2=\viiiptsy  \scriptfont2=\viptsy  \scriptscriptfont2=\vptsy
     \textfont3=\xptex     \scriptfont3=\xptex   \scriptscriptfont3=\xptex
     \textfont\itfam=\viiiptit
     \scriptfont\itfam=\viiptit
     \scriptscriptfont\itfam=\viiptit
     \textfont\bffam=\viiiptbf
     \scriptfont\bffam=\viptbf
     \scriptscriptfont\bffam=\vptbf
     \@fontstyleinit
     \def\rm{\viiiptrm\fam=\z@}%
     \def\it{\viiiptit\fam=\itfam}%
     \def\bf{\viiiptbf\fam=\bffam}%
     \def\oldstyle{\viiiptmit\fam=\@ne}%
     \rm\fi}


\def\getixpt{%
     \font\ixptrm=cmr9
     \font\ixptmit=cmmi9
     \font\ixptsy=cmsy9
     \font\ixptit=cmti9
     \font\ixptbf=cmbx9
     \skewchar\ixptmit='177 \skewchar\ixptsy='60
     \fontdimen16 \ixptsy=\the\fontdimen17 \ixptsy}

\def\ixpt{\ifmmode\err@badsizechange\else
     \@mathfontinit
     \textfont0=\ixptrm  \scriptfont0=\viiptrm  \scriptscriptfont0=\vptrm
     \textfont1=\ixptmit \scriptfont1=\viiptmit \scriptscriptfont1=\vptmit
     \textfont2=\ixptsy  \scriptfont2=\viiptsy  \scriptscriptfont2=\vptsy
     \textfont3=\xptex   \scriptfont3=\xptex    \scriptscriptfont3=\xptex
     \textfont\itfam=\ixptit
     \scriptfont\itfam=\viiptit
     \scriptscriptfont\itfam=\viiptit
     \textfont\bffam=\ixptbf
     \scriptfont\bffam=\viiptbf
     \scriptscriptfont\bffam=\vptbf
     \@fontstyleinit
     \def\rm{\ixptrm\fam=\z@}%
     \def\it{\ixptit\fam=\itfam}%
     \def\bf{\ixptbf\fam=\bffam}%
     \def\oldstyle{\ixptmit\fam=\@ne}%
     \rm\fi}


\font\xptrm=cmr10 \font\xptmit=cmmi10 \font\xptsy=cmsy10
\font\xptex=cmex10 \font\xptit=cmti10 \font\xptsl=cmsl10
\font\xptbf=cmbx10 \font\xpttt=cmtt10 \font\xptss=cmss10
\font\xptsc=cmcsc10 \font\xptbfs=cmb10 \font\xptbmit=cmmib10

\skewchar\xptmit='177 \skewchar\xptbmit='177 \skewchar\xptsy='60
\fontdimen16 \xptsy=\the\fontdimen17 \xptsy

\def\xpt{\ifmmode\err@badsizechange\else
     \@mathfontinit
     \textfont0=\xptrm  \scriptfont0=\viiptrm  \scriptscriptfont0=\vptrm
     \textfont1=\xptmit \scriptfont1=\viiptmit \scriptscriptfont1=\vptmit
     \textfont2=\xptsy  \scriptfont2=\viiptsy  \scriptscriptfont2=\vptsy
     \textfont3=\xptex  \scriptfont3=\xptex    \scriptscriptfont3=\xptex
     \textfont\itfam=\xptit
     \scriptfont\itfam=\viiptit
     \scriptscriptfont\itfam=\viiptit
     \textfont\bffam=\xptbf
     \scriptfont\bffam=\viiptbf
     \scriptscriptfont\bffam=\vptbf
     \textfont\bfsfam=\xptbfs
     \scriptfont\bfsfam=\viiptbf
     \scriptscriptfont\bfsfam=\vptbf
     \textfont\bmitfam=\xptbmit
     \scriptfont\bmitfam=\viiptmit
     \scriptscriptfont\bmitfam=\vptmit
     \@fontstyleinit
     \def\rm{\xptrm\fam=\z@}%
     \def\it{\xptit\fam=\itfam}%
     \def\sl{\xptsl}%
     \def\bf{\xptbf\fam=\bffam}%
     \def\tt{\xpttt}%
     \def\ss{\xptss}%
     \def\sc{\xptsc}%
     \def\bfs{\xptbfs\fam=\bfsfam}%
     \def\bmit{\fam=\bmitfam}%
     \def\oldstyle{\xptmit\fam=\@ne}%
     \rm\fi}


\def\getxipt{%
     \font\xiptrm=cmr10  scaled\magstephalf
     \font\xiptmit=cmmi10 scaled\magstephalf
     \font\xiptsy=cmsy10 scaled\magstephalf
     \font\xiptex=cmex10 scaled\magstephalf
     \font\xiptit=cmti10 scaled\magstephalf
     \font\xiptsl=cmsl10 scaled\magstephalf
     \font\xiptbf=cmbx10 scaled\magstephalf
     \font\xipttt=cmtt10 scaled\magstephalf
     \font\xiptss=cmss10 scaled\magstephalf
     \skewchar\xiptmit='177 \skewchar\xiptsy='60
     \fontdimen16 \xiptsy=\the\fontdimen17 \xiptsy}

\def\xipt{\ifmmode\err@badsizechange\else
     \@mathfontinit
     \textfont0=\xiptrm  \scriptfont0=\viiiptrm  \scriptscriptfont0=\viptrm
     \textfont1=\xiptmit \scriptfont1=\viiiptmit \scriptscriptfont1=\viptmit
     \textfont2=\xiptsy  \scriptfont2=\viiiptsy  \scriptscriptfont2=\viptsy
     \textfont3=\xiptex  \scriptfont3=\xptex     \scriptscriptfont3=\xptex
     \textfont\itfam=\xiptit
     \scriptfont\itfam=\viiiptit
     \scriptscriptfont\itfam=\viiptit
     \textfont\bffam=\xiptbf
     \scriptfont\bffam=\viiiptbf
     \scriptscriptfont\bffam=\viptbf
     \@fontstyleinit
     \def\rm{\xiptrm\fam=\z@}%
     \def\it{\xiptit\fam=\itfam}%
     \def\sl{\xiptsl}%
     \def\bf{\xiptbf\fam=\bffam}%
     \def\tt{\xipttt}%
     \def\ss{\xiptss}%
     \def\oldstyle{\xiptmit\fam=\@ne}%
     \rm\fi}


\font\xiiptrm=cmr12 \font\xiiptmit=cmmi12 \font\xiiptsy=cmsy10
scaled\magstep1 \font\xiiptex=cmex10  scaled\magstep1 \font\xiiptit=cmti12
\font\xiiptsl=cmsl12 \font\xiiptbf=cmbx12
\font\xiiptss=cmss12 \font\xiiptsc=cmcsc10 scaled\magstep1
\font\xiiptbfs=cmb10  scaled\magstep1 \font\xiiptbmit=cmmib10
scaled\magstep1

\skewchar\xiiptmit='177 \skewchar\xiiptbmit='177 \skewchar\xiiptsy='60
\fontdimen16 \xiiptsy=\the\fontdimen17 \xiiptsy

\def\xiipt{\ifmmode\err@badsizechange\else
     \@mathfontinit
     \textfont0=\xiiptrm  \scriptfont0=\viiiptrm  \scriptscriptfont0=\viptrm
     \textfont1=\xiiptmit \scriptfont1=\viiiptmit \scriptscriptfont1=\viptmit
     \textfont2=\xiiptsy  \scriptfont2=\viiiptsy  \scriptscriptfont2=\viptsy
     \textfont3=\xiiptex  \scriptfont3=\xptex     \scriptscriptfont3=\xptex
     \textfont\itfam=\xiiptit
     \scriptfont\itfam=\viiiptit
     \scriptscriptfont\itfam=\viiptit
     \textfont\bffam=\xiiptbf
     \scriptfont\bffam=\viiiptbf
     \scriptscriptfont\bffam=\viptbf
     \textfont\bfsfam=\xiiptbfs
     \scriptfont\bfsfam=\viiiptbf
     \scriptscriptfont\bfsfam=\viptbf
     \textfont\bmitfam=\xiiptbmit
     \scriptfont\bmitfam=\viiiptmit
     \scriptscriptfont\bmitfam=\viptmit
     \@fontstyleinit
     \def\rm{\xiiptrm\fam=\z@}%
     \def\it{\xiiptit\fam=\itfam}%
     \def\sl{\xiiptsl}%
     \def\bf{\xiiptbf\fam=\bffam}%
     \def\tt{\xiipttt}%
     \def\ss{\xiiptss}%
     \def\sc{\xiiptsc}%
     \def\bfs{\xiiptbfs\fam=\bfsfam}%
     \def\bmit{\fam=\bmitfam}%
     \def\oldstyle{\xiiptmit\fam=\@ne}%
     \rm\fi}


\def\getxiiipt{%
     \font\xiiiptrm=cmr12  scaled\magstephalf
     \font\xiiiptmit=cmmi12 scaled\magstephalf
     \font\xiiiptsy=cmsy9  scaled\magstep2
     \font\xiiiptit=cmti12 scaled\magstephalf
     \font\xiiiptsl=cmsl12 scaled\magstephalf
     \font\xiiiptbf=cmbx12 scaled\magstephalf
     \font\xiiipttt=cmtt12 scaled\magstephalf
     \font\xiiiptss=cmss12 scaled\magstephalf
     \skewchar\xiiiptmit='177 \skewchar\xiiiptsy='60
     \fontdimen16 \xiiiptsy=\the\fontdimen17 \xiiiptsy}

\def\xiiipt{\ifmmode\err@badsizechange\else
     \@mathfontinit
     \textfont0=\xiiiptrm  \scriptfont0=\xptrm  \scriptscriptfont0=\viiptrm
     \textfont1=\xiiiptmit \scriptfont1=\xptmit \scriptscriptfont1=\viiptmit
     \textfont2=\xiiiptsy  \scriptfont2=\xptsy  \scriptscriptfont2=\viiptsy
     \textfont3=\xivptex   \scriptfont3=\xptex  \scriptscriptfont3=\xptex
     \textfont\itfam=\xiiiptit
     \scriptfont\itfam=\xptit
     \scriptscriptfont\itfam=\viiptit
     \textfont\bffam=\xiiiptbf
     \scriptfont\bffam=\xptbf
     \scriptscriptfont\bffam=\viiptbf
     \@fontstyleinit
     \def\rm{\xiiiptrm\fam=\z@}%
     \def\it{\xiiiptit\fam=\itfam}%
     \def\sl{\xiiiptsl}%
     \def\bf{\xiiiptbf\fam=\bffam}%
     \def\tt{\xiiipttt}%
     \def\ss{\xiiiptss}%
     \def\oldstyle{\xiiiptmit\fam=\@ne}%
     \rm\fi}


\font\xivptrm=cmr12   scaled\magstep1 \font\xivptmit=cmmi12
scaled\magstep1 \font\xivptsy=cmsy10  scaled\magstep2
\font\xivptex=cmex10  scaled\magstep2 \font\xivptit=cmti12 scaled\magstep1
\font\xivptsl=cmsl12  scaled\magstep1 \font\xivptbf=cmbx12
scaled\magstep1
\font\xivptss=cmss12  scaled\magstep1 \font\xivptsc=cmcsc10
scaled\magstep2 \font\xivptbfs=cmb10  scaled\magstep2
\font\xivptbmit=cmmib10 scaled\magstep2

\skewchar\xivptmit='177 \skewchar\xivptbmit='177 \skewchar\xivptsy='60
\fontdimen16 \xivptsy=\the\fontdimen17 \xivptsy

\def\xivpt{\ifmmode\err@badsizechange\else
     \@mathfontinit
     \textfont0=\xivptrm  \scriptfont0=\xptrm  \scriptscriptfont0=\viiptrm
     \textfont1=\xivptmit \scriptfont1=\xptmit \scriptscriptfont1=\viiptmit
     \textfont2=\xivptsy  \scriptfont2=\xptsy  \scriptscriptfont2=\viiptsy
     \textfont3=\xivptex  \scriptfont3=\xptex  \scriptscriptfont3=\xptex
     \textfont\itfam=\xivptit
     \scriptfont\itfam=\xptit
     \scriptscriptfont\itfam=\viiptit
     \textfont\bffam=\xivptbf
     \scriptfont\bffam=\xptbf
     \scriptscriptfont\bffam=\viiptbf
     \textfont\bfsfam=\xivptbfs
     \scriptfont\bfsfam=\xptbfs
     \scriptscriptfont\bfsfam=\viiptbf
     \textfont\bmitfam=\xivptbmit
     \scriptfont\bmitfam=\xptbmit
     \scriptscriptfont\bmitfam=\viiptmit
     \@fontstyleinit
     \def\rm{\xivptrm\fam=\z@}%
     \def\it{\xivptit\fam=\itfam}%
     \def\sl{\xivptsl}%
     \def\bf{\xivptbf\fam=\bffam}%
     \def\tt{\xivpttt}%
     \def\ss{\xivptss}%
     \def\sc{\xivptsc}%
     \def\bfs{\xivptbfs\fam=\bfsfam}%
     \def\bmit{\fam=\bmitfam}%
     \def\oldstyle{\xivptmit\fam=\@ne}%
     \rm\fi}


\font\xviiptrm=cmr17 \font\xviiptmit=cmmi12 scaled\magstep2
\font\xviiptsy=cmsy10 scaled\magstep3 \font\xviiptex=cmex10
scaled\magstep3 \font\xviiptit=cmti12 scaled\magstep2
\font\xviiptbf=cmbx12 scaled\magstep2 \font\xviiptbfs=cmb10
scaled\magstep3

\skewchar\xviiptmit='177 \skewchar\xviiptsy='60 \fontdimen16
\xviiptsy=\the\fontdimen17 \xviiptsy

\def\xviipt{\ifmmode\err@badsizechange\else
     \@mathfontinit
     \textfont0=\xviiptrm  \scriptfont0=\xiiptrm  \scriptscriptfont0=\viiiptrm
     \textfont1=\xviiptmit \scriptfont1=\xiiptmit \scriptscriptfont1=\viiiptmit
     \textfont2=\xviiptsy  \scriptfont2=\xiiptsy  \scriptscriptfont2=\viiiptsy
     \textfont3=\xviiptex  \scriptfont3=\xiiptex  \scriptscriptfont3=\xptex
     \textfont\itfam=\xviiptit
     \scriptfont\itfam=\xiiptit
     \scriptscriptfont\itfam=\viiiptit
     \textfont\bffam=\xviiptbf
     \scriptfont\bffam=\xiiptbf
     \scriptscriptfont\bffam=\viiiptbf
     \textfont\bfsfam=\xviiptbfs
     \scriptfont\bfsfam=\xiiptbfs
     \scriptscriptfont\bfsfam=\viiiptbf
     \@fontstyleinit
     \def\rm{\xviiptrm\fam=\z@}%
     \def\it{\xviiptit\fam=\itfam}%
     \def\bf{\xviiptbf\fam=\bffam}%
     \def\bfs{\xviiptbfs\fam=\bfsfam}%
     \def\oldstyle{\xviiptmit\fam=\@ne}%
     \rm\fi}


\font\xxiptrm=cmr17  scaled\magstep1


\def\xxipt{\ifmmode\err@badsizechange\else
     \@mathfontinit
     \@fontstyleinit
     \def\rm{\xxiptrm\fam=\z@}%
     \rm\fi}


\font\xxvptrm=cmr17  scaled\magstep2


\def\xxvpt{\ifmmode\err@badsizechange\else
     \@mathfontinit
     \@fontstyleinit
     \def\rm{\xxvptrm\fam=\z@}%
     \rm\fi}




\message{Loading jyTeX macros...}

\message{modifications to plain.tex,}


\def\newcount{\alloc@0\count\countdef\insc@unt}
\def\newdimen{\alloc@1\dimen\dimendef\insc@unt}
\def\newskip{\alloc@2\skip\skipdef\insc@unt}
\def\newmuskip{\alloc@3\muskip\muskipdef\@cclvi}
\def\newbox{\alloc@4\box\chardef\insc@unt}
\def\newtoks{\alloc@5\toks\toksdef\@cclvi}
\def\newhelp#1#2{\newtoks#1\global#1\expandafter{\csname#2\endcsname}}
\def\newread{\alloc@6\read\chardef\sixt@@n}
\def\newwrite{\alloc@7\write\chardef\sixt@@n}
\def\newfam{\alloc@8\fam\chardef\sixt@@n}
\def\newinsert#1{\global\advance\insc@unt by\m@ne
     \ch@ck0\insc@unt\count
     \ch@ck1\insc@unt\dimen
     \ch@ck2\insc@unt\skip
     \ch@ck4\insc@unt\box
     \allocationnumber=\insc@unt
     \global\chardef#1=\allocationnumber
     \wlog{\string#1=\string\insert\the\allocationnumber}}
\def\newif#1{\count@\escapechar \escapechar\m@ne
     \expandafter\expandafter\expandafter
          \xdef\@if#1{true}{\let\noexpand#1=\noexpand\iftrue}%
     \expandafter\expandafter\expandafter
          \xdef\@if#1{false}{\let\noexpand#1=\noexpand\iffalse}%
     \global\@if#1{false}\escapechar=\count@}


\newlinechar=`\^^J
\overfullrule=0pt




\let\itfam=\undefined

\let\bffam=\undefined

\count18=3


\chardef\sharps="19


\mathchardef\alpha="710B \mathchardef\beta="710C
\mathchardef\gamma="710D \mathchardef\delta="710E
\mathchardef\epsilon="710F \mathchardef\zeta="7110
\mathchardef\eta="7111 \mathchardef\theta="7112 \mathchardef\iota="7113
\mathchardef\kappa="7114 \mathchardef\lambda="7115
\mathchardef\mu="7116 \mathchardef\nu="7117 \mathchardef\xi="7118
\mathchardef\pi="7119 \mathchardef\rho="711A \mathchardef\sigma="711B
\mathchardef\tau="711C \mathchardef\upsilon="711D
\mathchardef\phi="711E \mathchardef\chi="711F \mathchardef\psi="7120
\mathchardef\omega="7121 \mathchardef\varepsilon="7122
\mathchardef\vartheta="7123 \mathchardef\varpi="7124
\mathchardef\varrho="7125 \mathchardef\varsigma="7126
\mathchardef\varphi="7127 \mathchardef\imath="717B
\mathchardef\jmath="717C \mathchardef\ell="7160 \mathchardef\wp="717D
\mathchardef\partial="7140 \mathchardef\flat="715B
\mathchardef\natural="715C \mathchardef\sharp="715D



\def\angle{{\vbox{\ialign{$\m@th\scriptstyle##$\crcr
     \not\mathrel{\mkern14mu}\crcr
     \noalign{\nointerlineskip}
     \mkern2.5mu\leaders\hrule height.34\rp@\hfill\mkern2.5mu\crcr}}}}
\def\vdots{\vbox{\baselineskip4\rp@ \lineskiplimit\z@
     \kern6\rp@\hbox{.}\hbox{.}\hbox{.}}}
\def\ddots{\mathinner{\mkern1mu\raise7\rp@\vbox{\kern7\rp@\hbox{.}}\mkern2mu
     \raise4\rp@\hbox{.}\mkern2mu\raise\rp@\hbox{.}\mkern1mu}}
\def\overrightarrow#1{\vbox{\ialign{##\crcr
     \rightarrowfill\crcr
     \noalign{\kern-\rp@\nointerlineskip}
     $\hfil\displaystyle{#1}\hfil$\crcr}}}
\def\overleftarrow#1{\vbox{\ialign{##\crcr
     \leftarrowfill\crcr
     \noalign{\kern-\rp@\nointerlineskip}
     $\hfil\displaystyle{#1}\hfil$\crcr}}}
\def\overbrace#1{\mathop{\vbox{\ialign{##\crcr
     \noalign{\kern3\rp@}
     \downbracefill\crcr
     \noalign{\kern3\rp@\nointerlineskip}
     $\hfil\displaystyle{#1}\hfil$\crcr}}}\limits}
\def\underbrace#1{\mathop{\vtop{\ialign{##\crcr
     $\hfil\displaystyle{#1}\hfil$\crcr
     \noalign{\kern3\rp@\nointerlineskip}
     \upbracefill\crcr
     \noalign{\kern3\rp@}}}}\limits}
\def\big#1{{\hbox{$\left#1\vbox to8.5\rp@ {}\right.\n@space$}}}
\def\Big#1{{\hbox{$\left#1\vbox to11.5\rp@ {}\right.\n@space$}}}
\def\bigg#1{{\hbox{$\left#1\vbox to14.5\rp@ {}\right.\n@space$}}}
\def\Bigg#1{{\hbox{$\left#1\vbox to17.5\rp@ {}\right.\n@space$}}}
\def\@vereq#1#2{\lower.5\rp@\vbox{\baselineskip\z@skip\lineskip-.5\rp@
     \ialign{$\m@th#1\hfil##\hfil$\crcr#2\crcr=\crcr}}}
\def\rlh@#1{\vcenter{\hbox{\ooalign{\raise2\rp@
     \hbox{$#1\rightharpoonup$}\crcr
     $#1\leftharpoondown$}}}}
\def\bordermatrix#1{\begingroup\m@th
     \setbox\z@\vbox{%
          \def\cr{\crcr\noalign{\kern2\rp@\global\let\cr\endline}}%
          \ialign{$##$\hfil\kern2\rp@\kern\p@renwd
               &\thinspace\hfil$##$\hfil&&\quad\hfil$##$\hfil\crcr
               \omit\strut\hfil\crcr
               \noalign{\kern-\baselineskip}%
               #1\crcr\omit\strut\cr}}%
     \setbox\tw@\vbox{\unvcopy\z@\global\setbox\@ne\lastbox}%
     \setbox\tw@\hbox{\unhbox\@ne\unskip\global\setbox\@ne\lastbox}%
     \setbox\tw@\hbox{$\kern\wd\@ne\kern-\p@renwd\left(\kern-\wd\@ne
          \global\setbox\@ne\vbox{\box\@ne\kern2\rp@}%
          \vcenter{\kern-\ht\@ne\unvbox\z@\kern-\baselineskip}%
          \,\right)$}%
     \null\;\vbox{\kern\ht\@ne\box\tw@}\endgroup}
\def\endinsert{\egroup
     \if@mid\dimen@\ht\z@
          \advance\dimen@\dp\z@
          \advance\dimen@12\rp@
          \advance\dimen@\pagetotal
          \ifdim\dimen@>\pagegoal\@midfalse\p@gefalse\fi
     \fi
     \if@mid\bigskip\box\z@
          \bigbreak
     \else\insert\topins{\penalty100 \splittopskip\z@skip
               \splitmaxdepth\maxdimen\floatingpenalty\z@
               \ifp@ge\dimen@\dp\z@
                    \vbox to\vsize{\unvbox\z@\kern-\dimen@}%
               \else\box\z@\nobreak\bigskip
               \fi}%
     \fi
     \endgroup}


\def\cases#1{\left\{\,\vcenter{\m@th
     \ialign{$##\hfil$&\quad##\hfil\crcr#1\crcr}}\right.}
\def\matrix#1{\null\,\vcenter{\m@th
     \ialign{\hfil$##$\hfil&&\quad\hfil$##$\hfil\crcr
          \mathstrut\crcr
          \noalign{\kern-\baselineskip}
          #1\crcr
          \mathstrut\crcr
          \noalign{\kern-\baselineskip}}}\,}


\newif\ifraggedbottom

\def\raggedbottom{\ifraggedbottom\else
     \advance\topskip by\z@ plus60pt \raggedbottomtrue\fi}%
\def\normalbottom{\ifraggedbottom
     \advance\topskip by\z@ plus-60pt \raggedbottomfalse\fi}

\message{hacks,}


\toksdef\toks@i=1 \toksdef\toks@ii=2


\def\TeX{T\kern-.1667em \lower.5ex \hbox{E}\kern-.125em X\null}
\def\jyTeX{{\leavevmode
     \raise.587ex \hbox{\it\j}\kern-.1em \lower.048ex \hbox{\it y}\kern-.12em
     \TeX}}

\let\then=\iftrue
\def\ifnoarg#1\then{\def\hack@{#1}\ifx\hack@\empty}
\def\ifundefined#1\then{%
     \expandafter\ifx\csname\expandafter\blank\string#1\endcsname\relax}
\def\useif#1\then{\csname#1\endcsname}
\def\usename#1{\csname#1\endcsname}
\def\useafter#1#2{\expandafter#1\csname#2\endcsname}

\long\def\loop#1\repeat{\def\@iterate{#1\expandafter\@iterate\fi}\@iterate
     \let\@iterate=\relax}

\let\TeXend=\end
\def\begin#1{\begingroup\def\@@blockname{#1}\usename{begin#1}}
\def\end#1{\usename{end#1}\def\hack@{#1}%
     \ifx\@@blockname\hack@
          \endgroup
     \else\err@badgroup\hack@\@@blockname
     \fi}
\def\@@blockname{}

\def\defaultoption[#1]#2{%
     \def\hack@{\ifx\hack@ii[\toks@={#2}\else\toks@={#2[#1]}\fi\the\toks@}%
     \futurelet\hack@ii\hack@}

\def\markup#1{\let\@@marksf=\empty
     \ifhmode\edef\@@marksf{\spacefactor=\the\spacefactor\relax}\/\fi
     ${}^{\hbox{\subscriptfonts#1}}$\@@marksf}


\newtoks\shortyear
\newtoks\militaryhour
\newtoks\standardhour
\newtoks\minute
\newtoks\amorpm

\def\settime{\count@=\time\divide\count@ by60
     \militaryhour=\expandafter{\number\count@}%
     {\multiply\count@ by-60 \advance\count@ by\time
          \xdef\hack@{\ifnum\count@<10 0\fi\number\count@}}%
     \minute=\expandafter{\hack@}%
     \ifnum\count@<12
          \amorpm={am}
     \else\amorpm={pm}
          \ifnum\count@>12 \advance\count@ by-12 \fi
     \fi
     \standardhour=\expandafter{\number\count@}%
     \def\hack@19##1##2{\shortyear={##1##2}}%
          \expandafter\hack@\the\year}

\def\monthword#1{%
     \ifcase#1
          $\bullet$\err@badcountervalue{monthword}%
          \or January\or February\or March\or April\or May\or June%
          \or July\or August\or September\or October\or November\or December%
     \else$\bullet$\err@badcountervalue{monthword}%
     \fi}

\def\monthabbr#1{%
     \ifcase#1
          $\bullet$\err@badcountervalue{monthabbr}%
          \or Jan\or Feb\or Mar\or Apr\or May\or Jun%
          \or Jul\or Aug\or Sep\or Oct\or Nov\or Dec%
     \else$\bullet$\err@badcountervalue{monthabbr}%
     \fi}

\def\militarytime{\the\militaryhour:\the\minute}
\def\standardtime{\the\standardhour:\the\minute}


\def\@setnumstyle#1#2{\expandafter\global\expandafter\expandafter
     \expandafter\let\expandafter\expandafter
     \csname @\expandafter\blank\string#1style\endcsname
     \csname#2\endcsname}
\def\numstyle#1{\usename{@\expandafter\blank\string#1style}#1}
\def\ifblank#1\then{\useafter\ifx{@\expandafter\blank\string#1}\blank}

\def\blank#1{}

\def\Roman#1{\expandafter\uppercase\expandafter{\romannumeral#1}}
\def\alphabetic#1{%
     \ifcase#1
          $\bullet$\err@badcountervalue{alphabetic}%
          \or a\or b\or c\or d\or e\or f\or g\or h\or i\or j\or k\or l\or m%
          \or n\or o\or p\or q\or r\or s\or t\or u\or v\or w\or x\or y\or z%
     \else$\bullet$\err@badcountervalue{alphabetic}%
     \fi}
\def\Alphabetic#1{\expandafter\uppercase\expandafter{\alphabetic{#1}}}
\def\symbols#1{%
     \ifcase#1
          $\bullet$\err@badcountervalue{symbols}%
          \or*\or\dag\or\ddag\or\S\or$\|$%
          \or**\or\dag\dag\or\ddag\ddag\or\S\S\or$\|\|$%
     \else$\bullet$\err@badcountervalue{symbols}%
     \fi}


\catcode`\^^?=13 \def^^?{\relax}

\def\trimleading#1\to#2{\edef#2{#1}%
     \expandafter\@trimleading\expandafter#2#2^^?^^?}
\def\@trimleading#1#2#3^^?{\ifx#2^^?\def#1{}\else\def#1{#2#3}\fi}

\def\trimtrailing#1\to#2{\edef#2{#1}%
     \expandafter\@trimtrailing\expandafter#2#2^^? ^^?\relax}
\def\@trimtrailing#1#2 ^^?#3{\ifx#3\relax\toks@={}%
     \else\def#1{#2}\toks@={\trimtrailing#1\to#1}\fi
     \the\toks@}

\def\trim#1\to#2{\trimleading#1\to#2\trimtrailing#2\to#2}

\catcode`\^^?=15


\long\def\additemL#1\to#2{\toks@={\^^\{#1}}\toks@ii=\expandafter{#2}%
     \xdef#2{\the\toks@\the\toks@ii}}

\long\def\additemR#1\to#2{\toks@={\^^\{#1}}\toks@ii=\expandafter{#2}%
     \xdef#2{\the\toks@ii\the\toks@}}

\def\getitemL#1\to#2{\expandafter\@getitemL#1\hack@#1#2}
\def\@getitemL\^^\#1#2\hack@#3#4{\def#4{#1}\def#3{#2}}

\message{font macros,}


\newdimen\rp@
\newcount\@@sizeindex \@@sizeindex=0
\newcount\@@factori
\newcount\@@factorii
\newcount\@@factoriii
\newcount\@@factoriv

\countdef\maxfam=18
\newfam\itfam
\newfam\bffam
\newfam\bfsfam
\newfam\bmitfam

\def\@mathfontinit{\count@=4
     \loop\textfont\count@=\nullfont
          \scriptfont\count@=\nullfont
          \scriptscriptfont\count@=\nullfont
          \ifnum\count@<\maxfam\advance\count@ by\@ne
     \repeat}

\def\@fontstyleinit{%
     \def\it{\err@fontnotavailable\it}%
     \def\bf{\err@fontnotavailable\bf}%
     \def\bfs{\err@bfstobf}%
     \def\bmit{\err@fontnotavailable\bmit}%
     \def\sc{\err@fontnotavailable\sc}%
     \def\sl{\err@sltoit}%
     \def\ss{\err@fontnotavailable\ss}%
     \def\tt{\err@fontnotavailable\tt}}

\def\@parameterinit#1{\rm\rp@=.1em \@getscaling{#1}%
     \let\^^\=\@doscaling\scalingskipslist
     \setbox\strutbox=\hbox{\vrule
          height.708\baselineskip depth.292\baselineskip width\z@}}

\def\@getfactor#1#2#3#4{\@@factori=#1 \@@factorii=#2
     \@@factoriii=#3 \@@factoriv=#4}

\def\@getscaling#1{\count@=#1 \advance\count@ by-\@@sizeindex\@@sizeindex=#1
     \ifnum\count@<0
          \let\@mulordiv=\divide
          \let\@divormul=\multiply
          \multiply\count@ by\m@ne
     \else\let\@mulordiv=\multiply
          \let\@divormul=\divide
     \fi
     \edef\@@scratcha{\ifcase\count@                {1}{1}{1}{1}\or
          {1}{7}{23}{3}\or     {2}{5}{3}{1}\or      {9}{89}{13}{1}\or
          {6}{25}{6}{1}\or     {8}{71}{14}{1}\or    {6}{25}{36}{5}\or
          {1}{7}{53}{4}\or     {12}{125}{108}{5}\or {3}{14}{53}{5}\or
          {6}{41}{17}{1}\or    {13}{31}{13}{2}\or   {9}{107}{71}{2}\or
          {11}{139}{124}{3}\or {1}{6}{43}{2}\or     {10}{107}{42}{1}\or
          {1}{5}{43}{2}\or     {5}{69}{65}{1}\or    {11}{97}{91}{2}\fi}%
     \expandafter\@getfactor\@@scratcha}

\def\@doscaling#1{\@mulordiv#1by\@@factori\@divormul#1by\@@factorii
     \@mulordiv#1by\@@factoriii\@divormul#1by\@@factoriv}


\newskip\headskip
\newskip\footskip

\def\typesize=#1pt{\count@=#1 \advance\count@ by-10
     \ifcase\count@
          \@setsizex\or\err@badtypesize\or
          \@setsizexii\or\err@badtypesize\or
          \@setsizexiv
     \else\err@badtypesize
     \fi}

\def\@setsizex{\getixpt
     \def\subsubscriptfonts{\vpt}%
          \def\subsubscriptsize{\vpt\@parameterinit{-8}}%
     \def\subscriptfonts{\viipt}\def\subscriptsize{\viipt\@parameterinit{-4}}%
     \def\footnotefonts{\viiipt}\def\footnotesize{\viiipt\@parameterinit{-2}}%
     \def\smallfonts{\ixpt}\def\smallsize{\ixpt\@parameterinit{-1}}%
     \def\normalfonts{\xpt}\def\normalsize{\xpt\@parameterinit{0}}%
     \def\bigfonts{\xiipt}\def\bigsize{\xiipt\@parameterinit{2}}%
     \def\Bigfonts{\xivpt}\def\Bigsize{\xivpt\@parameterinit{4}}%
     \def\biggfonts{\xviipt}\def\biggsize{\xviipt\@parameterinit{6}}%
     \def\Biggfonts{\xxipt}\def\Biggsize{\xxipt\@parameterinit{8}}%
     \def\tinyfonts{\vpt}\def\tinysize{\vpt\@parameterinit{-8}}%
     \def\HUGEFONTS{\xxvpt}\def\HUGESIZE{\xxvpt\@parameterinit{10}}%
     \normalsize\fixedskipslist}

\def\@setsizexii{\getxipt
     \def\subsubscriptfonts{\vipt}%
          \def\subsubscriptsize{\vipt\@parameterinit{-6}}%
     \def\subscriptfonts{\viiipt}%
          \def\subscriptsize{\viiipt\@parameterinit{-2}}%
     \def\footnotefonts{\xpt}\def\footnotesize{\xpt\@parameterinit{0}}%
     \def\smallfonts{\xipt}\def\smallsize{\xipt\@parameterinit{1}}%
     \def\normalfonts{\xiipt}\def\normalsize{\xiipt\@parameterinit{2}}%
     \def\bigfonts{\xivpt}\def\bigsize{\xivpt\@parameterinit{4}}%
     \def\Bigfonts{\xviipt}\def\Bigsize{\xviipt\@parameterinit{6}}%
     \def\biggfonts{\xxipt}\def\biggsize{\xxipt\@parameterinit{8}}%
     \def\Biggfonts{\xxvpt}\def\Biggsize{\xxvpt\@parameterinit{10}}%
     \def\tinyfonts{\vpt}\def\tinysize{\vpt\@parameterinit{-8}}%
     \def\HUGEFONTS{\xxvpt}\def\HUGESIZE{\xxvpt\@parameterinit{10}}%
     \normalsize\fixedskipslist}

\def\@setsizexiv{\getxiiipt
     \def\subsubscriptfonts{\viipt}%
          \def\subsubscriptsize{\viipt\@parameterinit{-4}}%
     \def\subscriptfonts{\xpt}\def\subscriptsize{\xpt\@parameterinit{0}}%
     \def\footnotefonts{\xiipt}\def\footnotesize{\xiipt\@parameterinit{2}}%
     \def\smallfonts{\xiiipt}\def\smallsize{\xiiipt\@parameterinit{3}}%
     \def\normalfonts{\xivpt}\def\normalsize{\xivpt\@parameterinit{4}}%
     \def\bigfonts{\xviipt}\def\bigsize{\xviipt\@parameterinit{6}}%
     \def\Bigfonts{\xxipt}\def\Bigsize{\xxipt\@parameterinit{8}}%
     \def\biggfonts{\xxvpt}\def\biggsize{\xxvpt\@parameterinit{10}}%
     \def\Biggfonts{\err@sizetoolarge\Biggfonts\HUGEFONTS}%
          \def\Biggsize{\err@sizetoolarge\Biggsize\HUGESIZE}%
     \def\tinyfonts{\vpt}\def\tinysize{\vpt\@parameterinit{-8}}%
     \def\HUGEFONTS{\xxvpt}\def\HUGESIZE{\xxvpt\@parameterinit{10}}%
     \normalsize\fixedskipslist}

\def\subsubscriptfonts{\vpt} \def\subsubscriptsize{\vpt\@parameterinit{-8}}
\def\subscriptfonts{\viipt}  \def\subscriptsize{\viipt\@parameterinit{-4}}
\def\footnotefonts{\viiipt}  \def\footnotesize{\viiipt\@parameterinit{-2}}
\def\smallfonts{\err@sizenotavailable\smallfonts}
                             \def\smallsize{\ixpt\@parameterinit{-1}}
\def\normalfonts{\xpt}       \def\normalsize{\xpt\@parameterinit{0}}
\def\bigfonts{\xiipt}        \def\bigsize{\xiipt\@parameterinit{2}}
\def\Bigfonts{\xivpt}        \def\Bigsize{\xivpt\@parameterinit{4}}
\def\biggfonts{\xviipt}      \def\biggsize{\xviipt\@parameterinit{6}}
\def\Biggfonts{\xxipt}       \def\Biggsize{\xxipt\@parameterinit{8}}
\def\tinyfonts{\vpt}         \def\tinysize{\vpt\@parameterinit{-8}}
\def\HUGEFONTS{\xxvpt}       \def\HUGESIZE{\xxvpt\@parameterinit{10}}

\message{document layout,}


\newtoks\everyoutput \everyoutput={}
\newdimen\depthofpage
\newcount\pagenum \pagenum=0

\newdimen\oddtopmargin  \newdimen\eventopmargin
\newdimen\oddleftmargin \newdimen\evenleftmargin
\newtoks\oddhead        \newtoks\evenhead
\newtoks\oddfoot        \newtoks\evenfoot

\def\topmargin{\afterassignment\@seteventop\oddtopmargin}
\def\leftmargin{\afterassignment\@setevenleft\oddleftmargin}
\def\head{\afterassignment\@setevenhead\oddhead}
\def\foot{\afterassignment\@setevenfoot\oddfoot}

\def\@seteventop{\eventopmargin=\oddtopmargin}
\def\@setevenleft{\evenleftmargin=\oddleftmargin}
\def\@setevenhead{\evenhead=\oddhead}
\def\@setevenfoot{\evenfoot=\oddfoot}

\def\pagenumstyle#1{\@setnumstyle\pagenum{#1}}

\newif\ifdraft
\def\draft{\drafttrue\leftmargin=.5in \overfullrule=5pt }

\def\outputstyle#1{\global\expandafter\let\expandafter
          \@outputstyle\csname#1output\endcsname
     \usename{#1setup}}

\output={\@outputstyle}

\def\normaloutput{\the\everyoutput
     \global\advance\pagenum by\@ne
     \ifodd\pagenum
          \voffset=\oddtopmargin \hoffset=\oddleftmargin
     \else\voffset=\eventopmargin \hoffset=\evenleftmargin
     \fi
     \advance\voffset by-1in  \advance\hoffset by-1in
     \count0=\pagenum
     \expandafter\shipout\pagebox
     \ifnum\outputpenalty>-\@MM\else\dosupereject\fi}

\newdimen\fullhsize
\newbox\leftpage
\newcount\leftpagenum
\newcount\outputpagenum \outputpagenum=0
\let\leftorright=L

\def\twoupoutput{\the\everyoutput
     \global\advance\pagenum by\@ne
     \if L\leftorright
          \global\setbox\leftpage=\leftline{\pagebox}%
          \global\leftpagenum=\pagenum
          \global\let\leftorright=R%
     \else\global\advance\outputpagenum by\@ne
          \ifodd\outputpagenum
               \voffset=\oddtopmargin \hoffset=\oddleftmargin
          \else\voffset=\eventopmargin \hoffset=\evenleftmargin
          \fi
          \advance\voffset by-1in  \advance\hoffset by-1in
          \count0=\leftpagenum \count1=\pagenum
          \shipout\vbox{\hbox to\fullhsize
               {\box\leftpage\hfil\leftline{\pagebox}}}%
          \global\let\leftorright=L%
     \fi
     \ifnum\outputpenalty>-\@MM
     \else\dosupereject
          \if R\leftorright
               \globaldefs=\@ne\head={\hfil}\foot={\hfil}\globaldefs=\z@
               \null\newpage
          \fi
     \fi}

\def\pagebox{\vbox{\makeheadline\pagebody\makefootline}}

\def\makeheadline{%
     \vbox to\z@{\baselinestretch=\@m
          \vskip\topskip\vskip-.708\baselineskip\vskip-\headskip
          \line{\vbox to\ht\strutbox{}%
               \ifodd\pagenum\the\oddhead\else\the\evenhead\fi}%
          \vss}%
     \nointerlineskip}

\def\pagebody{\vbox to\vsize{%
     \boxmaxdepth\maxdepth
     \ifvoid\topins\else\unvbox\topins\fi
     \depthofpage=\dp255
     \unvbox255
     \ifraggedbottom\kern-\depthofpage\vfil\fi
     \ifvoid\footins
     \else\vskip\skip\footins
          \footnoterule
          \unvbox\footins
          \vskip-\footnoteskip
     \fi}}

\def\makefootline{\baselineskip=\footskip
     \line{\ifodd\pagenum\the\oddfoot\else\the\evenfoot\fi}}


\newskip\abovechapterskip
\newskip\belowchapterskip
\newskip\abovesectionskip
\newskip\belowsectionskip
\newskip\abovesubsectionskip
\newskip\belowsubsectionskip

\def\chapterstyle#1{\global\expandafter\let\expandafter\@chapterstyle
     \csname#1text\endcsname}
\def\sectionstyle#1{\global\expandafter\let\expandafter\@sectionstyle
     \csname#1text\endcsname}
\def\subsectionstyle#1{\global\expandafter\let\expandafter\@subsectionstyle
     \csname#1text\endcsname}

\def\chapter#1{%
     \ifdim\lastskip=17sp \else\chapterbreak\vskip\abovechapterskip\fi
     \@chapterstyle{\ifblank\chapternumstyle\then
          \else\newchapternum=\next\chapternumformat\ \fi#1}%
     \nobreak\vskip\belowchapterskip\vskip17sp }

\def\section#1{%
     \ifdim\lastskip=17sp \else\sectionbreak\vskip\abovesectionskip\fi
     \@sectionstyle{\ifblank\sectionnumstyle\then
          \else\newsectionnum=\next\sectionnumformat\ \fi#1}%
     \nobreak\vskip\belowsectionskip\vskip17sp }

\def\subsection#1{%
     \ifdim\lastskip=17sp \else\subsectionbreak\vskip\abovesubsectionskip\fi
     \@subsectionstyle{\ifblank\subsectionnumstyle\then
          \else\newsubsectionnum=\next\subsectionnumformat\ \fi#1}%
     \nobreak\vskip\belowsubsectionskip\vskip17sp }


\let\TeXunderline=\underline
\let\TeXoverline=\overline
\def\underline#1{\relax\ifmmode\TeXunderline{#1}\else
     $\TeXunderline{\hbox{#1}}$\fi}
\def\overline#1{\relax\ifmmode\TeXoverline{#1}\else
     $\TeXoverline{\hbox{#1}}$\fi}

\def\baselinestretch{\afterassignment\@baselinestretch\count@}
\def\@baselinestretch{\baselineskip=\normalbaselineskip
     \divide\baselineskip by\@m\baselineskip=\count@\baselineskip
     \setbox\strutbox=\hbox{\vrule
          height.708\baselineskip depth.292\baselineskip width\z@}%
     \bigskipamount=\the\baselineskip
          plus.25\baselineskip minus.25\baselineskip
     \medskipamount=.5\baselineskip
          plus.125\baselineskip minus.125\baselineskip
     \smallskipamount=.25\baselineskip
          plus.0625\baselineskip minus.0625\baselineskip}

\def\\{\ifhmode\ifnum\lastpenalty=-\@M\else\hfil\penalty-\@M\fi\fi
     \ignorespaces}
\def\newpage{\vfil\break}

\def\lefttext#1{\par{\@text\leftskip=\z@\rightskip=\centering
     \noindent#1\par}}
\def\righttext#1{\par{\@text\leftskip=\centering\rightskip=\z@
     \noindent#1\par}}
\def\centertext#1{\par{\@text\leftskip=\centering\rightskip=\centering
     \noindent#1\par}}
\def\@text{\parindent=\z@ \parfillskip=\z@ \everypar={}%
     \spaceskip=.3333em \xspaceskip=.5em
     \def\\{\ifhmode\ifnum\lastpenalty=-\@M\else\penalty-\@M\fi\fi
          \ignorespaces}}

\def\beginleft{\par\@text\leftskip=\z@ \rightskip=\centering}
     
\def\beginright{\par\@text\leftskip=\centering\rightskip=\z@ }
     
\def\begincenter{\par\@text\leftskip=\centering\rightskip=\centering}

\def\beginnarrow{\defaultoption[\parindent]\@beginnarrow}
\def\@beginnarrow[#1]{\par\advance\leftskip by#1\advance\rightskip by#1}

\begingroup
\catcode`\[=1 \catcode`\{=11 \gdef\beginignore[\endgroup\bgroup
     \catcode`\e=0 \catcode`\\=12 \catcode`\{=11 \catcode`\f=12 \let\or=\relax
     \let\nd{ignor=\fi \let\}=\egroup
     \iffalse}
\endgroup

\long\def\marginnote#1{\leavevmode
     \edef\@marginsf{\spacefactor=\the\spacefactor\relax}%
     \ifdraft\strut\vadjust{%
          \hbox to\z@{\hskip\hsize\hskip.1in
               \vbox to\z@{\vskip-\dp\strutbox
                    \marginnoteformat
                    \vskip-\ht\strutbox
                    \noindent\strut#1\par
                    \vss}%
               \hss}}%
     \fi
     \@marginsf}


\newtoks\everybye \everybye={\par\vfil}
\outer\def\bye{\the\everybye
     \footnotecheck
     \prelabelcheck
     \streamcheck
     \supereject
     \TeXend}

\message{footnotes,}

\newcount\footnotenum \footnotenum=0
\newskip\footnoteskip
\let\@footnotelist=\empty

\def\footnotenumstyle#1{\@setnumstyle\footnotenum{#1}%
     \useafter\ifx{@footnotenumstyle}\symbols
          \global\let\@footup=\empty
     \else\global\let\@footup=\markup
     \fi}

\def\footnote{\footnotecheck\defaultoption[]\@footnote}
\def\@footnote[#1]{\@footnotemark[#1]\@footnotetext}

\def\footnotemark{\defaultoption[]\@footnotemark}
\def\@footnotemark[#1]{\let\@footsf=\empty
     \ifhmode\edef\@footsf{\spacefactor=\the\spacefactor\relax}\/\fi
     \ifnoarg#1\then
          \global\advance\footnotenum by\@ne
          \@footup{\footnotenumformat}%
          \edef\@@foota{\footnotenum=\the\footnotenum\relax}%
          \expandafter\additemR\expandafter\@footup\expandafter
               {\@@foota\footnotenumformat}\to\@footnotelist
          \global\let\@footnotelist=\@footnotelist
     \else\markup{#1}%
          \additemR\markup{#1}\to\@footnotelist
          \global\let\@footnotelist=\@footnotelist
     \fi
     \@footsf}

\def\footnotetext{%
     \ifx\@footnotelist\empty\err@extrafootnotetext\else\@footnotetext\fi}
\def\@footnotetext{%
     \getitemL\@footnotelist\to\@@foota
     \global\let\@footnotelist=\@footnotelist
     \insert\footins\bgroup
     \footnoteformat
     \splittopskip=\ht\strutbox\splitmaxdepth=\dp\strutbox
     \interlinepenalty=\interfootnotelinepenalty\floatingpenalty=\@MM
     \noindent\llap{\@@foota}\strut
     \bgroup\aftergroup\@footnoteend
     \let\@@scratcha=}
\def\@footnoteend{\strut\par\vskip\footnoteskip\egroup}

\def\footnoterule{\normalfonts
     \kern-.3em \hrule width2in height.04em \kern .26em }

\def\footnotecheck{%
     \ifx\@footnotelist\empty
     \else\err@extrafootnotemark
          \global\let\@footnotelist=\empty
     \fi}

\message{labels,}

\let\@@labeldef=\xdef
\newif\if@labelfile
\newwrite\@labelfile
\let\@prelabellist=\empty

\def\label#1#2{\trim#1\to\@@labarg\edef\@@labtext{#2}%
     \edef\@@labname{lab@\@@labarg}%
     \useafter\ifundefined\@@labname\then\else\@yeslab\fi
     \useafter\@@labeldef\@@labname{#2}%
     \ifstreaming
          \expandafter\toks@\expandafter\expandafter\expandafter
               {\csname\@@labname\endcsname}%
          \immediate\write\streamout{\noexpand\label{\@@labarg}{\the\toks@}}%
     \fi}
\def\@yeslab{%
     \useafter\ifundefined{if\@@labname}\then
          \err@labelredef\@@labarg
     \else\useif{if\@@labname}\then
               \err@labelredef\@@labarg
          \else\global\usename{\@@labname true}%
               \useafter\ifundefined{pre\@@labname}\then
               \else\useafter\ifx{pre\@@labname}\@@labtext
                    \else\err@badlabelmatch\@@labarg
                    \fi
               \fi
               \if@labelfile
               \else\global\@labelfiletrue
                    \immediate\write\sixt@@n{--> Creating file \jobname.lab}%
                    \immediate\openout\@labelfile=\jobname.lab
               \fi
               \immediate\write\@labelfile
                    {\noexpand\prelabel{\@@labarg}{\@@labtext}}%
          \fi
     \fi}

\def\putlab#1{\trim#1\to\@@labarg\edef\@@labname{lab@\@@labarg}%
     \useafter\ifundefined\@@labname\then\@nolab\else\usename\@@labname\fi}
\def\@nolab{%
     \useafter\ifundefined{pre\@@labname}\then
          \undefinedlabelformat
          \err@needlabel\@@labarg
          \useafter\xdef\@@labname{\undefinedlabelformat}%
     \else\usename{pre\@@labname}%
          \useafter\xdef\@@labname{\usename{pre\@@labname}}%
     \fi
     \useafter\newif{if\@@labname}%
     \expandafter\additemR\@@labarg\to\@prelabellist}

\def\prelabel#1{\useafter\gdef{prelab@#1}}

\def\ifundefinedlabel#1\then{%
     \expandafter\ifx\csname lab@#1\endcsname\relax}
\def\useiflab#1\then{\csname iflab@#1\endcsname}

\def\prelabelcheck{{%
     \def\^^\##1{\useiflab{##1}\then\else\err@undefinedlabel{##1}\fi}%
     \@prelabellist}}

\message{equation numbering,}

\newcount\chapternum
\newcount\sectionnum
\newcount\subsectionnum
\newcount\equationnum
\newcount\subequationnum
\newcount\figurenum
\newcount\subfigurenum
\newcount\tablenum
\newcount\subtablenum

\newif\if@subeqncount
\newif\if@subfigcount
\newif\if@subtblcount

\def\newchapternum{\newsectionnum=\z@\@resetnum\chapternum}
\def\newsectionnum{\newsubsectionnum=\z@\@resetnum\sectionnum}
\def\newsubsectionnum{\newequationnum=\z@\newfigurenum=\z@\newtablenum=\z@
     \@resetnum\subsectionnum}
\def\newequationnum{\newsubequationnum=\z@\@resetnum\equationnum}
\def\newsubequationnum{\@resetnum\subequationnum}
\def\newfigurenum{\newsubfigurenum=\z@\@resetnum\figurenum}
\def\newsubfigurenum{\@resetnum\subfigurenum}
\def\newtablenum{\newsubtablenum=\z@\@resetnum\tablenum}
\def\newsubtablenum{\@resetnum\subtablenum}

\def\@resetnum#1{\global\advance#1by1 \edef\next{\the#1\relax}\global#1}

\newchapternum=0

\def\chapternumstyle#1{\@setnumstyle\chapternum{#1}}
\def\sectionnumstyle#1{\@setnumstyle\sectionnum{#1}}
\def\subsectionnumstyle#1{\@setnumstyle\subsectionnum{#1}}
\def\equationnumstyle#1{\@setnumstyle\equationnum{#1}}
\def\subequationnumstyle#1{\@setnumstyle\subequationnum{#1}%
     \ifblank\subequationnumstyle\then\global\@subeqncountfalse\fi
     \ignorespaces}
\def\figurenumstyle#1{\@setnumstyle\figurenum{#1}}
\def\subfigurenumstyle#1{\@setnumstyle\subfigurenum{#1}%
     \ifblank\subfigurenumstyle\then\global\@subfigcountfalse\fi
     \ignorespaces}
\def\tablenumstyle#1{\@setnumstyle\tablenum{#1}}
\def\subtablenumstyle#1{\@setnumstyle\subtablenum{#1}%
     \ifblank\subtablenumstyle\then\global\@subtblcountfalse\fi
     \ignorespaces}

\def\eqnlabel#1{%
     \if@subeqncount
          \newsubequationnum=\next
     \else\newequationnum=\next
          \ifblank\subequationnumstyle\then
          \else\global\@subeqncounttrue
               \newsubequationnum=\@ne
          \fi
     \fi
     \label{#1}{\puteqnformat}(\puteqn{#1})%
     \ifdraft\rlap{\hskip.1in{\tt#1}}\fi}

\let\puteqn=\putlab

\def\equation#1#2{\useafter\gdef{eqn@#1}{#2\eqno\eqnlabel{#1}}}
\def\Equation#1{\useafter\gdef{eqn@#1}}

\def\putequation#1{\useafter\ifundefined{eqn@#1}\then
     \err@undefinedeqn{#1}\else\usename{eqn@#1}\fi}

\def\eqnseriesstyle#1{\gdef\@eqnseriesstyle{#1}}
\def\begineqnseries{\subequationnumstyle{\@eqnseriesstyle}%
     \defaultoption[]\@begineqnseries}
\def\@begineqnseries[#1]{\edef\@@eqnname{#1}}
\def\endeqnseries{\subequationnumstyle{blank}%
     \expandafter\ifnoarg\@@eqnname\then
     \else\label\@@eqnname{\puteqnformat}%
     \fi
     \aftergroup\ignorespaces}

\def\figlabel#1{%
     \if@subfigcount
          \newsubfigurenum=\next
     \else\newfigurenum=\next
          \ifblank\subfigurenumstyle\then
          \else\global\@subfigcounttrue
               \newsubfigurenum=\@ne
          \fi
     \fi
     \label{#1}{\putfigformat}\putfig{#1}%
     {\def\marginnoteformat{\tt}\marginnote{#1}}}

\let\putfig=\putlab

\def\figseriesstyle#1{\gdef\@figseriesstyle{#1}}
\def\beginfigseries{\subfigurenumstyle{\@figseriesstyle}%
     \defaultoption[]\@beginfigseries}
\def\@beginfigseries[#1]{\edef\@@figname{#1}}
\def\endfigseries{\subfigurenumstyle{blank}%
     \expandafter\ifnoarg\@@figname\then
     \else\label\@@figname{\putfigformat}%
     \fi
     \aftergroup\ignorespaces}

\def\tbllabel#1{%
     \if@subtblcount
          \newsubtablenum=\next
     \else\newtablenum=\next
          \ifblank\subtablenumstyle\then
          \else\global\@subtblcounttrue
               \newsubtablenum=\@ne
          \fi
     \fi
     \label{#1}{\puttblformat}\puttbl{#1}%
     {\def\marginnoteformat{\tt}\marginnote{#1}}}

\let\puttbl=\putlab

\def\tblseriesstyle#1{\gdef\@tblseriesstyle{#1}}
\def\begintblseries{\subtablenumstyle{\@tblseriesstyle}%
     \defaultoption[]\@begintblseries}
\def\@begintblseries[#1]{\edef\@@tblname{#1}}
\def\endtblseries{\subtablenumstyle{blank}%
     \expandafter\ifnoarg\@@tblname\then
     \else\label\@@tblname{\puttblformat}%
     \fi
     \aftergroup\ignorespaces}

\message{reference numbering,}

\newcount\referencenum \referencenum=0
\newcount\@@prerefcount \@@prerefcount=0
\newcount\@@thisref
\newcount\@@lastref
\newcount\@@loopref
\newcount\@@refseq
\newdimen\refnumindent
\let\@undefreflist=\empty

\def\referencenumstyle#1{\@setnumstyle\referencenum{#1}}

\def\referencestyle#1{\usename{@ref#1}}

\def\@refsequential{%
     \gdef\@refpredef##1{\global\advance\referencenum by\@ne
          \let\^^\=0\label{##1}{\^^\{\the\referencenum}}%
          \useafter\gdef{ref@\the\referencenum}{{##1}{\undefinedlabelformat}}}%
     \gdef\@reference##1##2{%
          \ifundefinedlabel##1\then
          \else\def\^^\####1{\global\@@thisref=####1\relax}\putlab{##1}%
               \useafter\gdef{ref@\the\@@thisref}{{##1}{##2}}%
          \fi}%
     \gdef\endputreferences{%
          \loop\ifnum\@@loopref<\referencenum
                    \advance\@@loopref by\@ne
                    \expandafter\expandafter\expandafter\@printreference
                         \csname ref@\the\@@loopref\endcsname
          \repeat
          \par}}

\def\@refpreordered{%
     \gdef\@refpredef##1{\global\advance\referencenum by\@ne
          \additemR##1\to\@undefreflist}%
     \gdef\@reference##1##2{%
          \ifundefinedlabel##1\then
          \else\global\advance\@@loopref by\@ne
               {\let\^^\=0\label{##1}{\^^\{\the\@@loopref}}}%
               \@printreference{##1}{##2}%
          \fi}
     \gdef\endputreferences{%
          \def\^^\####1{\useiflab{####1}\then
               \else\reference{####1}{\undefinedlabelformat}\fi}%
          \@undefreflist
          \par}}

\def\beginprereferences{\par
     \def\reference##1##2{\global\advance\referencenum by1\@ne
          \let\^^\=0\label{##1}{\^^\{\the\referencenum}}%
          \useafter\gdef{ref@\the\referencenum}{{##1}{##2}}}}
\def\endprereferences{\global\@@prerefcount=\the\referencenum\par}

\def\beginputreferences{\par
     \refnumindent=\z@\@@loopref=\z@
     \loop\ifnum\@@loopref<\referencenum
               \advance\@@loopref by\@ne
               \setbox\z@=\hbox{\referencenum=\@@loopref
                    \referencenumformat\enskip}%
               \ifdim\wd\z@>\refnumindent\refnumindent=\wd\z@\fi
     \repeat
     \putreferenceformat
     \@@loopref=\z@
     \loop\ifnum\@@loopref<\@@prerefcount
               \advance\@@loopref by\@ne
               \expandafter\expandafter\expandafter\@printreference
                    \csname ref@\the\@@loopref\endcsname
     \repeat
     \let\reference=\@reference}

\def\@printreference#1#2{\ifx#2\undefinedlabelformat\err@undefinedref{#1}\fi
     \noindent\ifdraft\rlap{\hskip\hsize\hskip.1in \tt#1}\fi
     \llap{\referencenum=\@@loopref\referencenumformat\enskip}#2\par}

\def\reference#1#2{{\par\refnumindent=\z@\putreferenceformat\noindent#2\par}}

\def\putref#1{\trim#1\to\@@refarg
     \expandafter\ifnoarg\@@refarg\then
          \toks@={\relax}%
     \else\@@lastref=-\@m\def\@@refsep{}\def\@more{\@nextref}%
          \toks@={\@nextref#1,,}%
     \fi\the\toks@}
\def\@nextref#1,{\trim#1\to\@@refarg
     \expandafter\ifnoarg\@@refarg\then
          \let\@more=\relax
     \else\ifundefinedlabel\@@refarg\then
               \expandafter\@refpredef\expandafter{\@@refarg}%
          \fi
          \def\^^\##1{\global\@@thisref=##1\relax}%
          \global\@@thisref=\m@ne
          \setbox\z@=\hbox{\putlab\@@refarg}%
     \fi
     \advance\@@lastref by\@ne
     \ifnum\@@lastref=\@@thisref\advance\@@refseq by\@ne\else\@@refseq=\@ne\fi
     \ifnum\@@lastref<\z@
     \else\ifnum\@@refseq<\thr@@
               \@@refsep\def\@@refsep{,}%
               \ifnum\@@lastref>\z@
                    \advance\@@lastref by\m@ne
                    {\referencenum=\@@lastref\putrefformat}%
               \else\undefinedlabelformat
               \fi
          \else\def\@@refsep{--}%
          \fi
     \fi
     \@@lastref=\@@thisref
     \@more}

\message{streaming,}

\newif\ifstreaming

\def\streamto{\defaultoption[\jobname]\@streamto}
\def\@streamto[#1]{\global\streamingtrue
     \immediate\write\sixt@@n{--> Streaming to #1.str}%
     \newwrite\streamout\immediate\openout\streamout=#1.str }

\def\streamfrom{\defaultoption[\jobname]\@streamfrom}
\def\@streamfrom[#1]{\newread\streamin\openin\streamin=#1.str
     \ifeof\streamin
          \expandafter\err@nostream\expandafter{#1.str}%
     \else\immediate\write\sixt@@n{--> Streaming from #1.str}%
          \let\@@labeldef=\gdef
          \ifstreaming
               \edef\@elc{\endlinechar=\the\endlinechar}%
               \endlinechar=\m@ne
               \loop\read\streamin to\@@scratcha
                    \ifeof\streamin
                         \streamingfalse
                    \else\toks@=\expandafter{\@@scratcha}%
                         \immediate\write\streamout{\the\toks@}%
                    \fi
                    \ifstreaming
               \repeat
               \@elc
               \input #1.str
               \streamingtrue
          \else\input #1.str
          \fi
          \let\@@labeldef=\xdef
     \fi}

\def\streamcheck{\ifstreaming
     \immediate\write\streamout{\pagenum=\the\pagenum}%
     \immediate\write\streamout{\footnotenum=\the\footnotenum}%
     \immediate\write\streamout{\referencenum=\the\referencenum}%
     \immediate\write\streamout{\chapternum=\the\chapternum}%
     \immediate\write\streamout{\sectionnum=\the\sectionnum}%
     \immediate\write\streamout{\subsectionnum=\the\subsectionnum}%
     \immediate\write\streamout{\equationnum=\the\equationnum}%
     \immediate\write\streamout{\subequationnum=\the\subequationnum}%
     \immediate\write\streamout{\figurenum=\the\figurenum}%
     \immediate\write\streamout{\subfigurenum=\the\subfigurenum}%
     \immediate\write\streamout{\tablenum=\the\tablenum}%
     \immediate\write\streamout{\subtablenum=\the\subtablenum}%
     \immediate\closeout\streamout
     \fi}


\def\err@badtypesize{%
     \errhelp={The limited availability of certain fonts requires^^J%
          that the base type size be 10pt, 12pt, or 14pt.^^J}%
     \errmessage{--> Illegal base type size}}

\def\err@badsizechange{\immediate\write\sixt@@n
     {--> Size change not allowed in math mode, ignored}}

\def\err@sizetoolarge#1{\immediate\write\sixt@@n
     {--> \noexpand#1 too big, substituting HUGE}}

\def\err@sizenotavailable#1{\immediate\write\sixt@@n
     {--> Size not available, \noexpand#1 ignored}}

\def\err@fontnotavailable#1{\immediate\write\sixt@@n
     {--> Font not available, \noexpand#1 ignored}}

\def\err@sltoit{\immediate\write\sixt@@n
     {--> Style \noexpand\sl not available, substituting \noexpand\it}%
     \it}

\def\err@bfstobf{\immediate\write\sixt@@n
     {--> Style \noexpand\bfs not available, substituting \noexpand\bf}%
     \bf}

\def\err@badgroup#1#2{%
     \errhelp={The block you have just tried to close was not the one^^J%
          most recently opened.^^J}%
     \errmessage{--> \noexpand\end{#1} doesn't match \noexpand\begin{#2}}}

\def\err@badcountervalue#1{\immediate\write\sixt@@n
     {--> Counter (#1) out of bounds}}

\def\err@extrafootnotemark{\immediate\write\sixt@@n
     {--> \noexpand\footnotemark command
          has no corresponding \noexpand\footnotetext}}

\def\err@extrafootnotetext{%
     \errhelp{You have given a \noexpand\footnotetext command without first
          specifying^^Ja \noexpand\footnotemark.^^J}%
     \errmessage{--> \noexpand\footnotetext command has no corresponding
          \noexpand\footnotemark}}

\def\err@labelredef#1{\immediate\write\sixt@@n
     {--> Label "#1" redefined}}

\def\err@badlabelmatch#1{\immediate\write\sixt@@n
     {--> Definition of label "#1" doesn't match value in \jobname.lab}}

\def\err@needlabel#1{\immediate\write\sixt@@n
     {--> Label "#1" cited before its definition}}

\def\err@undefinedlabel#1{\immediate\write\sixt@@n
     {--> Label "#1" cited but never defined}}

\def\err@undefinedeqn#1{\immediate\write\sixt@@n
     {--> Equation "#1" not defined}}

\def\err@undefinedref#1{\immediate\write\sixt@@n
     {--> Reference "#1" not defined}}

\def\err@nostream#1{%
     \errhelp={You have tried to input a stream file that doesn't exist.^^J}%
     \errmessage{--> Stream file #1 not found}}

\message{jyTeX initialization}

\everyjob{\immediate\write16{--> jyTeX version \fmtversion}%
     \edef\@@jobname{\jobname}%
     \edef\jobname{\@@jobname}%
     \settime
     \openin0=\jobname.lab
     \ifeof0
     \else\closein0
          \immediate\write16{--> Getting labels from file \jobname.lab}%
          \input\jobname.lab
     \fi}


\def\fixedskipslist{%
     \^^\{\topskip}%
     \^^\{\splittopskip}%
     \^^\{\maxdepth}%
     \^^\{\skip\topins}%
     \^^\{\skip\footins}%
     \^^\{\headskip}%
     \^^\{\footskip}}

\def\scalingskipslist{%
     \^^\{\p@renwd}%
     \^^\{\delimitershortfall}%
     \^^\{\nulldelimiterspace}%
     \^^\{\scriptspace}%
     \^^\{\jot}%
     \^^\{\normalbaselineskip}%
     \^^\{\normallineskip}%
     \^^\{\normallineskiplimit}%
     \^^\{\baselineskip}%
     \^^\{\lineskip}%
     \^^\{\lineskiplimit}%
     \^^\{\bigskipamount}%
     \^^\{\medskipamount}%
     \^^\{\smallskipamount}%
     \^^\{\parskip}%
     \^^\{\parindent}%
     \^^\{\abovedisplayskip}%
     \^^\{\belowdisplayskip}%
     \^^\{\abovedisplayshortskip}%
     \^^\{\belowdisplayshortskip}%
     \^^\{\abovechapterskip}%
     \^^\{\belowchapterskip}%
     \^^\{\abovesectionskip}%
     \^^\{\belowsectionskip}%
     \^^\{\abovesubsectionskip}%
     \^^\{\belowsubsectionskip}}


\def\twoupsetup{
     \topmargin=.75in
     \leftmargin=.5in
     \vsize=6.9in
     \hsize=4.75in
     \fullhsize=10in
     \let\draft=\relax}

\outputstyle{normal}                             

\def\marginnoteformat{\subscriptsize             
     \hsize=1in \baselinestretch=1000 \everypar={}%
     \tolerance=5000 \hbadness=5000 \parskip=0pt \parindent=0pt
     \leftskip=0pt \rightskip=0pt \raggedright}

\head={\ifdraft\normalfonts\it\hfil DRAFT\hfil   
     \llap{\number\day\ \monthword\month\ \militarytime}\else\hfil\fi}
\foot={\hfil\normalfonts\numstyle\pagenum\hfil}  

\normalbaselineskip=12pt                         
\normallineskip=0pt                              
\normallineskiplimit=0pt                         
\normalbaselines                                 

\topskip=.85\baselineskip \splittopskip=\topskip \headskip=2\baselineskip
\footskip=\headskip

\pagenumstyle{arabic}                            

\parskip=0pt                                     
\parindent=20pt                                  

\baselinestretch=1000                            


\chapterstyle{left}                              
\chapternumstyle{blank}                          
\def\chapterbreak{\newpage}                      
\abovechapterskip=0pt                            
\belowchapterskip=1.5\baselineskip               
     plus.38\baselineskip minus.38\baselineskip
\def\chapternumformat{\numstyle\chapternum.}     

\sectionstyle{left}                              
\sectionnumstyle{blank}                          
\def\sectionbreak{\vskip0pt plus4\baselineskip\penalty-100
     \vskip0pt plus-4\baselineskip}              
\abovesectionskip=1.5\baselineskip               
     plus.38\baselineskip minus.38\baselineskip
\belowsectionskip=\the\baselineskip              
     plus.25\baselineskip minus.25\baselineskip
\def\sectionnumformat{
     \ifblank\chapternumstyle\then\else\numstyle\chapternum.\fi
     \numstyle\sectionnum.}

\subsectionstyle{left}                           
\subsectionnumstyle{blank}                       
\def\subsectionbreak{\vskip0pt plus4\baselineskip\penalty-100
     \vskip0pt plus-4\baselineskip}              
\abovesubsectionskip=\the\baselineskip           
     plus.25\baselineskip minus.25\baselineskip
\belowsubsectionskip=.75\baselineskip            
     plus.19\baselineskip minus.19\baselineskip
\def\subsectionnumformat{
     \ifblank\chapternumstyle\then\else\numstyle\chapternum.\fi
     \ifblank\sectionnumstyle\then\else\numstyle\sectionnum.\fi
     \numstyle\subsectionnum.}


\footnotenumstyle{symbols}                       
\footnoteskip=0pt                                
\def\footnotenumformat{\numstyle\footnotenum}    
\def\footnoteformat{\footnotesize                
     \everypar={}\parskip=0pt \parfillskip=0pt plus1fil
     \leftskip=1em \rightskip=0pt
     \spaceskip=0pt \xspaceskip=0pt
     \def\\{\ifhmode\ifnum\lastpenalty=-10000
          \else\hfil\penalty-10000 \fi\fi\ignorespaces}}


\def\undefinedlabelformat{$\bullet$}             


\equationnumstyle{arabic}                        
\subequationnumstyle{blank}                      
\figurenumstyle{arabic}                          
\subfigurenumstyle{blank}                        
\tablenumstyle{arabic}                           
\subtablenumstyle{blank}                         

\eqnseriesstyle{alphabetic}                      
\figseriesstyle{alphabetic}                      
\tblseriesstyle{alphabetic}                      

\def\puteqnformat{\hbox{
     \ifblank\chapternumstyle\then\else\numstyle\chapternum.\fi
     \ifblank\sectionnumstyle\then\else\numstyle\sectionnum.\fi
     \ifblank\subsectionnumstyle\then\else\numstyle\subsectionnum.\fi
     \numstyle\equationnum
     \numstyle\subequationnum}}
\def\putfigformat{\hbox{
     \ifblank\chapternumstyle\then\else\numstyle\chapternum.\fi
     \ifblank\sectionnumstyle\then\else\numstyle\sectionnum.\fi
     \ifblank\subsectionnumstyle\then\else\numstyle\subsectionnum.\fi
     \numstyle\figurenum
     \numstyle\subfigurenum}}
\def\puttblformat{\hbox{
     \ifblank\chapternumstyle\then\else\numstyle\chapternum.\fi
     \ifblank\sectionnumstyle\then\else\numstyle\sectionnum.\fi
     \ifblank\subsectionnumstyle\then\else\numstyle\subsectionnum.\fi
     \numstyle\tablenum
     \numstyle\subtablenum}}


\referencestyle{sequential}                      
\referencenumstyle{arabic}                       
\def\putrefformat{\numstyle\referencenum}        
\def\referencenumformat{\numstyle\referencenum.} 
\def\putreferenceformat{
     \everypar={\hangindent=1em \hangafter=1 }%
     \def\\{\hfil\break\null\hskip-1em \ignorespaces}%
     \leftskip=\refnumindent\parindent=0pt \interlinepenalty=1000 }


\normalsize


\def\fmtversion{2.6M (June 1992)}

\catcode`\@=12

\typesize=10pt \magnification=1200 \baselineskip17truept
\footnotenumstyle{arabic} \hsize=6truein\vsize=8.5truein
\sectionnumstyle{blank}
\chapternumstyle{blank}
\chapternum=1
\sectionnum=1
\pagenum=0

\def\begintitle{\pagenumstyle{blank}\parindent=0pt
\begin{narrow}[0.4in]}
\def\endtitle{\end{narrow}\newpage\pagenumstyle{arabic}}


\def\beginexercise{\vskip 20truept\parindent=0pt\begin{narrow}[10
truept]}
\def\endexercise{\vskip 10truept\end{narrow}}


\def\eql#1{\eqno\eqnlabel{#1}}
\def\ref{\reference}
\def\peq{\puteqn}
\def\pref{\putref}

\def\mgn{\marginnote}
\def\bex{\begin{exercise}}
\def\eex{\end{exercise}}


\font\open=msbm10 

 
\def\StretchRtArr#1{{\count255=0\loop\relbar\joinrel\advance\count255 by1
\ifnum\count255<#1\repeat\rightarrow}}
\def\StretchLtArr#1{\,{\leftarrow\!\!\count255=0\loop\relbar
\joinrel\advance\count255 by1\ifnum\count255<#1\repeat}}

\def\StretchLRtArr#1{\,{\leftarrow\!\!\count255=0\loop\relbar\joinrel\advance
\count255 by1\ifnum\count255<#1\repeat\rightarrow\,\,}}

\def\mbox#1{{\leavevmode\hbox{#1}}}

\def\hspace#1{{\phantom{\mbox#1}}}
\def\oZ{\mbox{\open\char90}}

\def\al{\alpha}

\def\bbe{{\bmit\beta}} 
\def\bka{{\bmit\kappa}}
\def\be{\beta}

\def\de{\delta}
\def\Ga{\Gamma}

\def\ep{\epsilon}

\def\ka{\kappa}
\def\la{\lambda}

\def\Om{\Omega}

\def\si{\sigma}

\def\th{\theta}
\def\Th{\Theta}
\def\ze{\zeta}

\def\De{\Delta}

\def\sc{{\rm sc }}

\def\zf{$\zeta$--function}


\def\frac#1/#2{\leavevmode\kern.1em
\raise.5ex\hbox{\the\scriptfont0 #1}\kern-.1em/\kern-.15em
\lower.25ex\hbox{\the\scriptfont0 #2}}
\def\sfrac#1/#2{\leavevmode\kern.1em
\raise.5ex\hbox{\the\scriptscriptfont0 #1}\kern-.1em/\kern-.15em
\lower.25ex\hbox{\the\scriptscriptfont0 #2}}

\def\gtorder{\mathrel{\raise.3ex\hbox{$>$}\mkern-14mu
             \lower0.6ex\hbox{$\sim$}}}
\def\ltorder{\mathrel{\raise.3ex\hbox{$<$}\mkern-14mu
             \lower0.6ex\hbox{$\sim$}}}

\def\semidirprod{\rlap{\ss C}\raise1pt\hbox{$\mkern.75mu\times$}}
\def\for{\lower6pt\hbox{$\Big|$}}
\def\fish{\kern-.25em{\phantom{abcde}\over \phantom{abcde}}\kern-.25em}


\def\boxit#1{\vbox{\hrule\hbox{\vrule\kern3pt
        \vbox{\kern3pt#1\kern3pt}\kern3pt\vrule}\hrule}}
\def\dalemb#1#2{{\vbox{\hrule height .#2pt
        \hbox{\vrule width.#2pt height#1pt \kern#1pt \vrule
                width.#2pt} \hrule height.#2pt}}}

\def\ol{\overline}
\def\frac#1#2{{{#1}\over{#2}}}

\def\noin{\noindent}

\def\comb#1#2{{\left(#1\atop#2\right)}}

\def\cot{{\rm cot\,}}

\def\etc{{\it etc. }}

\def\eg{{\it e.g.}}
\def\ie{{\it i.e. }}
\def\cf{{\it cf }}
\def\pa{\partial}



\def\3j#1#2#3#4#5#6{\left\lgroup\matrix{#1&#2&#3\cr#4&#5&#6\cr}
\right\rgroup}

\def\m?{\mgn{?}}

\def\pa{\partial}

\def\beq{\begin{eqnarray}}
\def\eeq{\end{eqnarray}}


\def\aop#1#2#3{{\it Ann. Phys.} {\bf {#1}} ({#2}) #3}
\def\cjp#1#2#3{{\it Can. J. Phys.} {\bf {#1}} ({#2}) #3}
\def\cmp#1#2#3{{\it Comm. Math. Phys.} {\bf {#1}} ({#2}) #3}
\def\cqg#1#2#3{{\it Class. Quant. Grav.} {\bf {#1}} ({#2}) #3}

\def\ijmp#1#2#3{{\it Int. J. Mod. Phys.} {\bf {#1}} ({#2}) #3}

\def\jmp#1#2#3{{\it J. Math. Phys.} {\bf {#1}} ({#2}) #3}
\def\jpa#1#2#3{{\it J. Phys.} {\bf A{#1}} ({#2}) #3}
\def\jpc#1#2#3{{\it J. Phys.} {\bf C{#1}} ({#2}) #3}
\def\lnm#1#2#3{{\it Lect. Notes Math.} {\bf {#1}} ({#2}) #3}

\def\np#1#2#3{{\it Nucl. Phys.} {\bf B{#1}} ({#2}) #3}
\def\npa#1#2#3{{\it Nucl. Phys.} {\bf A{#1}} ({#2}) #3}
\def\pl#1#2#3{{\it Phys. Lett.} {\bf {#1}} ({#2}) #3}

\def\prp#1#2#3{{\it Phys. Rep.} {\bf {#1}} ({#2}) #3}
\def\pr#1#2#3{{\it Phys. Rev.} {\bf {#1}} ({#2}) #3}
\def\prA#1#2#3{{\it Phys. Rev.} {\bf A{#1}} ({#2}) #3}

\def\prD#1#2#3{{\it Phys. Rev.} {\bf D{#1}} ({#2}) #3}
\def\prE#1#2#3{{\it Phys. Rev.} {\bf E{#1}} ({#2}) #3}
\def\prl#1#2#3{{\it Phys. Rev. Lett.} {\bf #1} ({#2}) #3}

\def\rmp#1#2#3{{\it Rev. Mod. Phys.} {\bf {#1}} ({#2}) #3}

\def\zfp#1#2#3{{\it Z. f. Phys.} {\bf {#1}} ({#2}) #3}

\def\cras#1#2#3{{\it Comptes Rend. Acad. Sci. (Paris)} {\bf{#1}} (#2) #3}
\def\prs#1#2#3{{\it Proc. Roy. Soc.} {\bf A{#1}} ({#2}) #3}
\def\pcps#1#2#3{{\it Proc. Camb. Phil. Soc.} {\bf{#1}} ({#2}) #3}
\def\mpcps#1#2#3{{\it Math. Proc. Camb. Phil. Soc.} {\bf{#1}} ({#2}) #3}

\def\amsh#1#2#3{{\it Abh. Math. Sem. Ham.} {\bf {#1}} ({#2}) #3}
\def\am#1#2#3{{\it Acta Mathematica} {\bf {#1}} ({#2}) #3}
\def\aim#1#2#3{{\it Adv. in Math.} {\bf {#1}} ({#2}) #3}
\def\ajm#1#2#3{{\it Am. J. Math.} {\bf {#1}} ({#2}) #3}
\def\amm#1#2#3{{\it Am. Math. Mon.} {\bf {#1}} ({#2}) #3}

\def\aom#1#2#3{{\it Ann. of Math.} {\bf {#1}} ({#2}) #3}
\def\cjm#1#2#3{{\it Can. J. Math.} {\bf {#1}} ({#2}) #3}
\def\bams#1#2#3{{\it Bull.Am.Math.Soc.} {\bf {#1}} ({#2}) #3}

\def\cmh#1#2#3{{\it Comm. Math. Helv.} {\bf {#1}} ({#2}) #3}

\def\dmj#1#2#3{{\it Duke Math. J.} {\bf {#1}} ({#2}) #3}
\def\invm#1#2#3{{\it Invent. Math.} {\bf {#1}} ({#2}) #3}

\def\jdg#1#2#3{{\it J. Diff. Geom.} {\bf {#1}} ({#2}) #3}

\def\joa#1#2#3{{\it J. of Algebra} {\bf {#1}} ({#2}) #3}
\def\jram#1#2#3{{\it J. f. Reine u. Angew. Math.} {\bf {#1}} ({#2}) #3}
\def\jims#1#2#3{{\it J. Indian. Math. Soc.} {\bf {#1}} ({#2}) #3}
\def\jlms#1#2#3{{\it J. Lond. Math. Soc.} {\bf {#1}} ({#2}) #3}
\def\jmpa#1#2#3{{\it J. Math. Pures. Appl.} {\bf {#1}} ({#2}) #3}
\def\ma#1#2#3{{\it Math. Ann.} {\bf {#1}} ({#2}) #3}

\def\mz#1#2#3{{\it Math. Zeit.} {\bf {#1}} ({#2}) #3}
\def\ojm#1#2#3{{\it Osaka J.Math.} {\bf {#1}} ({#2}) #3}

\def\pems#1#2#3{{\it Proc. Edin. Math. Soc.} {\bf {#1}} ({#2}) #3}

\def\plb#1#2#3{{\it Phys. Letts.} {\bf {B#1}} ({#2}) #3}
\def\pla#1#2#3{{\it Phys. Letts.} {\bf {A#1}} ({#2}) #3}
\def\plms#1#2#3{{\it Proc. Lond. Math. Soc.} {\bf {#1}} ({#2}) #3}
\def\pgma#1#2#3{{\it Proc. Glasgow Math. Ass.} {\bf {#1}} ({#2}) #3}
\def\qjm#1#2#3{{\it Quart. J. Math.} {\bf {#1}} ({#2}) #3}
\def\qjpam#1#2#3{{\it Quart. J. Pure and Appl. Math.} {\bf {#1}} ({#2}) #3}

\def\rmjm#1#2#3{{\it Rocky Mountain J. Math.} {\bf {#1}} ({#2}) #3}

\def\tams#1#2#3{{\it Trans.Am.Math.Soc.} {\bf {#1}} ({#2}) #3}

\begin{title}
\vglue 0.5truein
\vskip15truept
\centertext {\Bigfonts \bf On Sylvester waves and } \vskip2truept
\vskip10truept\centertext{\Bigfonts \bf restricted partitions }
 \vskip 20truept
\centertext{J.S.Dowker\footnote{dowker@man.ac.uk; dowkeruk@yahoo.co.uk}} \vskip
7truept \centertext{\it Theory Group,} \centertext{\it School of Physics and Astronomy,}
\centertext{\it The University of Manchester,} \centertext{\it Manchester, England} \vskip
7truept \centertext{}

\vskip 7truept

\vskip40truept
\begin{narrow}
The higher Sylvester waves are discussed. Techniques used involve finite difference
operators. For example,  using Herschel's theorem, elegant expressions for Euler's rational
functions and the Todd operator are found. Derivative expansions are also rapidly treated by
the same method.

A general form for the wave is obtained using multiplicative series, and comments made on
its further reduction. As is known, Dedekind sums arise in the case of coprime components
and it is pointed out that Brioschi had this result, but not the terminology,  very early on. His
proof is repeated.

Adding a set of ones to the components of the denumerant corresponds to a succession of
(discrete) smoothings and is a Cesaro sum. Using a spectral vocabulary, I take the
opportunity to exhibit the finite difference counterparts of some continuous, distributional
properties of Riesz typical means.

\end{narrow}
\vskip 5truept
\vskip 60truept
\vfil
\end{title}
\pagenum=0
\newpage

\section{\bf 1. Introduction}
This paper is a continuation of an earlier one, [\pref{Dowsyl}], concerned with Sylvester
waves, Ehrhart polynomials and degeneracies in spectral problems.\footnote{ The
appearance of the present work has been delayed by factors beyond my control.} In that
work I gave explicit formulae for the first two waves, resurrecting Sylvester's expressions.
Here I consider the higher waves but not so far as to produce similarly clear--cut results.
However, during the analysis, several amusing pieces of information and technique arose
and will be exposed here.

A classic number problem is that of restricted partitions. Given a set of non--negative
integers ${\bf d}_i$, $(i=1,2,\ldots,d)$, in how many ways can a given non--negative
integer, say $l$, be expressed as a (non-negative) integer linear combination of the
$d_i$? I call the $d_i$, the `components'  and  write this combination
  $$
 l = m_i\,d_i= {\bf m.d}\,.
 \eql{dioph}
  $$
The number of ways is referred to as the denumerant of $l$, following Sylvester, who was
one of the first to study this problem in any detail. Euler gives the number in terms of a
generating function
  $$
  \sum_{l=0}^\infty {l;\over {\bf d;}}\,\si^l=\prod_{i=1}^d{1\over1-\si^{d_i}}\,,
  \eql{eugen}
  $$
where the coefficient of the $\si^l$ is the denumerant. (I use the notation of
[\pref{Sylvester4}] to avoid typographical confusion.)

It is frequently possible, in discussions involving the combination (\peq{dioph}), to work
throughout, more or less, with this generating function. For example, in spectral
problems, where the denumerant might come up as a degeneracy, one could set
$\si=e^{-t}$ and interpret the left--hand side as a cylinder heat--kernel associated,
perhaps, with the square root of a Laplace--like operator, after some slight adjustments.
A Mellin transform would then yield the corresponding \zf\ which, in this case, would be a
Barnes \zf,
  $$\eqalign{
  \ze(s,a\mid{\bf d})&\equiv\sum_{l\ge0}{l;\over {\bf d;}}{1\over(a+l)^s}
  ={1\over\Ga(s)}\int_0^\infty dt\, t^{s-1}{e^{-at}\over \prod_i(1-e^{-d_it})}\cr
  &={i\Ga(1-s)\over2\pi}\int_L dt\,(- t)^{s-1}{e^{-at}\over \prod_i(1-e^{-d_it})}\,.
  }
  \eql{barnes}
  $$
(On this definition, the $d_i$ are not restricted to be positive integers.)

I will not pursue this global quantity further now because there is interest in studying the
denumerant on its own, as attested to by the long history.

Sylvester proves the theorem that the denumerant can be expressed as a finite sum of
quantities periodic in $l$, called waves, each of which is associated with a root of unity.
That associated with the root 1 is not periodic (has period infinity) and is a normal
polynomial, $W_1$,  in $l$. Hence one can set,
  $$
   {l;\over {\bf d};}=W_1+U
  $$
where $U$ is the properly periodic part and is a sum of waves, generally written
  $$
  U=W_2+\sum_{q>2} W_q\,.
  \eql{und}
  $$
The whole denumerant is a quasipolynomial in $l$ \ie a polynomial with periodic
coefficients, \eg\ [\pref{Stanley2}]. A simple proof of this was given by Wright,
[\pref{Wright}].

As a function of the augmented variable, $\ol l\equiv l+\sum_id_i/2$, the denumerant ${l;/
{\bf d};}$ satisfies a parity reciprocity under $\ol l\to -\ol l$, [\pref{Sylvester3}]. The
choice between $l$ and $\ol l$ comes up later.

The second wave $W_2$, is that associated with the root $-1$, and has also been
separated in (\peq{und})  because Sylvester has given an explicit form for $W_1$ and a
more or less explicit form for $W_2$. A simpler expression for $W_2$, like that for
$W_1$, was presented in our earlier work, [\pref{Dowsyl}], and I wish here to treat the
other waves by a similar technique.

The evaluation of the denumerant is a linear Diophantine question  related to the Frobenius
problem and to the Ehrhart polynomial and has thus been the subject of extensive analysis.
I refer to two books, [\pref{BandR}], [\pref{Alfonsin}], for some history and motivation.
However, Sylvester's specific formula is not referred to very often. The papers by Rubinstein
and Fel, [\pref{RandF}], and Rubinstein, [\pref{Rubinstein}], describe his method, but the
detailed techniques are different. A discussion of Sylvester's basic method can be found in
the book by Netto, [\pref{Netto}] published first in 1902, although he does not give
Sylvester's final form for $W_1$, which seems to be largely ignored, apart from the 1909
paper of Glaisher, [\pref{Glaisher2}], who extends the working to the other waves. The
present paper could be looked upon as an independent, partial recasting of Glaisher's
computations in as small a compass as I could manage, plus some comments about the
literature.\footnote{ Unfortunately, I have not been able to see some relevant early Italian
work by  Trudi.}

There is some very interesting recent work on computing these partition numbers
[\pref{SandZ}], [\pref{Munagi}].

\section{\bf2.  Sylvester's waves}

Sylvester's theorem leads to a prescription for the waves $W_q$ which is the following.
Write out all the factors of the components, $d_i$. Let there be $\mu$ such and call a
typical one, $q$. Separate the components into two groups -- those divisible by $q$, call
them $\al_i$, $i=1,2,\ldots\al$  ($\al$ is the `frequency') and those that are not, say
$\be_j$, $j=1,2,\ldots\be$ with $\al+\be=\mu$. Sylvester, [\pref{Sylvester2}], then
says that,
  $$
 W_q={\rm co}_{-1}\,\,\sum_\rho{\rho^{\,\ol l}\, e^{\,{\ol l}t}\over \prod_i
 \big(e^{\al_it/2}-e^{-\al_it/2}\big)\,
 \prod_j\big(\rho^{\be_j/2}e^{\be_jt/2}-\rho^{-\be_j/2}e^{-\be_jt/2}\big)}
 \eql{genfun3}
  $$
which stands for the coefficient of $1/t$ in the indicated generating function. The sum is
over all {\it prime} $q$th roots of unity, $\rho$, \ie $\rho^q=1$, $\rho=e^{2\pi ip/q}$,
$(p,q)=1$. The variable $\ol l$ is the augmented one, $\ol l\equiv l+\sum \al_i/2+\sum
\be_j/2=l+\sum d_i/2$.

The necessary inverse power(s) of $t$ are provided only by the $\al_i$ group. Therefore,
expanding the exponential, the coefficient of ${\ol l^n}/n!$ is
   $$\eqalign{
  &{\rm co}_{-1}\,\,\sum_\rho{\rho^{\, l}\, t^n\over \prod_i
 \big(e^{\al_it/2}-e^{-\al_it/2}\big)\,
 \prod_j\big(e^{\be_jt/2}-\rho^{-\be_j}e^{-\be_jt/2}\big)}\cr
 &={1\over\prod\al_i}\,
 \sum_\rho{\rho^{\,l}\over \prod_j(1-\rho^{-\be_j})}
 {\rm co}_{\al-n-1} \prod_i{ \al_it/2\over
 \sinh \al_it/2}\,
 \prod_j{\sin\be_j\pi p/q\over\sin\be_j(\pi p/q-it/2)}\,.
 }
 \eql{genfun4}
  $$
So far as the exponent of $\rho$ is concerned, it is best to retain the integer $l$, as Netto
[\pref{Netto}] suggests,  also $\rho^{q-1}+\rho^{q-2}+\ldots+1=0$.

As in my earlier work, I invoke some of Euler's old products. Specifically, after taking logs
  $$\eqalign{
  \log{z\over\sinh z}&=-\sum_{n=1}^\infty\log\bigg(1+{z^2\over n^2\pi^2}\bigg)=-
  {\ze(2)\over\pi^2}\,z^2-{\ze(4)\over2\pi^4}z^4-{\ze(6)\over3\pi^6}z^6-\ldots\cr
}
\eql{eulerp}
  $$
and
  $$\eqalign{
\log{\sin\pi a\over\sin\pi(z+a)}&=-
\sum_{n=1}^\infty\log\bigg(1-{z\over n-a}\bigg)+-\sum_{n=0}^\infty
\log\bigg(1+{z\over n+a}\bigg)\cr
&=-\lim\big(\ze(1,a)-\ze(1,1-a)\big)z+{1\over2}\big(\ze(2,a)+\ze(2,1-a)\big)z^2
  -\cr
  &\hspace{**}{1\over3}\big(\ze(3,a)-\ze(3,1-a)\big)z^3+
  {1\over4}\big(\ze(4,a)+\ze(4,1-a)\big)z^4-\cr
  &=-z\pi\cot\pi a-{1\over2}z^2{d\over da}\,\pi\cot\pi a-{1\over3}
  {1\over2!}z^3{d^2\over da^2}\,\pi\cot\pi a-\cr
  &=-\bigg(\int_0^z dz\,e^{z{d\over da}}\bigg)\,\pi\,\cot\pi a={1-e^{z\pa_a}\over\pa_a}
  \,\pi\,\cot\pi a\cr
  &=-\int_0^z dz\,\pi\,\cot\pi(z+ a)
}
\eql{eulerp2}
  $$
which is really only a check because it follows directly by school calculus.

Reversing the argument, the step from the second to the third line, which is the standard
reflection formulae of the Hurwitz \zf, (or, equivalently of the polygamma function),
   $$
   \ze(n,a)+(-1)^n\ze(n,1-a)={(-1)^{n-1}\over (n-1)!}{d^{n-1}
   \over da^{n-1}}\,\pi\,\cot \pi a\,,
   $$
can  be derived. For this nexus of notions consult Hoffman, [\pref{Hoffman,Hoffman2}],
who has a different organisation, and a little history.

I also note that one could expand in terms of the Euler rational functions, \eg\ Carlitz
[\pref{Carlitz2}], [\pref{RandF}], defined by,
  $$
  {1-\la\over\la -e^x}=\sum_n H_n(\la)\,{x^n\over n!}\,.
  \eql{eulerp3}
  $$
As a practical arithmetical means of calculation in any specific numerical case when the
number of components and the frequencies are not too large, this is probably as good a
route as any.  It was used by Cayley in his treatment of partitions, [\pref{Cayley5}]. He
gives two (related) relevant expansions, based on Herschel's theorem, \eg\
[\pref{Herschel}], [\pref{Boole}], [\pref{Milne-Thomson}]. Firstly, for small
denumerants,
  $$\eqalign{
  {1\over1-c\,e^{-t}}&=\sum_{f=0}^\infty{(-1)^f\over f!}\,t^f\,{1\over1-c(1+\De)}\,0^f\cr
  &={1\over1-c(1+\De)}\,e^{-0.t}
  }
  \eql{cayley}
  $$
which gives a direct expression for the Euler functions (\peq{eulerp3}) in finite
difference terms. The quantity  acting on $0^f$ is a finite polynomial in the difference
operator, $\De$, because $\De ^n\,0^m=0$ if $n>m$ and, generally,
  $$\eqalign{
  \De ^n\,0^m=(-1)^n\sum_{k=0}^n (-1)^k k^m\comb nk
  =n!\bigg\{{m\atop n}\bigg\}\,,
  }
  \eql{del0}
  $$
in terms of Stirling numbers of the second kind, a very old result. Thus I find the elegant
expression,
  $$
  H_n(\la)=-\bigg(1-{1\over\la-1}\De\bigg)^{-1}\,0^n\,,
  \eql{hn}
  $$
which I do not take any further here but note that it is equivalent to a result of Vandiver,
[\pref{Vandiver}]. See later.

Moreover, a form for the Todd operator, Todd$(c,\pa_h)$, of Brion and Vergne,
[\pref{BandV}], is
   $$\eqalign{
  {\pa\over1-c\,e^{-\pa_h}}={\log(1+\De)\over c(1+\De)-1}\,e^{-0.\pa_h}\,,
  }
  \eql{cayley5}
  $$
so that
  $$
  {\rm Todd}(c,\pa_h)\,f(h)={\log(1+\De)\over c(1+\De)-1}\,f(h-0)\,.
  \eql{Todd1}
  $$

For complicated denumerants, the logarithm is best expanded using (\peq{cayley}) and
  $$\eqalign{
  \log{1-c\over1-c\,e^{-t}}&=\int_0^t dt\bigg(1-{1\over1- c\,e^{-t}}\bigg)\cr
  &={1\over1-c(1+\De)}\,{1\over 0}\,\big(1-e^{-0.t}\big)\,.
  }
  $$

These particular expressions lead to the denumerant expressed in terms of the
unaugmented variable, $l$, and this is what Cayley produces after some calculation.
Interesting work on relieving some of this effort is given by Munagi, [\pref{Munagi}],
and Sills and Zeilberger, [\pref{SandZ}].

To encounter the (preferred) augmented variable, $\ol l$, a more symmetrical expansion
is required, corresponding to the cotangent in (\peq{eulerp2}). The relevant expansion
follows as
  $$
 {c e^{-t}\over1-c e^{-t}}={c (1+\De)\over1-c (1+\De)}\,e^{-0.t}\,.
 \eql{cotexp}
  $$

Set, now, $c=e^{-s}$ so that the left--hand side becomes
  $$
  {e^{-s-t}\over1- e^{-s-t}}
  $$
and the derivatives with respect to $t$, at $t=0$, which are explicit in (\peq{cotexp}), are
just the derivatives of $e^{-s}/(1-e^{-s})=(\coth s/2-1)/2$ with respect to $s$. This,
incidentally, provides an elegant difference proof of these higher derivatives,

\begin{ignore}
  $$\eqalign{
   {d^n\over ds^n}\coth s/2&=(-1)^n{1+e^{-s} (1+\De)\over1-e^{-s}(1+\De)}\,0^n\,\cr
   &=(-1)^n \bigg({2\over1-e^{-s} (1+\De)}-1\bigg)\,0^n\cr
    &=(-1)^n {2\over1-e^{-s} (1+\De)}\,0^n\,,\quad n>0\,\cr
    &= (-1)^n {2\over1-e^{-s}}\bigg(1-{e^{-s}\over1-e^{-s}} \De\bigg)^{-1}\,0^n\cr
    &=(-1)^n {2\over1-e^{-s}}\bigg(1+{e^{-s}\over1-e^{-s}} \De
    +\big({e^{-s}\over1-e^{-s}} \De\big)^2+\ldots\bigg)\,0^n\cr
    &=(-1)^n \big(\coth s/2+1\big)\sum_{k=1}^n {1\over2^k}\big(\coth s/2-1\big)^k
    \De^k\,0^n\cr
    &=(-1)^n \big(\coth s/2+1\big)\sum_{k=1}^n {k!\over2^k}\bigg\{{n\atop k}\bigg\}
    \big(\coth s/2-1\big)^k
   }
   \eql{dcoth}
  $$
  \end{ignore}
  
  $$\eqalign{
   {1\over2}{d^n\over ds^n}\bigg(\coth s/2-1\bigg)=
   {d^n\over ds^n}{1\over e^s-1}&=(-1)^n{1+\De
   \over e^s-(1+\De)}\,0^n\,\cr
   &=-(-1)^n \,\log\bigg(1-{\De\over e^s-1}\bigg)0^{n+1}\cr
    &=(-1)^n \sum_{k=1}^{n+1} {1\over k}\,{1\over (e^s-1)^k}\,\De^k\,0^{n+1}\cr
    &=(-1)^n \sum_{k=1}^{n+1} (k-1)!\bigg\{{n+1\atop k}\bigg\}
    \,{1\over (e^s-1)^k}\,,
   }
   \eql{dcoth2}
  $$
using (\peq{del0}) to give an explicit expression. The passage from the first to the second
line is effected by the identity due to Herschel (Boole, [\pref{Boole2}], Ch.II, \S13),
  $$
   \phi(\De)\,0^{n+1}=(1+\De)\,\phi'(\De)\,0^n\,,
  $$
whose essential content is the recurrence relation for the Stirling numbers.

The expansion, (\peq{dcoth2}) and its equivalents, surface from time to time. Agoh and
Dilcher, [\pref{AgandD}], prove it by induction and use it, in several papers, to derive
various Bernoulli number identities and Sterling convolutions.
  
It is equivalent to a formula in Adamchik, [\pref{Adamchik}]. See also Cvijovi\'{c},
[\pref{Cvijovic}], Knopf, [\pref{Knopf}], Hoffman, [\pref{Hoffman}], Boyadzhiev,
[\pref{Boyadzhiev}]. The technical relation used in these works is,
  $$
   \bigg(x{d\over dx}\bigg)^n f(x)\equiv {\ol {d^n}\over dx^n}f(x)=\sum_{k=1}^n
   \bigg\{{n\atop k}\bigg\}\,x^k\,{d^k\over dx^k}f(x)\,,
   \eql{knopf}
  $$
which Knopf, attributes, in essence, to Scherk in 1824 who, so it seems to me, was
addressing the same question solved by Herschel in 1816, arriving at the same result
\footnote{ Gould gives further interesting detail and earlier history in [\pref{Gould}].}.
Equation (\peq{knopf}) is contained in equation (4) in Herschel, [\pref{Herschel}], and
explicitly exhibited as Exercise 6, in Boole, [\pref{Boole}], p.26, {\it viz.}
  $$
   {\ol {d^n}\over dx^n}=\sum_{k=1}^n{1\over k!}\De^k0^n\,x^k{d^k\over dx^k}\,.
   \eql{her}
  $$

 For amusement, I interject a proof of this. From Herschel's theorem one has,
   $$
   f\big(e^{t+s}\big)=f\big(e^s(1+\De)\big)e^{\,0.t}\,,
   \eql{ht}
   $$
and note the oft used device,
  $$
  {d^n\over dt^n} f\big(e^{t+s}\big)\bigg|_{t=0}={d^n\over ds^n} f\big(e^{s}\big)
  ={\ol {d^n}\over dx^n} f(x)\,,\quad{\rm where}\quad x=e^s\,,
  $$
the left--hand side of which can be picked out of (\peq{ht}) as the coefficient of
$t^n/n!$.  The expansion of $f\big(x(1+\De)\big)$ in powers of $\De$,
  $$
   f\big(x(1+\De)\big)=\sum_{k=0}^\infty {1\over k!}\De^k0^n\, x^k{d^k\over dx^k}\,f(x)\,,
  $$
after retaining powers of $\De$ no greater than $n$ yields (\peq{her}). (The first term
also goes away because $0^n=0$.)

These considerations also show the equivalence of (\peq{hn}) with a formula in
Vandiver, [\pref{Vandiver}], who also employs the operator $\ol d/dx$.

After this digression, I return to the computation of the denumerant, \ie (\peq{genfun4}),
where one has to deal with the products. The first product has been encountered in
[\pref{Dowsyl}] and so I turn attention to the second one, to which (\peq{eulerp2}) applies
after the identifications, $a=a_j=\be_j p/q$ and $z=z_j=-i\be_j t/2$. I remark that $\pi
a_j$ is half the argument of $\rho^{\be_j}$, $\rho^{\be_j}=\exp(2\pi i\,a_j)$

Denoting by $\Om(z,a)$ the argument of the logarithm in (\peq{eulerp2}), the
corresponding multiplicative sequence follows on first constructing the product,
  $$\eqalign{
  \Om_1\Om_2\ldots=
  &\exp\bigg(-{t\over2}\,\Xi_1 +
  {1\over2}{t^2\over2^2}\,\Xi_2-
   {1\over3}{t^3\over2^3}{1\over2!}\,\Xi_3+
  {1\over4}{t^4\over2^4}{1\over3!}\,\Xi_4+\ldots\bigg)
  }
  \eql{prod}
 $$
where the $\Xi_n$ are defined by (\peq{eulerp2}).
  $$\eqalign{
  \Xi_n&=\sum_j i^n\,\be_j^n {d^{n-1}\over da^{n-1}}\,\cot\pi a\bigg|_{a=a_j}\,\cr
       &=2\sum_j (-1)^n\,\be_j^n\,{1\over1-\rho^{-\be_j}
       (1+\De)}\,0^{n-1}\,,\quad n>1\cr
       \Xi_1&=\sum_j\be_j\,{1+\rho^{\be_j}\over1-\rho^{\be_j}}
  }
  \eql{xin}
  $$
after employing (\peq{dcoth2}) or its equivalent.

\begin{ignore}

which I now re--express in terms of $\rho$ and note the usual relations
  $$
  i\cot\pi a ={1+\rho^\be\over1-\rho^\be}=2{\rho^\be\over1-\rho^\be}+1=
  {2\over1-\rho^\be}-1\,.
  \eql{cotexp}
  $$
and also
  $$
   {d\over da}=2\pi i \rho^\be{d\over d\rho^\be}\equiv 2\pi i{\ol d\over d\rho^\be}\,.
  $$

Therefore, for example,
  $$\eqalign{
   {d^{n}\over da^{n}}i\,\cot\pi a&= \bigg(2\pi i{\ol d\over d\rho^\be}\bigg)^{n}
   \bigg({2\over1-\rho^\be}-1\bigg)=\cr
   &=2(2\pi i)^n{1\over1-\rho^\be}\sum_{k=1}^n k!\,\bigg\{{n\atop k}\bigg\}
   \bigg({\rho^\be\over1-\rho^\be}\bigg)^k\,,\quad n>1\,,
   }
  $$
\end{ignore}

In terms of Stirling numbers, \cf\ (\peq{dcoth2}), (\peq{xin}) reads,
  $$\eqalign{
  &\Xi_n=2(-2\pi)^{n-1}\sum_j \be_j^n {1\over1-\rho^{\be_j}}
  \sum_{k=1}^{n-1} k!\,\bigg\{{{n-1}\atop k}\bigg\}
   \bigg({\rho^{\be_j}\over1-\rho^{\be_j}}\bigg)^k\,,\quad n>1\,.\cr
   }
  \eql{xin4}
  $$

Things can be expressed equivalently in terms of known Euler functions. Vandiver has
shown that, generally, \mgn{check history}
  $$
  {\ol{d^n}\over dx^n}\bigg({x\over 1-x}\bigg)={x\over1-x}\,H_n(x)
  $$
in terms of the functions defined by (\peq{eulerp3}) a few examples being $H_0=1$
and,
  $$
  H_1(x)={1\over x-1}\,,\quad
  H_2(x)={1-x\over(x-1)^2}\,,\quad
  H_3(x)={1+4x+x^2\over(x-1)^3}\,.
  \eql{eulerrat}
  $$
Setting $x=\rho^\be$ one has that,
   $$
  \bigg({\ol{d}\over d\rho^\be}\bigg)^n
  \bigg(i\,\cot\pi a-1\bigg)={2\rho^\be\over1-\rho^\be}\,H_n(\rho^\be)
  $$

Then,
  $$
  \Xi_n=2(2\pi)^{n-1}\sum_j \be_j^{n}{\rho^{\be_j}\over(1-\rho^{\be_j})^{n}}
  R_{n-1}\big(\rho^{\be_j}\big)\,,\quad n>1\,,
  \eql{xin3}
  $$
where $R_n(x)$ are Euler polynomials, the numerators of (\peq{eulerrat}). This is just
(\peq{xin}). A small simplification occurs if the inversion property,
  $$
  R_n(x)=(-1)^{n-1}x^{n-1}R_n(1/x)
  $$
is used. Then,
  $$
  \Xi_n=2(2\pi)^{n-1}\sum_j \be_j^{n}{1\over(1-\rho^{-\be_j})^{n}}
  R_{n-1}\big(\rho^{-\be_j}\big)\,,\quad n>1\,,
  \eql{xin2}
  $$
which can be obtained directly.

Any systematic analysis of Sylvester's theorem is bound to lead to similar quantities and
Glaisher, in his extensive treatment, encounters similar polynomials. See
[\pref{Glaisher2}], especially \S\S79,80,93-100. He also employs the difference operator
$\De$.

Equation (\peq{xin2}) yields an explicit, but unsimplified, formula for the product
(\peq{prod}), occurring in (\peq{genfun4}), which now has to be combined with the
product over the $\al_i$ coming via (\peq{eulerp}) as
   $$\eqalign{
   Q_1Q_2\ldots&=\exp\bigg(-s_2
  {\ze(2)\over\pi^2}\,{t^2\over2^2}+s_4{\ze(4)\over2\pi^4}{t^4\over2^4}
  -s_6{\ze(6)\over3\pi^6}{t^6\over2^6}+\ldots\bigg)\cr
  &\equiv\exp\bigg(-{\ol\tau_2\over2}\,{t^2}+{\ol\tau_4\over4}{t^4}
  -{\ol\tau_6\over6}{t^6}+\ldots\bigg)
  }
   $$
where $s_n$ is the sum of the $n$th powers of the $\al_i$ set. \footnote{ The relation
between the constants, $\ol \tau$, here and those, $\tau$, in [\pref{Dowsyl}] and
[\pref{Sylvester3}], is $\ol\tau_{2n}=2\tau_n$. This reflects the fact that the series
now contains all powers of $t$.}

Combining the two series gives,
  $$
  \Om_1\Om_2\ldots Q_1Q_2\ldots=\exp\bigg(-\ka_1t+
  {1\over2}\ka_2 t^2-{1\over3}\ka_3t^3+\ldots\bigg)
  $$
where
  $$\eqalign{
  \ka_1&={1\over2}\,\Xi_1\cr
  \ka_2&={1\over2^2}\,\Xi_2-\ol\tau_2\cr
  \ka_3&={1\over2^3\,2!}\,\Xi_3\cr
  \ka_4&={1\over2^4\,3!}\,\Xi_4+\ol\tau_4\cr
  &etc.
  }
  \eql{kappas}
  $$

The final step expands the exponential as a power series
  $$
  \Om_1\Om_2\ldots Q_1Q_2\ldots=1-\Th_1(\ka_1)\,t+
  \Th_2(\ka_1,\ka_2) t^2-\Th_3(\ka_1,\ka_2,\ka_3)\,t^3+\ldots
  $$
where the $\Th_r\big([{\bka}]_r\big)$ are all the homogeneous products of the
quantities of which the $\ka_i$ would be sums of powers and are classic functions of the
$\ka_i$, (see [\pref{Sylvester3}], [\pref{Dowsyl}]). Some examples are
  $$\eqalign{
  \Th_0=1\,,\quad \Th_1=\ka_1\,,\quad \Th_2={1\over2}\big(\ka_1^2+\ka_2\big)\,,
  \quad \Th_3={1\over6}\big(\ka_1^3+3\ka_1\ka_3+2\ka_3)\,.
   }
   \eql{aitchs}
  $$

Leaving the expressions as they are, returning to (\peq{genfun4})  for the polynomial
coefficient, the polynomial for the wave $W_q$ takes the finite form
  $$
  W_q={1\over\prod\al_i}\,
 \sum_\rho{\rho^{\,l}\over \prod_j(1-\rho^{-\be_j})}\bigg({{\ol l}^{\al-1}\over (\al-1)!}\,\Th_0
 -{{\ol l}^{\al-2}\over (\al-2)!}\,\Th_1+{{\ol l}^{\al-3}\over (\al-3)!}\,\Th_2-\bigg)\,.
 \eql{wq}
  $$

A trivial check is provided by  $W_2$, for which $\rho=-1$ and so, since the $\be_j$ are
all odd, the preliminary factor in (\peq{wq}) is just $1/(2^\be\prod\al) $. All the
$\Xi_n$ are zero and the result, using (\peq{kappas}), reduces to that in
[\pref{Dowsyl,Glaisher2}].

In this case the roots of unity dependence is trivial, which is not true for the general wave
and this is the remaining computational issue. A strategic decision has to be taken
concerning the `final' form for the periodic polynomial. Sylvester in [\pref{Sylvester4}]
writes it in terms of elementary denumerants, $(\ol l\pm r) ;/q;$, which first entails  a
reduction into a polynomial in the prime roots, $\rho$. This reduction also occurs in the
calculations of Cayley, [\pref{Cayley}] p.50, who expresses the final answer, equivalently,
in terms of prime circulators, as does Glaisher, whose calculation of $W_5$,
[\pref{Glaisher2}] \S\S 81-87, which he takes to the third term in (\peq{wq}) involving the
square of $\Xi_1$, brings out the attendant complications. Andrews, [\pref{Andrews}],
rewrites things in terms of the greatest integer function, which makes the integrality more
obvious, it being obscured in the other formulations. In, \eg\ Beck and Robins,
[\pref{BandR}], the roots of unity expressions are not taken further but are mostly left, and
analysed, as (generalised) Dedekind sums.\footnote{ Beck and Robins note that Dedekind
sums occur, implicitly, in the work of Sylvester and explicitly in Israilov, [\pref{Israilov}] in
1979.}

\section{\bf3. The simplest case of coprime components. Dedekind sums}
There is no need for complicated bookkeeping when the components, $d_i$, are mutually
prime. There is a subset of waves belonging to each factor on the denominator,
separately.

The expression for $W_1$ was given by Sylvester, [\pref{Sylvester3}], and was
rederived in [\pref{Glaisher2,Dowsyl}] so I consider the higher waves, $W_q$, where
$q$ divides the typical component $d_i$, $i=1,\ldots,d$, also treated in
[\pref{Glaisher2}]. Rather than write out this case specially, it is easier to refer to the
general form (\peq{genfun3}), set $\al=1$, $\al_1=d_i$. There is only one term on the
denominator that is proportional to $t$ so that $t$ can be set equal to $0$ in the
remaining terms. This easily gives,
   $$\eqalign{
 W_{q}&={1\over q}\sum_\rho{\rho^{ l}\over
 \prod_j\big(1-\rho^{-\be_j}\big)}={1\over q}\sum_\rho{\rho^{- l}\over
 \prod_j\big(1-\rho^{\be_j}\big)}\,,
 }
 \eql{genfun5}
  $$
by setting $\rho\to\rho^{-1}$. The $\be_i$ are the remaining components and the
$\rho$ are the non--trivial $q$--prime roots of unity, of which there are $q-1$, if $q$
happens to be prime when $\rho_m=e^{2\pi m/q}$ and also, in this case,
  $$\eqalign{
 W_{q}=
 {1\over q}\sum_{m=1}^{q-1}{\rho^{ l}_m\over
 \prod_j\big(1-\rho_m^{-\be_j}\big)}
 ={1\over q}\sum_{m=1}^{q-1}{\rho^{ -l}_m\over
 \prod_j\big(1-\rho_m^{\,\be_j}\big)}\,.
 }
 \eql{prime}
  $$

These $W_q$ are just the Fourier--Dedekind sums studied in [\pref{BandR}], and earlier
references therein. The actual definition there used is,
  $$
  s_l\big({\bbe};q\big)\equiv{1\over q}\sum_{m=1}^{q-1}{\rho_m^{\,l}\over\prod_j
  \big(1-\rho_m^{\,\be_j}\big)}\,,
  $$
so that it is only notation to write the wave, (\peq{prime}), as,
  $$
  W_{q}=s_{-l}\big(\bbe;q\big)\,.
  $$

The special case expression, (\peq{genfun5}), occurs in the early short note by Brioschi,
[\pref{Brioschi}], equn.(8), and so one should add his name to those who had already
encountered Dedekind sums. He also gives a contour proof of the general Sylvester
theorem and I repeat it in the Appendix for the coprime case.

Leaving the waves as Dedekind sums, (\peq{genfun5}) cannot be considered as a
complete determination and some energy has to be expended in computing the sums.
This is the purpose of the manipulations of Cayley, Glaisher and Sylvester, amongst
others. Brioschi, [\pref{Brioschi}], works out the particular example ${\bf d}=(2,3,5)$.
Sylvester, [\pref{Sylvester4}], does $(1,2,3), (1,4,7)$ and $(1,3,5)$, for those who like
variety. Remarks on the computability of the Dedekind sums are also made in
[\pref{BandR}] and \cf\ Exercise 8.2.

I consider the factor $\prod_j (1-\rho^{-\be_j})^{-1}$ where, to remind, $\rho$ is a
primitive $q$th root of unity and the $\be_j$ are the components of the denumerant not
divisible by $q$, $\be$ in number. It can be advantageous, as suggested by Cayley,
[\pref{Cayley}], and done by Glaisher, [\pref{Glaisher2}], to reduce the $\be_j$ mod
$q$. Then
  $$
  \prod_j (1-\rho^{-\be_j})^{-1}=\prod_{m=1}^{q-1}(1-\rho^{-m})^{-h_m}
  $$
where the non--negative integers $h_m$ depend on the set of the $\be_j$ and, in
general, have no pattern.

\section{\bf4. Roots of unity and prime circulators}
No point would be served by continuing the calculation along Glaisher's lines, which rapidly
becomes unwieldy.  Something more universal and automatic is required. Unfortunately I
cannot provide it here.

For example, for $W_q$, Glaisher reduces, mod $q$, the powers, $\be_j$ of $\rho$
before performing any combinations. In a general method this would be premature and
involve unnecessary extra labour.

Looking at the general structure (\peq{wq}), with (\peq{aitchs}), (\peq{kappas}) and
(\peq{xin}) or (\peq{xin2}), it is sufficient to write each term in the summation over $j$
in (\peq{xin2}) as a polynomial in $\rho$ of order $q-1$ (using $\rho^q=1$) and then
reduce any products of these (coming from the powers in (\peq{kappas})) to similar
polynomials. Ultimately, the $\Th_k$ are then also such polynomials and one then
converts $\Th_k/\prod_j (1-\rho^{-\be_j})$ to a like polynomial, which can then be
converted into Cayley's prime circulators, if desired.

The basic  algebraic problem is to reduce the ratio of polynomials of a primitive $q$th
root of unity, $\rho_\mu$, to a similar polynomial,
  $$
   {\sum_{m=0}^{q-1} A_m \,\rho_\mu^m\over\sum_{m=0}^{q-1} B_m \,\rho_\mu^m}
   =\sum_{m=0}^{q-1} C_m \,\rho_\mu^m\,,
  $$
\cf Battaglini (1857). Unfortunately I have not been able to see this work. It is mentioned
with a few details in Dickson, [\pref{Dickson3}] p.121.

At this rather unsatisfactory point I leave this detailed algebraic aspect.

\section{\bf5. Denumerants and Ces\`{a}ro sums}

For a given set of components, the numerical computation of a denumerant can proceed in
several ways. The classic special case when $d_1=1,d_2=2$ \etc, is treated at length in
Gupta [\pref{Gupta}], where some history is also given. The expressions derived by Cayley
(which are essentially the same as Sylvester's as enlarged by Glaisher, [\pref{Glaisher2}])
are considered. Large $l$ and largish $d$ are discussed and a typical example is detailed.
Cancellations and factorisations occur.

Euler computed many values using recursion, but for smallish $l$ and $d$ it is, perhaps,
just as easy to employ a convolution--iteration technique which yields an expression in
terms of simple denumerants only.\footnote{ The convolution corresponds to the product
form of the generating function, (\peq{eugen}), and the computation to the
time--honoured one of expanding each factor and collecting terms.}

The last component $d_d$ can be separated using the convolution
  $$
  {l;\over d_1,d_2,\ldots, d_d;}=\sum_{l'=0}^l\,{(l-l');\over d_d;}
   {l';\over d_1,d_2,\ldots, d_{d-1};}\,,
   \eql{denconv}
  $$
which is a special case of the more general,
  $$
  {l;\over d_1,d_2,\ldots, d_d;}=\sum_{l'=0}^l\,{(l-l');\over d_j,\ldots, d_d;}\,
   {l';\over d_1,d_2,\ldots, d_{j-1};}\,.
   \eql{denconv2}
  $$

Equation (\peq{denconv}) can be iterated to the intermediate form,
  $$
  {l;\over d_1,d_2,\ldots, d_d;}=
  \sum_{l_d=0}^l\,{(l-l_d);\over d_d;}\!\!
  \sum_{l_{d-1}=0}^{l_d}\,{(l_d-l_{d-1});\over d_{d-1};}
  \ldots\!\!\!\sum_{l_{d_{i+1}}=0}^{l_{d_{i+2}}}
  {(l_{d_{i+2}}-l_{d_{i+1}});\over d_{i+1};}
  {l_{d_{i+1}};\over d_1,\ldots, d_i;}
  \eql{intconv}
  $$
and, completely, down to
  $$
  {l;\over d_1,d_2,\ldots, d_d;}=
  \sum_{l_d=0}^l\,{(l-l_d);\over d_d;}
  \sum_{l_{d-1}=0}^{l_d}\,{(l_d-l_{d-1});\over d_{d-1};}
  \ldots\sum_{l_2=0}^{l_3}\,{(l_3-l_2);\over d_2;}{l_2;\over d_1;}\,.
  \eql{herden}
  $$

The simple denumerant, $l;/q;$, is 1 if $q$ divides $l$ and zero otherwise. I refer to this as
Herschel's function as it is just the average of the $l$th powers of all of the $q$th roots of
unity, which Herschel used when introducing his circulating functions, [\pref{Herschel2}]. In
terms of the fractional part and the Kronecker delta,
  $$
   {l;\over q;}=\de_{\{l/q\},0}\,,
  $$
and (\peq{herden}) is very easily programmed but, being recursive, is not very efficient.
However it does present the denumerant as an obvious integer.\footnote{ By this ancient
brute force method, I computed $100;/1,2,3,4,5;$ as 46262 in 2 minutes using DERIVE
and an Athlon IIx4 610e processor.}

From the spectral aspect, $l;/q;$, $l=0,1,\ldots$ is the Laplacian (Neumann) degeneracy
on the $2q$--divided circle (or, equivalently, on a $\pi/q$ interval) (eigenfunctions,
$\cos(lq\th)$). This can be envisaged as a one--dimensional lune.

The convolution, (\peq{denconv}), can be thought of as the addition of another
component. If the component 1 is added to the simple denumerant, $l;/q;$, one obtains
$l;/1,q;$ which is the degeneracy on the two--dimensional lune, $\lfloor l/q\rfloor+1$.

This can be repeated to give the degeneracy on the $d$--lune and I pursue this
particular process, but in the context of the intermediate form (\peq{intconv}) which
yields,
  $$
  {l;\over d_1,\ldots, d_{j-1},{\bf 1}_{i-j+1};}=
  \sum_{l_i=0}^l\,\,
  \sum_{l_{i-1}=0}^{l_i}
  \ldots\sum_{l_j=0}^{l_{j+1}}
  {l_j;\over d_1,\ldots, d_{j-1};}\,,
  \eql{intconv2}
  $$
which I think of as a succession of smoothings of the summand or as a nested series of
accumulated degeneracies.

All summations can be performed except the last, and one finds
  $$
  {l;\over d_1,\ldots, d_{j-1},{\bf 1}_{i-j+1};}=
  \sum_{l'=0}^l\comb{i-j+l-l'}{l-l'}
  {l';\over d_1,\ldots, d_{j-1};}\,,
  \eql{intconv3}
  $$
which is, to check, the same as (\peq{intconv3}), \ie,
  $$
  {l;\over d_1,\ldots, d_{j-1},{\bf 1}_{i-j+1};}=
  \sum_{l'=0}^l\,{(l-l');\over {\bf 1}_{i-j+1};}\,
   {l';\over d_1,d_2,\ldots, d_{j-1};}\,,
   \eql{denconv3}
  $$
in view of the classic value of the unit denumerant,
  $$
  {l;\over{\bf 1_d};}=\comb{d-1+l}{l}\,,
  $$
as follows, for example, from the expansion of the generating function,
  $$
  \sum_{l=0}^\infty {l;\over {\bf 1}_d;}\,\si^l={1\over1-\si^d}\,.
  \eql{eugen2}
  $$
This denumerant gives the degeneracy on the $d$--hemisphere.

The generating function version of (\peq{denconv3}) is the rather trivial splitting,
  $$
  {1\over (1-\si)^{i-j+1}(1-\si^{d_1})\ldots(1-\si^{d_{j-1}})}=
  {1\over (1-\si)^{i-j+1}}{1\over(1-\si^{d_1})\ldots(1-\si^{d_{j-1}})}\,.
  \eql{gfunsp}
  $$

It will be recognised that (\peq{intconv3}) is a Ces\`{a}ro sum originally introduced to
deal with divergent series. The standard situation (see \eg\ Hobson, [\pref{Hobson2}],
Bromwich, [\pref{Bromwich2}], Knopp, [\pref{Knopp3}]) is that for an infinite
sequence, $g_\nu$, one defines the finite sums
  $$
   S_n^{(r)}=\sum_{\nu=0}^n \comb{r+n-\nu}r\,g_{\nu}
   \eql{cesum}
  $$
and investigates the $r$th `arithmetic mean'
  $$
 \comb {r+n}n^{-1} S_n^{(r)}
  $$
as $n\to\infty$. In this expression,  the index, $r$, although it originates as the number
of smoothing summations, can assume any value, real or complex. Negative integers,
however, are usually excluded (but see later).

In terms of generating functions, (\peq{cesum})  translates into
  $$
  \sum_{\nu=0}^\infty  S_\nu^{(r)}\,\si^\nu=
  {1\over(1-\si)^{r+1}}\sum_{\nu=0}^\infty g_\nu \si^\nu
  \eql{cesum2}
  $$
of which (\peq{gfunsp}) is an example.

The Ces\`{a}ro mean is a discrete analogue of,  and a motivation for, the more powerful
`typical mean' of Riesz, which introduces a handy continuous variable into the analysis, \eg\
Hardy and Riesz, [\pref{HaandR}], [\pref{Hobson2}]. I give only the briefest details in a
more general setting.

Typically, in an eigenproblem, $\nu$ would be an eigen{\it level} label, and $g_\nu$ its
degeneracy. A function $N(\la)$ of the continuous variable, $\la$, is defined as a
counting function, encoding the spectrum $\la_\nu$,
  $$
  N(\la)=\sum_{\nu,\la\le\la_\nu}\,g_\nu\,.
  $$
In the special situation above, $\nu$ would be $l$, $g_\nu$ the denumerant, $l;/{\bf
d};$ and the eigenvalue, $\la_\nu$, a function of $l$, typically $a+l$, where $a$ is a
constant.

Defining the (first) accumulated degeneracy by,
  $$
   G_\nu=\sum_{\nu'=0}^\nu g_{\nu'}
  $$
or,
  $$
   \sum_{\nu=0}^\infty G_\nu\,\si^\nu=
   {1\over 1-\si}\sum_{\nu=0}^\infty g_{\nu}\si^\nu\,,
  $$
one has the connection,
  $$\eqalign{
  N(\la)&=G_{\nu-1},\quad \la_{\nu-1}<\la<\la_\nu\cr
    &=G_{\nu-1}+{1\over2}g_\nu=G_{\nu+1}-{1\over2}g_\nu
    \,,\quad \la=\la_\nu\,,
  }
  $$
and, of course, $G_\nu$ is the first Ce\`saro sum $S^{(0)}_\nu$, (\peq{cesum}).

Unlike the accumulated degeneracies, $N(\la)$ depends on the actual form of the
eigenvalues, $\la_\nu$.

The Ces\`aro generating function, (\peq{cesum2}) in terms of $G$, is
  $$
  \sum_{\nu=0}^\infty  S_\nu^{(r)}\,\si^\nu=
  {1\over(1-\si)^{r}}\sum_{\nu=0}^\infty G_\nu \si^\nu\,.
  \eql{cesum3}
  $$
For the {\it Riesz} mean, this iterated summation of $G_\nu$ is replaced by an
$(r+1)$--fold iterated integration of $N(\la)$ which equals the Cauchy convolution, (see
Knopp, [\pref{Knopp3}]),
  $$
  N_{r+1}(\la)=\int_0^\infty d\la' {(\la-\la')^{r}\over \Ga(r+1)}\,N(\la')\,,
  $$
employed specifically as a smoothing of $N$ by Fedosov, [\pref{Fedosov}], see Baltes
and Hilf, [\pref{BaandH}], Balian and Bloch [\pref{BaandB}].

Introducing the distribution $\Phi$, (see Gel'fand and Shilov, [\pref{GeandS}] \S5.5),
  $$
   \Phi_\al(x)={x_+^{\al-1}\over\Ga(\al)}\,,
  $$
the convolution is neatly written,
  $$
  N_\al=N*\Phi_\al\,,
  $$
 at the same time extending $r+1$ to an arbitrary variable, $\al$, by continuation.

 Comparing this convolution with the Ces\`aro discrete version, (\peq{cesum})
 or (\peq{intconv3}), leads to the analogy
   $$
  \comb{l+\al-1}{l}={\Ga(l+\al)\over\Ga(l+1)\Ga(\al)}\sim \Phi_\al(l)\,.
  \eql{anal}
   $$
For typographical reasons, the left--hand side is often denoted by $A_l^{\al-1}$.

It is interesting to enlarge upon this analogy and to exhibit some continuous relations
together with their finite difference counterparts.

The distribution $\Phi$ has a number of basic properties, [\pref{GeandS}].

\vskip 5truept

(i) The concentration of $\Phi_\al(x)$ on the positive $x$ axis.

(ii) The convolution,
  $$
  \Phi_\al*\Phi_\be=\Phi_{\al+\be}\,.
  $$

(iii) The Laplace transform,
  $$
  \int_0^\infty dx\, \Phi_\al(x)\,e^{-tx}={1\over t^\al}\,,\quad t>0\,.
  $$

(iv) The singularities,
  $$
  \Phi_{-k}(x)=\de^{(k)}(x)\,.
  $$

Property (i) corresponds to the statement that $A_l^{\al-1}=0$ if $l$ is a negative
integer (because of the poles in the $\Ga(l+1)$--function) unless $\al$ is a positive
integer, when it vanishes at a finite number of negative integers.  $\al$ a negative
integer comes under property (iv).

The convolution property (ii) is a simple consequence of the definition and analytic
continuation. The corresponding discrete equation {\it viz}, the well known binomial
identity,
  $$
  \sum_{\nu=0}^lA_\nu^{\al-1}\,A_{l-\nu}^{\be-1}=A_l^{\al+\be-1}\,,
  $$
follows likewise, most easily from factorisation of the generating function. This is just
Cauchy's product.

Property (iii) translates into the generating function definition,
  $$
  \sum_{l=0}^\infty A_l^{\al-1}\,\si^l={1\over\big(1-\si\big)^\al}\,,\quad \si=e^{-t},\,\,t>0.
  \eql{adef}
  $$

Property (iv) is the most interesting one as it concerns the case when $\al$ is a negative
integer which is usually excluded (\eg\ Chapman, [\pref{Chapman}]). Looking at
(\peq{anal}), the quantity of interest is,
  $$
  \lim_{\al\to -k=0,-1,\ldots}{\Ga(l-\al)\over\Ga(l+1)\Ga(\al)}=
  (-1)^l{\Ga(k+1)\over\Ga(l+1)\Ga(k-l+1)}=(-1)^l\comb kl\,.
  \eql{negal}
  $$
This also follows quickly from the binomial expansion of (\peq{adef}) and one could use
generating functions systematically, \eg\ Jordan, [\pref{Jordan}], to re--express the
following remarks. \footnote{ The isomorphic algebraic  scheme of generalised
sequences, \eg\ [\pref{Moore}], [\pref{Traub}], could also be employed. This forms an
operator calculus in the field of finite differences.}

I recall now the expression for the $k$th difference, \eg\ Boole, [\pref{Boole}], applied
to the Kronecker delta, $\de_l^{\,\,l'}$, considered as a function of $l$ (\cf\ Tauber and
Dean, [\pref{TandD}], Traub, [\pref{Traub}] \footnote{ There appears to be an overall
sign error in equations (3.6) and (3.7) of [\pref{TandD}].}). The expressions are,
  $$\eqalign{
  {\overrightarrow\De}^k\,\de_l^{\,\,l'}=\sum_{m=0}^k (-1)^m
  \comb km\,\de_{k-m+l}^{\,\,l'}=(-1)^{l'-l-k}\comb k{l'-l}\cr
   {\overleftarrow\De}^k\,\de_l^{\,\,l'}=\sum_{m=0}^k (-1)^{m}
   \comb km\,\de_{-m+l}^{\,\,l'}=(-1)^{l-l'}\comb k{l-l'}
  }
  \eql{diffdel}
  $$
where ${\overrightarrow\De}$ is the forward difference and ${\overleftarrow\De}$ the
backwards one. Comparing with (\peq{negal}) the conclusion is that,
  $$
  \lim_{\al\to -k} A^{\al-1}_{l-l'}={\overleftarrow\De}^k\,\de_{l}^{\,\,l'}=
  (-1)^k{\overrightarrow\De}^k\,\de_{l'}^{\,\,l}
  $$
which is the finite difference counterpart of property (iv), as I wished to show.

Unlike the continuous version, the higher difference
${\overleftarrow\De}^k\,\de_{l}^{\,\,0}$ is not concentrated at the origin, $l=0$, but
extends to the right for a `distance' $k$, emphasising the non--locality, or fuzziness, in
its construction. As a simple feature, I remark that, $\sum_{l=0}^k
{\overleftarrow\De}^k\,\de_{l}^{\,\,0}=\de_k^0$, which just states that the strength
of the (oscillating) `curve' of a higher derivative is zero. The corresponding forwards
quantity, ${\overrightarrow\De}^k\,\de_{l}^{\,\,0}$, extends a distance $k$ to the {\it
left} of the origin.

As usual, employing forwards or backwards differences involves a manifest loss of
symmetry which can be restored by shifting the origin of the graphs to their mid points.
Analytically this is accomplished by the translation operator, conventionally written ${\bf
E}^{\pm1/2}$, and corresponds to using the central difference, $\de$, so that
$\de^k\de_l^0$ is the closer analogue of the Dirac derivative $\de^{(k)}(x)$.\footnote{
There is, of course, an {\it intrinsic} non--locality in the continuous higher derivative.}

Pictures of some discontinuous approximations to $\de$--function derivatives obtained by
central differences of the (piecewise) continuous step function can be seen in van der Pol and
Bremmer, [\pref{PandB}], p.83.\footnote{ If the blocks there depicted are squashed to
their midpoints, and $\ep=1/2$, then these graphs yield precisely the symmetrical
constructions, $\de^{k}\de_l^0$, here.}

\section{\bf 6. Discussion}
I have given several disparate pieces of analysis based on the explicit calculation of a
denumerant, or restricted partition. On the way a neat finite difference expression,
(\peq{hn}), for Euler's functions was found and also one for the Todd operator,
(\peq{Todd1}). An elegant derivation of a derivative expansion is likewise given.

An expression for the general wave $W_q$ is given, (\peq{wq}), but this is not necessarily
in its final form.

The expansion of the denumerant into waves can be sustituted into the Barnes \zf,
(\peq{barnes}), although would not aid its specific computation.

Brioschi's early 1857 calculation of the simplest case of prime components is resurrected and
leads to Dedekind sums, which is known more recently.

\section{\bf Appendix  A. Proof of Brioschi's formula}

For completeness, I give a pedagogic derivation of equation (\peq{genfun5}) for the wave
in the simplest case when the components are all prime. There is nothing new in this, \cf\
Brioschi in 1857, [\pref{Brioschi}].

The starting point is always Euler's generating function for the denumerant,
  $$
   \sum_{l=0}^\infty {l;\over{\bf d;}}\,z^l
   ={1\over\prod_i(1-z^{d_i})}={1\over(1-z^\al)\prod_i(1-z^{\be_j})}\equiv
   {1\over(1-z^\al)\,f(z)}\,.
   \eql{egenfun}
  $$
where $\al$ is one component selected from the $d_i$ and $\be_j$ the rest. Notationally
I set $\al=q\in\oZ$. Extracting the power $z^l$ using residue calculus gives
  $$
  {l;\over{\bf d;}}={1\over2\pi i}\oint_C dz{1\over z^{l+1}}{1\over(1-z^q)\,f(z)}
  \eql{contour}
  $$
where $C$ circles the origin. Now blow up $C$ to wrap round the other poles, which all
lie on the unit circle and are given by the vanishing places of the denominator in
(\peq{egenfun}). In this sum, I consider only those arising from from the vanishing of
the separated typical factor $(1-z^q)$ at all the {\it non--trivial} $q$th roots of unity,
$\rho_1,\ldots,\rho_{q-1}$. This particular sum equals, by definition, the wave, $W_q$.
The contour then continues on to infinity, which contributes nothing.

The residue of $1/(1-z^q)$ at $z=\rho_j$ is easily found to be $-\rho_j/q$ so that the
residue of the integrand is $\rho_j^{-l}/q f(\rho_j)$, and I thence obtain\mgn{check
minus signs} (\peq{genfun5}).

The multiple pole at $z=1$, coming from all factors in the denominator would give the first
(non--periodic) wave, $W_1$, a polynomial, essentially a generalised Bernoulii polynomial,
\eg\ [\pref{Rubinstein}]. There are no other multiple poles, which accounts for the
simplicity of this evaluation.\footnote{ Because of the separation of the roots, one might say
that there is no interference between the waves.} The total denumerant is then the sum
  $$
  {l;\over{\bf d;}}=W_1+\sum_{i=1}^{d}W_{d_i}\,.
  \eql{totden}
  $$

If the components are only coprime, the roots $\rho_i$, for {\it each} component, separate
into the prime roots for every factor of that component so giving a subset of waves, one
wave, of the form (\peq{genfun5}), for each divisor. The total denumerant is still given by
(\peq{totden}) with the subsum,
  $$
  W_{d_i}=\sum_{q/d_i} W_q\,.
  \eql{coprime}
  $$

The general case, when some components possess common factors, involves coincident
poles on the unit circle originating from different factors on the denominator. These, like
the poles at $z=1$, give rise to polynomials. Analysing and  organising this situation
leads to Sylvester's theorem, [\pref{Brioschi}].

Ehrhart, [\pref{Ehrhart}], Theorem 9.2, also gives the coprime expression for the
denumerant (`compteur') which he seems to derive independently. He uses it
numerically to compute the denumerant by calculating the expression to an adequate
approximation. The bulk of the value comes from the first wave, polynomial part.

 \vglue 15truept

\noin{\bf References.} \vskip5truept

\begin{putreferences}
   \ref{BaandB}{Balian,R. and Bloch,C. \aop{60}{1970}{401}.}
   \ref{AgandD}{Agoh,T. and Dilcher,K. {\it J. Number Theory} {\bf 124} (2007) 105.}
   \ref{Andrews}{Andrews,G.E. {\it Ramaujan J.} {\bf 7} (2003) 385.}
    \ref{Glaisher2}{Glaisher,J.W.L. \qjm {40}{1909}{275}.}
    \ref{RandF}{Rubinstein, B.Y. and Fel,L.G., {\it Ramanujan J.} {\bf11}(2006)331.}
    \ref{Rubinstein}{Rubinstein, B.Y. {\it Ramanujan J.} {\bf15}(2008)177.}
   \ref{Dickson3}{Dickson,L.E. {\it History of the Theory of Numbers} vol.2.
   (Carnegie Institute, Washington, 1920).}
   \ref{Brioschi}{Brioschi,F. {\it Annali di sc. mat. e fis.} {\bf8} (1857) 5.}
   \ref{Fedosov}{Fedosov,B.V. {\it Sov.Math.Dokl.} {\bf 5} (1963) 1992}
   \ref{Chapman}{Chapman}
    \ref{Jordan}{Jordan,C. {\it Calculus of Finite Differences}, (Budapest, 1939).}
   \ref{PandB}{van der Pol, B. and Bremmer,H. {\it Operational Calculus} (CUP,
   Cambridge, 1964).}
   \ref{TandD}{Tauber,S. and Dean,D. {\it J.SIAM} {\bf 8} (1960) 174.}
   \ref{Traub}{Traub,J.F. {\it Math. of Comp.} {\bf 19} (1965) 177.}
   \ref{Chapman}{Chapman,S. \plms{9}{1911}{369}. }
   \ref{Moore}{Moore,D.H. \amm{69}{1962}{132}.}
   \ref{HaandR}{Hardy,G.H. and Riesz,M. {\it General theory of Dirichlet's series}
   (CUP, Cambridge, 1915).}
   \ref{Gupta}{Gupta,H. {\it Tables of partitions} (Royal.Soc.Math.Tables 4)
   (Cambridge, 1958).}
   \ref{Bromwich}{Bromwich, T.J.I'A. {\it Infinite Series},
  (Macmillan, 1947).}
  \ref{Hobson2}{Hobson,E.W. {\it Theory of functions of a real variable}, vol.2.
  (C.U.P., Cambridge,1907).}
   \ref{Bromwich2}{Bromwich, T.J.I'A. {\it Infinite Series}, 1st Edn.
  (Macmillan, London,1907).}
   \ref{Knopp3}{Knopp,K. {\it Theory of Infinite Series} (Blackie, London, 1928).}
   \ref{Vandiver}{Vandiver,H.S. \tams{51}{1942}{502}.}
   \ref{Adamchik}{Adamchik,V.S.{\it Appl.Math. and Comp.} {\bf 187} (2007) 3.}
   \ref{Cvijovic}{Cvijovi\'{c},D. {\it Appl.Math. and Comp.} {\bf 215} (2009) 3002.}
   \ref{Carlitz2}{Carlitz,L. {\it Math.Magazine} {\bf 32} (1959) 247.}
   \ref{Israilov}{Israilov.M.I. {\it Sibirsk Mat. Zh.} {\bf 22} (1981) 121.}
   \ref{SandZ}{Sills,A.V. and Zeilberger,D. {\it Formulae for the number of partitions of
    $n$ into at most $m$ parts (using the quasi-polynomial ansatz)} ArXiv:1108.4391.}
   \ref{BandV}{Brion,M. and Vergne,M. {\it J.Am.Math.Soc.} {\bf 10} (1997) 371.}
   \ref{Knopf}{Knopf,P.M. {\it Math. Magazine} {\bf 76} (2003) 364.}
   \ref{Gould}{Gould,H.W. {\it Am. Math. Monthly} {\bf 86} (1978) 450.}
   \ref{Hoffman}{Hoffman,M.E. \amm{102}{1995}{23}.}
   \ref{Hoffman2}{Hoffman,M.E. {\it Electronic J.Comb.} {\bf 6} (1999) \#R21.}
   \ref{Boyadzhiev}{Boyadzhiev, K.N. {\it Fibonacci Quarterly} {\bf45} (2008) 291.}
   \ref{Sills}{Sills}
   \ref{Munagi}{Munagi,A.O. {\it Electronic J. Comb. Number Th.} {\bf 7} (2007) \#A25.}
   \ref{Dowsyl}{Dowker,J.S. {\it Relations between Ehrhart polynomials, the heat kernel
   and \break Sylvester waves} ArXiv:1108.1760}
   \ref{Alfonsin}{Alfonsin,J.L.R. {\it The Diophantine Frobenius Problem} (O.U.P., Oxford, 2005).}
   \ref{GeandS}{Gel'fand,I.M. and Shilov,G.E. {\it Generalized Functions} Vol.1
   (Academic Press, N.Y., 1964).}
    \ref{Boole}{Boole, G. {\it Calculus of Finite Differences}, (MacMillan, Cambridge,
1860).}
   \ref{Boole2}{Boole, G. {\it Calculus of Finite Differences}, 2nd Edn. (MacMillan, Cambridge,
1872).}
     \ref{Cayley5}{Cayley, A.  {\it Phil. Trans. Roy. Soc. Lond.} {\bf 148} (1858) 47.}
     \ref{Milne-Thomson}{Milne-Thomson, L.M. {\it The Calculus of Finite Differences},
     (MacMillan, London, 1933).}
    \ref{Herschel}{Herschel, J.F.W.  {\it Phil. Trans. Roy. Soc. Lond.} {\bf 106} (1816) 25.}
    \ref{Herschel2}{Herschel, J.F.W.  {\it Phil. Trans. Roy. Soc. Lond.} {\bf 108} (1818) 144.}
    \ref{Littlewood2}{Littlewood,D.E. {\it The Theory of Group Characters}
    (Clarendon Press, Oxford, 1950).}
    \ref{Wright}{Wright,E.M. \amm{68}{1961}{144}.}
\ref{Carlitz}{Carlitz,L. \dmj{27}{1960}{401}.}
     \ref{Netto}{Netto,E. {\it Lehrbuch der Combinatorik} 2nd Edn. (Teubner, Leipzig, 1927).}
    \ref{FdeB}{Fa\`{a} de Bruno, F. {\it Th\'eorie des Formes Binaires} (Brero, Turin,1876).}
    \ref{Ehrhart}{Ehrhart,E. \jram{227}{25}{1967}.}
    \ref{Bell}{Bell,E.T.\ajm{65}{1943}{382}.}
    \ref{BandR}{Beck,M. and Robins,S. {\it Computing the Continuous Discretely,}
    (Springer, New York, 2007).}
    \ref{BandR2}{Beck, M. and Robins,S. {\it Discrete and Comp. Geom.} {\bf 27}(2002) 443.}
    \ref{Harmer}{Harmer,M. {\it J.Australian Math.Soc.} {\bf 84}(2008)217.}
    \ref{RandF}{Rubinstein, B.Y. and Fel,L.G., {\it Ramanujan J.} {\bf11}(2006)331.}
    \ref{BGK}{Beck,M., Gessel, I.M. and Komatsu,T. {\it Electronic Journal of Combinatorics}
    {\bf8}(2001) 1.}
    \ref{Sylvester}{Sylvester,J.J. \qjpam{1}{1858}{81}.}
    \ref{Sylvester2}{Sylvester,J.J. \qjpam{1}{1858}{142}.}
    \ref{Sylvester3}{Sylvester,J.J. \ajm{5}{1882}{79}.}
    \ref{Sylvester4}{Sylvester,J.J. \plms{28}{1896}{33}.}
    \ref{Dowgta}{J.S.Dowker, {\it Group theory aspects of spectral problems on spherical factors},
   ArXiv.Math.DG: 0907.1309.}
    \ref{BDR}{Beck,M., Diaz and Robins,S. {\it J.Numb.Theory} {\bf 96} (2002) 1.}
    \ref{PandS}{P\'{o}lya, G. and Szeg\H{o},G. {\it Aufgaben und Lehrs\"atze aus der Analysis}
    (Springer--Verlag, Berlin, 1925).}
    \ref{EOS}{Elizalde,E., Odintsov, S.D. and Saharian, A.A. \prD{79}{2009}{065023}.}
    \ref{Cavalcanti}{Cavalcanti,R.M. \prD{69}{2004}{065015}.}
    \ref{MWK}{Milton, K.A., Wagner,J. and Kirsten,K. \prD{80}{2009}{125028}.}
    \ref{EBM2}{Ellingsen,S.A., Brevik,I. and Milton,K.A. \prE{81}{2010}{065031}.}
    \ref{EBM}{Ellingsen,S.A., Brevik,I. and Milton,K.A. \prE{80}{2009}{021125}.}
    \ref{BEM}{Brevik,I., Ellingsen,S.A. and Milton,K.A. \prE{79}{2009}{041120}.}
    \ref{FKW}{Fulling,S.A, Kaplan L. and Wilson,J.H. \prA{76}{2007}{012118}.}
    \ref{Lukosz}{Lukosz,W, {\it Physica} {\bf 56} (1971) 109; \zfp{258}{1973}{99}
    ;\zfp{262}{1973}{327}.}
    \ref{Gromes}{Gromes, D. \mz{94}{1966}{110}.}
    \ref{FandK1}{Kirsten,K. and Fulling,S.A. \prD{79}{2009}{065019} .}
    \ref{FandK2}{Fucci,G. and Kirsten,K, JHEP (2011), 1103:016.}
    \ref{dowgjms}{Dowker,J.S. {\it Determinants and conformal anomalies
    of GJMS operators on spheres}, ArXiv: 1007.3865.}
    \ref{Dowcascone}{dowker,J.S. \prD{36}{1987}{3095}.}
    \ref{Dowcos}{dowker,J.S. \prD{36}{1987}{3742}.}
    \ref{Dowtherm}{Dowker,J.S. \prD{18}{1978}{1856}.}
    \ref{Dowgeo}{Dowker,J.S. \cqg{11}{1994}{L55}.}
    \ref{ApandD2}{Dowker,J.S. and Apps,J.S. \cqg{12}{1995}{1363}.}
   \ref{HandW}{Hertzberg,M.P. and Wilczek,F. {\it Some calculable contributions to
   Entanglement Entropy}, ArXiv:1007.0993.}
   \ref{KandB}{Kamela,M. and Burgess,C.P. \cjp{77}{1999}{85}.}
   \ref{Dowhyp}{Dowker,J.S. \jpa{43}{2010}{445402}; ArXiv:1007.3865.}
   \ref{LNST}{Lohmayer,R., Neuberger,H, Schwimmer,A. and Theisen,S.
   \plb{685}{2010}{222}.}
   \ref{Allen2}{Allen,B. PhD Thesis, University of Cambridge, 1984.}
   \ref{MyandS}{Myers,R.C. and Sinha,A. {\it Seeing a c-theorem with
   holography}, ArXiv:1006.1263}
   \ref{MyandS2}{Myers,R.C. and Sinha,A. {\it Holographic c-theorems in
   arbitrary dimensions},\break ArXiv: 1011.5819.}
   \ref{RyandT}{Ryu,S. and Takayanagi,T. JHEP {\bf 0608}(2006)045.}
   \ref{CaandH}{Casini,H. and Huerta,M. {\it Entanglement entropy
   for the n--sphere},\break arXiv:1007.1813.}
   \ref{CaandH3}{Casini,H. and Huerta,M. \jpa {42}{2009}{504007}.}
   \ref{Solodukhin}{Solodukhin,S.N. \plb{665}{2008}{305}.}
   \ref{Solodukhin2}{Solodukhin,S.N. \plb{693}{2010}{605}.}
   \ref{CaandW}{Callan,C.G. and Wilczek,F. \plb{333}{1994}{55}.}
   \ref{FandS1}{Fursaev,D.V. and Solodukhin,S.N. \plb{365}{1996}{51}.}
   \ref{FandS2}{Fursaev,D.V. and Solodukhin,S.N. \prD{52}{1995}{2133}.}
   \ref{Fursaev}{Fursaev,D.V. \plb{334}{1994}{53}.}
   \ref{Donnelly2}{Donnelly,H. \ma{224}{1976}{161}.}
   \ref{ApandD}{Apps,J.S. and Dowker,J.S. \cqg{15}{1998}{1121}.}
   \ref{FandM}{Fursaev,D.V. and Miele,G. \prD{49}{1994}{987}.}
   \ref{Dowker2}{Dowker,J.S.\cqg{11}{1994}{L137}.}
   \ref{Dowker1}{Dowker,J.S.\prD{50}{1994}{6369}.}
   \ref{FNT}{Fujita,M.,Nishioka,T. and Takayanagi,T. JHEP {\bf 0809}
   (2008) 016.}
   \ref{Hund}{Hund,F. \zfp{51}{1928}{1}.}
   \ref{Elert}{Elert,W. \zfp {51}{1928}{8}.}
   \ref{Poole2}{Poole,E.G.C. \qjm{3}{1932}{183}.}
   \ref{Bellon}{Bellon,M.P. {\it On the icosahedron: from two to three
   dimensions}, arXiv:0705.3241.}
   \ref{Bellon2}{Bellon,M.P. \cqg{23}{2006}{7029}.}
   \ref{McLellan}{McLellan,A,G. \jpc{7}{1974}{3326}.}
   \ref{Boiteaux}{Boiteaux, M. \jmp{23}{1982}{1311}.}
   \ref{HHandK}{Hage Hassan,M. and Kibler,M. {\it On Hurwitz
   transformations} in {Le probl\`eme de factorisation de Hurwitz}, Eds.,
   A.Ronveaux and D.Lambert (Fac.Univ.N.D. de la Paix, Namur, 1991),
   pp.1-29.}
   \ref{Weeks2}{Weeks,Jeffrey \cqg{23}{2006}{6971}.}
   \ref{LandW}{Lachi\`eze-Rey,M. and Weeks,Jeffrey, {\it Orbifold construction of
   the modes on the Poincar\'e dodecahedral space}, arXiv:0801.4232.}
   \ref{Cayley4}{Cayley,A. \qjpam{58}{1879}{280}.}
   \ref{JMS}{Jari\'c,M.V., Michel,L. and Sharp,R.T. {\it J.Physique}
   {\bf 45} (1984) 1. }
   \ref{AandB}{Altmann,S.L. and Bradley,C.J.  {\it Phil. Trans. Roy. Soc. Lond.}
   {\bf 255} (1963) 199.}
   \ref{CandP}{Cummins,C.J. and Patera,J. \jmp{29}{1988}{1736}.}
   \ref{Sloane}{Sloane,N.J.A. \amm{84}{1977}{82}.}
   \ref{Gordan2}{Gordan,P. \ma{12}{1877}{147}.}
   \ref{DandSh}{Desmier,P.E. and Sharp,R.T. \jmp{20}{1979}{74}.}
   \ref{Kramer}{Kramer,P., \jpa{38}{2005}{3517}.}
   \ref{Klein2}{Klein, F.\ma{9}{1875}{183}.}
   \ref{Hodgkinson}{Hodgkinson,J. \jlms{10}{1935}{221}.}
   \ref{ZandD}{Zheng,Y. and Doerschuk, P.C. {\it Acta Cryst.} {\bf A52}
   (1996) 221.}
   \ref{EPM}{Elcoro,L., Perez--Mato,J.M. and Madariaga,G.
   {\it Acta Cryst.} {\bf A50} (1994) 182.}
    \ref{PSW2}{Prandl,W., Schiebel,P. and Wulf,K.
   {\it Acta Cryst.} {\bf A52} (1999) 171.}
    \ref{FCD}{Fan,P--D., Chen,J--Q. and Draayer,J.P.
   {\it Acta Cryst.} {\bf A55} (1999) 871.}
   \ref{FCD2}{Fan,P--D., Chen,J--Q. and Draayer,J.P.
   {\it Acta Cryst.} {\bf A55} (1999) 1049.}
   \ref{Honl}{H\"onl,H. \zfp{89}{1934}{244}.}
   \ref{PSW}{Patera,J., Sharp,R.T. and Winternitz,P. \jmp{19}{1978}{2362}.}
   \ref{LandH}{Lohe,M.A. and Hurst,C.A. \jmp{12}{1971}{1882}.}
   \ref{RandSA}{Ronveaux,A. and Saint-Aubin,Y. \jmp{24}{1983}{1037}.}
   \ref{JandDeV}{Jonker,J.E. and De Vries,E. \npa{105}{1967}{621}.}
   \ref{Rowe}{Rowe, E.G.Peter. \jmp{19}{1978}{1962}.}
   \ref{KNR}{Kibler,M., N\'egadi,T. and Ronveaux,A. {\it The Kustaanheimo-Stiefel
   transformation and certain special functions} \lnm{1171}{1985}{497}.}
   \ref{GLP}{Gilkey,P.B., Leahy,J.V. and Park,J-H, \jpa{29}{1996}{5645}.}
   \ref{Kohler}{K\"ohler,K.: Equivariant Reidemeister torsion on
   symmetric spaces. Math.Ann. {\bf 307}, 57-69 (1997)}
   \ref{Kohler2}{K\"ohler,K.: Equivariant analytic torsion on ${\bf P^nC}$.
   Math.Ann.{\bf 297}, 553-565 (1993) }
   \ref{Kohler3}{K\"ohler,K.: Holomorphic analytic torsion on Hermitian
   symmetric spaces. J.Reine Angew.Math. {\bf 460}, 93-116 (1995)}
   \ref{Zagierzf}{Zagier,D. {\it Zetafunktionen und Quadratische
   K\"orper}, (Springer--Verlag, Berlin, 1981).}
   \ref{Stek}{Stekholschkik,R. {\it Notes on Coxeter transformations and the McKay
   correspondence.} (Springer, Berlin, 2008).}
   \ref{Pesce}{Pesce,H. \cmh {71}{1996}{243}.}
   \ref{Pesce2}{Pesce,H. {\it Contemp. Math} {\bf 173} (1994) 231.}
   \ref{Sutton}{Sutton,C.J. {\it Equivariant isospectrality
   and isospectral deformations on spherical orbifolds}, ArXiv:math/0608567.}
   \ref{Sunada}{Sunada,T. \aom{121}{1985}{169}.}
   \ref{GoandM}{Gornet,R, and McGowan,J. {\it J.Comp. and Math.}
   {\bf 9} (2006) 270.}
   \ref{Suter}{Suter,R. {\it Manusc.Math.} {\bf 122} (2007) 1-21.}
   \ref{Lomont}{Lomont,J.S. {\it Applications of finite groups} (Academic
   Press, New York, 1959).}
   \ref{DandCh2}{Dowker,J.S. and Chang,Peter {\it Analytic torsion on
   spherical factors and tessellations}, arXiv:math.DG/0904.0744 .}
   \ref{Mackey}{Mackey,G. {\it Induced representations}
   (Benjamin, New York, 1968).}
   \ref{Koca}{Koca, {\it Turkish J.Physics}.}
   \ref{Brylinski}{Brylinski, J-L., {\it A correspondence dual to McKay's}
    ArXiv alg-geom/9612003.}
   \ref{Rossmann}{Rossman,W. {\it McKay's correspondence
   and characters of finite subgroups of\break SU(2)} {\it Progress in Math.}
      Birkhauser  (to appear) .}
   \ref{JandL}{James, G. and Liebeck, M. {\it Representations and
   characters of groups} (CUP, Cambridge, 2001).}
   \ref{IandR}{Ito,Y. and Reid,M. {\it The Mckay correspondence for finite
   subgroups of SL(3,C)} Higher dimensional varieties, (Trento 1994),
   221-240, (Berlin, de Gruyter 1996).}
   \ref{BandF}{Bauer,W. and Furutani, K. {\it J.Geom. and Phys.} {\bf
   58} (2008) 64.}
   \ref{Luck}{L\"uck,W. \jdg{37}{1993}{263}.}
   \ref{LandR}{Lott,J. and Rothenberg,M. \jdg{34}{1991}{431}.}
   \ref{DoandKi} {Dowker.J.S. and Kirsten, K. {\it Analysis and Appl.}
   {\bf 3} (2005) 45.}
   \ref{dowtess1}{Dowker,J.S. \cqg{23}{2006}{1}.}
   \ref{dowtess2}{Dowker,J.S. {\it J.Geom. and Phys.} {\bf 57} (2007) 1505.}
   \ref{MHS}{De Melo,T., Hartmann,L. and Spreafico,M. {\it Reidemeister
   Torsion and analytic torsion of discs}, ArXiv:0811.3196.}
   \ref{Vertman}{Vertman, B. {\it Analytic Torsion of a  bounded
   generalized cone}, ArXiv:0808.0449.}
   \ref{WandY} {Weng,L. and You,Y., {\it Int.J. of Math.}{\bf 7} (1996)
   109.}
   \ref{ScandT}{Schwartz, A.S. and Tyupkin,Yu.S. \np{242}{1984}{436}.}
   \ref{AAR}{Andrews, G.E., Askey,R. and Roy,R. {\it Special functions}
  (CUP, Cambridge, 1999).}
   \ref{Tsuchiya}{Tsuchiya, N.: R-torsion and analytic torsion for spherical
   Clifford-Klein manifolds.: J. Fac.Sci., Tokyo Univ. Sect.1 A, Math.
   {\bf 23}, 289-295 (1976).}
   \ref{Tsuchiya2}{Tsuchiya, N. J. Fac.Sci., Tokyo Univ. Sect.1 A, Math.
   {\bf 23}, 289-295 (1976).}
  \ref{Lerch}{Lerch,M. \am{11}{1887}{19}.}
  \ref{Lerch2}{Lerch,M. \am{29}{1905}{333}.}
  \ref{TandS}{Threlfall, W. and Seifert, H. \ma{104}{1930}{1}.}
  \ref{RandS}{Ray, D.B., and Singer, I. \aim{7}{1971}{145}.}
  \ref{RandS2}{Ray, D.B., and Singer, I. {\it Proc.Symp.Pure Math.}
  {\bf 23} (1973) 167.}
  \ref{Jensen}{Jensen,J.L.W.V. \aom{17}{1915-1916}{124}.}
  \ref{Rosenberg}{Rosenberg, S. {\it The Laplacian on a Riemannian Manifold}
  (CUP, Cambridge, 1997).}
  \ref{Nando2}{Nash, C. and O'Connor, D-J. {\it Int.J.Mod.Phys.}
  {\bf A10} (1995) 1779.}
  \ref{Fock}{Fock,V. \zfp{98}{1935}{145}.}
  \ref{Levy}{Levy,M. \prs {204}{1950}{145}.}
  \ref{Schwinger2}{Schwinger,J. \jmp{5}{1964}{1606}.}
  \ref{Muller}{M\"uller, \lnm{}{}{}.}
  \ref{VMK}{Varshalovich.}
  \ref{DandWo}{Dowker,J.S. and Wolski, A. \prA{46}{1992}{6417}.}
  \ref{Zeitlin1}{Zeitlin,V. {\it Physica D} {\bf 49} (1991).  }
  \ref{Zeitlin0}{Zeitlin,V. {\it Nonlinear World} Ed by
   V.Baryakhtar {\it et al},  Vol.I p.717,  (World Scientific, Singapore, 1989).}
  \ref{Zeitlin2}{Zeitlin,V. \prl{93}{2004}{264501}. }
  \ref{Zeitlin3}{Zeitlin,V. \pla{339}{2005}{316}. }
  \ref{Groenewold}{Groenewold, H.J. {\it Physica} {\bf 12} (1946) 405.}
  \ref{Cohen}{Cohen, L. \jmp{7}{1966}{781}.}
  \ref{AandW}{Argawal G.S. and Wolf, E. \prD{2}{1970}{2161,2187,2206}.}
  \ref{Jantzen}{Jantzen,R.T. \jmp{19}{1978}{1163}.}
  \ref{Moses2}{Moses,H.E. \aop{42}{1967}{343}.}
  \ref{Carmeli}{Carmeli,M. \jmp{9}{1968}{1987}.}
  \ref{SHS}{Siemans,M., Hancock,J. and Siminovitch,D. {\it Solid State
  Nuclear Magnetic Resonance} {\bf 31}(2007)35.}
 \ref{Dowk}{Dowker,J.S. \prD{28}{1983}{3013}.}
 \ref{Heine}{Heine, E. {\it Handbuch der Kugelfunctionen}
  (G.Reimer, Berlin. 1878, 1881).}
  \ref{Pockels}{Pockels, F. {\it \"Uber die Differentialgleichung $\De
  u+k^2u=0$} (Teubner, Leipzig. 1891).}
  \ref{Hamermesh}{Hamermesh, M., {\it Group Theory} (Addison--Wesley,
  Reading. 1962).}
  \ref{Racah}{Racah, G. {\it Group Theory and Spectroscopy}
  (Princeton Lecture Notes, 1951). }
  \ref{Gourdin}{Gourdin, M. {\it Basics of Lie Groups} (Editions
  Fronti\'eres, Gif sur Yvette. 1982.)}
  \ref{Clifford}{Clifford, W.K. \plms{2}{1866}{116}.}
  \ref{Story2}{Story, W.E. \plms{23}{1892}{265}.}
  \ref{Story}{Story, W.E. \ma{41}{1893}{469}.}
  \ref{Poole}{Poole, E.G.C. \plms{33}{1932}{435}.}
  \ref{Dickson}{Dickson, L.E. {\it Algebraic Invariants} (Wiley, N.Y.
  1915).}
  \ref{Dickson2}{Dickson, L.E. {\it Modern Algebraic Theories}
  (Sanborn and Co., Boston. 1926).}
  \ref{Hilbert2}{Hilbert, D. {\it Theory of algebraic invariants} (C.U.P.,
  Cambridge. 1993).}
  \ref{Olver}{Olver, P.J. {\it Classical Invariant Theory} (C.U.P., Cambridge.
  1999.)}
  \ref{AST}{A\v{s}erova, R.M., Smirnov, J.F. and Tolsto\v{i}, V.N. {\it
  Teoret. Mat. Fyz.} {\bf 8} (1971) 255.}
  \ref{AandS}{A\v{s}erova, R.M., Smirnov, J.F. \np{4}{1968}{399}.}
  \ref{Shapiro}{Shapiro, J. \jmp{6}{1965}{1680}.}
  \ref{Shapiro2}{Shapiro, J.Y. \jmp{14}{1973}{1262}.}
  \ref{NandS}{Noz, M.E. and Shapiro, J.Y. \np{51}{1973}{309}.}
  \ref{Cayley2}{Cayley, A. {\it Phil. Trans. Roy. Soc. Lond.}
  {\bf 144} (1854) 244.}
  \ref{Cayley3}{Cayley, A. {\it Phil. Trans. Roy. Soc. Lond.}
  {\bf 146} (1856) 101.}
  \ref{Wigner}{Wigner, E.P. {\it Gruppentheorie} (Vieweg, Braunschweig. 1931).}
  \ref{Sharp}{Sharp, R.T. \ajop{28}{1960}{116}.}
  \ref{Laporte}{Laporte, O. {\it Z. f. Naturf.} {\bf 3a} (1948) 447.}
  \ref{Lowdin}{L\"owdin, P-O. \rmp{36}{1964}{966}.}
  \ref{Ansari}{Ansari, S.M.R. {\it Fort. d. Phys.} {\bf 15} (1967) 707.}
  \ref{SSJR}{Samal, P.K., Saha, R., Jain, P. and Ralston, J.P. {\it
  Testing Isotropy of Cosmic Microwave Background Radiation},
  astro-ph/0708.2816.}
  \ref{Lachieze}{Lachi\'eze-Rey, M. {\it Harmonic projection and
  multipole Vectors}. astro- \break ph/0409081.}
  \ref{CHS}{Copi, C.J., Huterer, D. and Starkman, G.D.
  \prD{70}{2003}{043515}.}
  \ref{Jaric}{Jari\'c, J.P. {\it Int. J. Eng. Sci.} {\bf 41} (2003) 2123.}
  \ref{RandD}{Roche, J.A. and Dowker, J.S. \jpa{1}{1968}{527}.}
  \ref{KandW}{Katz, G. and Weeks, J.R. \prD{70}{2004}{063527}.}
  \ref{Waerden}{van der Waerden, B.L. {\it Die Gruppen-theoretische
  Methode in der Quantenmechanik} (Springer, Berlin. 1932).}
  \ref{EMOT}{Erdelyi, A., Magnus, W., Oberhettinger, F. and Tricomi, F.G. {
  \it Higher Transcendental Functions} Vol.1 (McGraw-Hill, N.Y. 1953).}
  \ref{Dowzilch}{Dowker, J.S. {\it Proc. Phys. Soc.} {\bf 91} (1967) 28.}
  \ref{DandD}{Dowker, J.S. and Dowker, Y.P. {\it Proc. Phys. Soc.}
  {\bf 87} (1966) 65.}
  \ref{DandD2}{Dowker, J.S. and Dowker, Y.P. \prs{}{}{}.}
  \ref{Dowk3}{Dowker,J.S. \cqg{7}{1990}{1241}.}
  \ref{Dowk5}{Dowker,J.S. \cqg{7}{1990}{2353}.}
  \ref{CoandH}{Courant, R. and Hilbert, D. {\it Methoden der
  Mathematischen Physik} vol.1 \break (Springer, Berlin. 1931).}
  \ref{Applequist}{Applequist, J. \jpa{22}{1989}{4303}.}
  \ref{Torruella}{Torruella, \jmp{16}{1975}{1637}.}
  \ref{Weinberg}{Weinberg, S.W. \pr{133}{1964}{B1318}.}
  \ref{Meyerw}{Meyer, W.F. {\it Apolarit\"at und rationale Curven}
  (Fues, T\"ubingen. 1883.) }
  \ref{Ostrowski}{Ostrowski, A. {\it Jahrsb. Deutsch. Math. Verein.} {\bf
  33} (1923) 245.}
  \ref{Kramers}{Kramers, H.A. {\it Grundlagen der Quantenmechanik}, (Akad.
  Verlag., Leipzig, 1938).}
  \ref{ZandZ}{Zou, W.-N. and Zheng, Q.-S. \prs{459}{2003}{527}.}
  \ref{Weeks1}{Weeks, J.R. {\it Maxwell's multipole vectors
  and the CMB}.  astro-ph/0412231.}
  \ref{Corson}{Corson, E.M. {\it Tensors, Spinors and Relativistic Wave
  Equations} (Blackie, London. 1950).}
  \ref{Rosanes}{Rosanes, J. \jram{76}{1873}{312}.}
  \ref{Salmon}{Salmon, G. {\it Lessons Introductory to the Modern Higher
  Algebra} 3rd. edn. \break (Hodges,  Dublin. 1876.)}
  \ref{Milnew}{Milne, W.P. {\it Homogeneous Coordinates} (Arnold. London. 1910).}
  \ref{Niven}{Niven, W.D. {\it Phil. Trans. Roy. Soc.} {\bf 170} (1879) 393.}
  \ref{Scott}{Scott, C.A. {\it An Introductory Account of
  Certain Modern Ideas and Methods in Plane Analytical Geometry,}
  (MacMillan, N.Y. 1896).}
  \ref{Bargmann}{Bargmann, V. \rmp{34}{1962}{300}.}
  \ref{Maxwell}{Maxwell, J.C. {\it A Treatise on Electricity and
  Magnetism} 2nd. edn. (Clarendon Press, Oxford. 1882).}
  \ref{BandL}{Biedenharn, L.C. and Louck, J.D.
  {\it Angular Momentum in Quantum Physics} (Addison-Wesley, Reading. 1981).}
  \ref{Weylqm}{Weyl, H. {\it The Theory of Groups and Quantum Mechanics}
  (Methuen, London. 1931).}
  \ref{Robson}{Robson, A. {\it An Introduction to Analytical Geometry} Vol I
  (C.U.P., Cambridge. 1940.)}
  \ref{Sommerville}{Sommerville, D.M.Y. {\it Analytical Conics} 3rd. edn.
   (Bell, London. 1933).}
  \ref{Coolidge}{Coolidge, J.L. {\it A Treatise on Algebraic Plane Curves}
  (Clarendon Press, Oxford. 1931).}
  \ref{SandK}{Semple, G. and Kneebone. G.T. {\it Algebraic Projective
  Geometry} (Clarendon Press, Oxford. 1952).}
  \ref{AandC}{Abdesselam A., and Chipalkatti, J. {\it The Higher
  Transvectants are redundant}, arXiv:0801.1533 [math.AG] 2008.}
  \ref{Elliott}{Elliott, E.B. {\it The Algebra of Quantics} 2nd edn.
  (Clarendon Press, Oxford. 1913).}
  \ref{Elliott2}{Elliott, E.B. \qjpam{48}{1917}{372}.}
  \ref{Howe}{Howe, R. \tams{313}{1989}{539}.}
  \ref{Clebsch}{Clebsch, A. \jram{60}{1862}{343}.}
  \ref{Prasad}{Prasad, G. \ma{72}{1912}{136}.}
  \ref{Dougall}{Dougall, J. \pems{32}{1913}{30}.}
  \ref{Penrose}{Penrose, R. \aop{10}{1960}{171}.}
  \ref{Penrose2}{Penrose, R. \prs{273}{1965}{171}.}
  \ref{Burnside}{Burnside, W.S. \qjm{10}{1870}{211}. }
  \ref{Lindemann}{Lindemann, F. \ma{23} {1884}{111}.}
  \ref{Backus}{Backus, G. {\it Rev. Geophys. Space Phys.} {\bf 8} (1970) 633.}
  \ref{Baerheim}{Baerheim, R. {\it Q.J. Mech. appl. Math.} {\bf 51} (1998) 73.}
  \ref{Lense}{Lense, J. {\it Kugelfunktionen} (Akad.Verlag, Leipzig. 1950).}
  \ref{Littlewood}{Littlewood, D.E. \plms{50}{1948}{349}.}
  \ref{Fierz}{Fierz, M. {\it Helv. Phys. Acta} {\bf 12} (1938) 3.}
  \ref{Williams}{Williams, D.N. {\it Lectures in Theoretical Physics} Vol. VII,
  (Univ.Colorado Press, Boulder. 1965).}
  \ref{Dennis}{Dennis, M. \jpa{37}{2004}{9487}.}
  \ref{Pirani}{Pirani, F. {\it Brandeis Lecture Notes on
  General Relativity,} edited by S. Deser and K. Ford. (Brandeis, Mass. 1964).}
  \ref{Sturm}{Sturm, R. \jram{86}{1878}{116}.}
  \ref{Schlesinger}{Schlesinger, O. \ma{22}{1883}{521}.}
  \ref{Askwith}{Askwith, E.H. {\it Analytical Geometry of the Conic
  Sections} (A.\&C. Black, London. 1908).}
  \ref{Todd}{Todd, J.A. {\it Projective and Analytical Geometry}.
  (Pitman, London. 1946).}
  \ref{Glenn}{Glenn. O.E. {\it Theory of Invariants} (Ginn \& Co, N.Y. 1915).}
  \ref{DowkandG}{Dowker, J.S. and Goldstone, M. \prs{303}{1968}{381}.}
  \ref{Turnbull}{Turnbull, H.A. {\it The Theory of Determinants,
  Matrices and Invariants} 3rd. edn. (Dover, N.Y. 1960).}
  \ref{MacMillan}{MacMillan, W.D. {\it The Theory of the Potential}
  (McGraw-Hill, N.Y. 1930).}
   \ref{Hobson}{Hobson, E.W. {\it The Theory of Spherical
   and Ellipsoidal Harmonics} (C.U.P., Cambridge. 1931).}
  \ref{Hobson1}{Hobson, E.W. \plms {24}{1892}{55}.}
  \ref{GandY}{Grace, J.H. and Young, A. {\it The Algebra of Invariants}
  (C.U.P., Cambridge, 1903).}
  \ref{FandR}{Fano, U. and Racah, G. {\it Irreducible Tensorial Sets}
  (Academic Press, N.Y. 1959).}
  \ref{TandT}{Thomson, W. and Tait, P.G. {\it Treatise on Natural Philosophy}
   (Clarendon Press, Oxford. 1867).}
  \ref{Brinkman}{Brinkman, H.C. {\it Applications of spinor invariants in
atomic physics}, North Holland, Amsterdam 1956.}
  \ref{Kramers1}{Kramers, H.A. {\it Proc. Roy. Soc. Amst.} {\bf 33} (1930) 953.}
  \ref{DandP2}{Dowker,J.S. and Pettengill,D.F. \jpa{7}{1974}{1527}}
  \ref{Dowk1}{Dowker,J.S. \jpa{}{}{45}.}
  \ref{Dowk2}{Dowker,J.S. \aop{71}{1972}{577}}
  \ref{DandA}{Dowker,J.S. and Apps, J.S. \cqg{15}{1998}{1121}.}
  \ref{Weil}{Weil,A., {\it Elliptic functions according to Eisenstein
  and Kronecker}, Springer, Berlin, 1976.}
  \ref{Ling}{Ling,C-H. {\it SIAM J.Math.Anal.} {\bf5} (1974) 551.}
  \ref{Ling2}{Ling,C-H. {\it J.Math.Anal.Appl.}(1988).}
 \ref{BMO}{Brevik,I., Milton,K.A. and Odintsov, S.D. \aop{302}{2002}{120}.}
 \ref{KandL}{Kutasov,D. and Larsen,F. {\it JHEP} 0101 (2001) 1.}
 \ref{KPS}{Klemm,D., Petkou,A.C. and Siopsis {\it Entropy
 bounds, monoticity properties and scaling in CFT's}. hep-th/0101076.}
 \ref{DandC}{Dowker,J.S. and Critchley,R. \prD{15}{1976}{1484}.}
 \ref{AandD}{Al'taie, M.B. and Dowker, J.S. \prD{18}{1978}{3557}.}
 \ref{Dow1}{Dowker,J.S. \prD{37}{1988}{558}.}
 \ref{Dow30}{Dowker,J.S. \prD{28}{1983}{3013}.}
 \ref{DandK}{Dowker,J.S. and Kennedy,G. \jpa{}{1978}{895}.}
 \ref{Dow2}{Dowker,J.S. \cqg{1}{1984}{359}.}
 \ref{DandKi}{Dowker,J.S. and Kirsten, K. {\it Comm. in Anal. and Geom.
 }{\bf7} (1999) 641.}
 \ref{DandKe}{Dowker,J.S. and Kennedy,G.\jpa{11}{1978}{895}.}
 \ref{Gibbons}{Gibbons,G.W. \pl{60A}{1977}{385}.}
 \ref{Cardy}{Cardy,J.L. \np{366}{1991}{403}.}
 \ref{ChandD}{Chang,P. and Dowker,J.S. \np{395}{1993}{407}.}
 \ref{DandC2}{Dowker,J.S. and Critchley,R. \prD{13}{1976}{224}.}
 \ref{Camporesi}{Camporesi,R. \prp{196}{1990}{1}.}
 \ref{BandM}{Brown,L.S. and Maclay,G.J. \pr{184}{1969}{1272}.}
 \ref{CandD}{Candelas,P. and Dowker,J.S. \prD{19}{1979}{2902}.}
 \ref{Unwin1}{Unwin,S.D. Thesis. University of Manchester. 1979.}
 \ref{Unwin2}{Unwin,S.D. \jpa{13}{1980}{313}.}
 \ref{DandB}{Dowker,J.S. and Banach,R. \jpa{11}{1978}{2255}.}
 \ref{Obhukov}{Obhukov,Yu.N. \pl{109B}{1982}{195}.}
 \ref{Kennedy}{Kennedy,G. \prD{23}{1981}{2884}.}
 \ref{CandT}{Copeland,E. and Toms,D.J. \np {255}{1985}{201}.}
  \ref{CandT2}{Copeland,E. and Toms,D.J. \cqg {3}{1986}{431}.}
 \ref{ELV}{Elizalde,E., Lygren, M. and Vassilevich,
 D.V. \jmp {37}{1996}{3105}.}
 \ref{Malurkar}{Malurkar,S.L. {\it J.Ind.Math.Soc} {\bf16} (1925/26) 130.}
 \ref{Glaisher}{Glaisher,J.W.L. {\it Messenger of Math.} {\bf18} (1889) 1.}
  \ref{Anderson}{Anderson,A. \prD{37}{1988}{536}.}
 \ref{CandA}{Cappelli,A. and D'Appollonio,\pl{487B}{2000}{87}.}
 \ref{Wot}{Wotzasek,C. \jpa{23}{1990}{1627}.}
 \ref{RandT}{Ravndal,F. and Tollesen,D. \prD{40}{1989}{4191}.}
 \ref{SandT}{Santos,F.C. and Tort,A.C. \pl{482B}{2000}{323}.}
 \ref{FandO}{Fukushima,K. and Ohta,K. {\it Physica} {\bf A299} (2001) 455.}
 \ref{GandP}{Gibbons,G.W. and Perry,M. \prs{358}{1978}{467}.}
 \ref{Dow4}{Dowker,J.S..}
  \ref{Rad}{Rademacher,H. {\it Topics in analytic number theory,}
Springer-Verlag,  Berlin,1973.}
  \ref{Halphen}{Halphen,G.-H. {\it Trait\'e des Fonctions Elliptiques},
  Vol 1, Gauthier-Villars, Paris, 1886.}
  \ref{CandW}{Cahn,R.S. and Wolf,J.A. {\it Comm.Mat.Helv.} {\bf 51}
  (1976) 1.}
  \ref{Berndt}{Berndt,B.C. \rmjm{7}{1977}{147}.}
  \ref{Hurwitz}{Hurwitz,A. \ma{18}{1881}{528}.}
  \ref{Hurwitz2}{Hurwitz,A. {\it Mathematische Werke} Vol.I. Basel,
  Birkhauser, 1932.}
  \ref{Berndt2}{Berndt,B.C. \jram{303/304}{1978}{332}.}
  \ref{RandA}{Rao,M.B. and Ayyar,M.V. \jims{15}{1923/24}{150}.}
  \ref{Hardy}{Hardy,G.H. \jlms{3}{1928}{238}.}
  \ref{TandM}{Tannery,J. and Molk,J. {\it Fonctions Elliptiques},
   Gauthier-Villars, Paris, 1893--1902.}
  \ref{schwarz}{Schwarz,H.-A. {\it Formeln und
  Lehrs\"atzen zum Gebrauche..},Springer 1893.(The first edition was 1885.)
  The French translation by Henri Pad\'e is {\it Formules et Propositions
  pour L'Emploi...},Gauthier-Villars, Paris, 1894}
  \ref{Hancock}{Hancock,H. {\it Theory of elliptic functions}, Vol I.
   Wiley, New York 1910.}
  \ref{watson}{Watson,G.N. \jlms{3}{1928}{216}.}
  \ref{MandO}{Magnus,W. and Oberhettinger,F. {\it Formeln und S\"atze},
  Springer-Verlag, Berlin 1948.}
  \ref{Klein}{Klein,F. {\it Lectures on the Icosohedron}
  (Methuen, London. 1913).}
  \ref{AandL}{Appell,P. and Lacour,E. {\it Fonctions Elliptiques},
  Gauthier-Villars,
  Paris. 1897.}
  \ref{HandC}{Hurwitz,A. and Courant,C. {\it Allgemeine Funktionentheorie},
  Springer,
  Berlin. 1922.}
  \ref{WandW}{Whittaker,E.T. and Watson,G.N. {\it Modern analysis},
  Cambridge. 1927.}
  \ref{SandC}{Selberg,A. and Chowla,S. \jram{227}{1967}{86}. }
  \ref{zucker}{Zucker,I.J. {\it Math.Proc.Camb.Phil.Soc} {\bf 82 }(1977)
  111.}
  \ref{glasser}{Glasser,M.L. {\it Maths.of Comp.} {\bf 25} (1971) 533.}
  \ref{GandW}{Glasser, M.L. and Wood,V.E. {\it Maths of Comp.} {\bf 25}
  (1971)
  535.}
  \ref{greenhill}{Greenhill,A,G. {\it The Applications of Elliptic
  Functions}, MacMillan. London, 1892.}
  \ref{Weierstrass}{Weierstrass,K. {\it J.f.Mathematik (Crelle)}
{\bf 52} (1856) 346.}
  \ref{Weierstrass2}{Weierstrass,K. {\it Mathematische Werke} Vol.I,p.1,
  Mayer u. M\"uller, Berlin, 1894.}
  \ref{Fricke}{Fricke,R. {\it Die Elliptische Funktionen und Ihre Anwendungen},
    Teubner, Leipzig. 1915, 1922.}
  \ref{Konig}{K\"onigsberger,L. {\it Vorlesungen \"uber die Theorie der
 Elliptischen Funktionen},  \break Teubner, Leipzig, 1874.}
  \ref{Milne}{Milne,S.C. {\it The Ramanujan Journal} {\bf 6} (2002) 7-149.}
  \ref{Schlomilch}{Schl\"omilch,O. {\it Ber. Verh. K. Sachs. Gesell. Wiss.
  Leipzig}  {\bf 29} (1877) 101-105; {\it Compendium der h\"oheren
  Analysis}, Bd.II, 3rd Edn, Vieweg, Brunswick, 1878.}
  \ref{BandB}{Briot,C. and Bouquet,C. {\it Th\`eorie des Fonctions
  Elliptiques}, Gauthier-Villars, Paris, 1875.}
  \ref{Dumont}{Dumont,D. \aim {41}{1981}{1}.}
  \ref{Andre}{Andr\'e,D. {\it Ann.\'Ecole Normale Superior} {\bf 6} (1877)
  265;
  {\it J.Math.Pures et Appl.} {\bf 5} (1878) 31.}
  \ref{Raman}{Ramanujan,S. {\it Trans.Camb.Phil.Soc.} {\bf 22} (1916) 159;
 {\it Collected Papers}, Cambridge, 1927}
  \ref{Weber}{Weber,H.M. {\it Lehrbuch der Algebra} Bd.III, Vieweg,
  Brunswick 190  3.}
  \ref{Weber2}{Weber,H.M. {\it Elliptische Funktionen und algebraische
  Zahlen},
  Vieweg, Brunswick 1891.}
  \ref{ZandR}{Zucker,I.J. and Robertson,M.M.
  {\it Math.Proc.Camb.Phil.Soc} {\bf 95 }(1984) 5.}
  \ref{JandZ1}{Joyce,G.S. and Zucker,I.J.
  {\it Math.Proc.Camb.Phil.Soc} {\bf 109 }(1991) 257.}
  \ref{JandZ2}{Zucker,I.J. and Joyce.G.S.
  {\it Math.Proc.Camb.Phil.Soc} {\bf 131 }(2001) 309.}
  \ref{zucker2}{Zucker,I.J. {\it SIAM J.Math.Anal.} {\bf 10} (1979) 192,}
  \ref{BandZ}{Borwein,J.M. and Zucker,I.J. {\it IMA J.Math.Anal.} {\bf 12}
  (1992) 519.}
  \ref{Cox}{Cox,D.A. {\it Primes of the form $x^2+n\,y^2$}, Wiley,
  New York, 1989.}
  \ref{BandCh}{Berndt,B.C. and Chan,H.H. {\it Mathematika} {\bf42} (1995)
  278.}
  \ref{EandT}{Elizalde,R. and Tort.hep-th/}
  \ref{KandS}{Kiyek,K. and Schmidt,H. {\it Arch.Math.} {\bf 18} (1967) 438.}
  \ref{Oshima}{Oshima,K. \prD{46}{1992}{4765}.}
  \ref{greenhill2}{Greenhill,A.G. \plms{19} {1888} {301}.}
  \ref{Russell}{Russell,R. \plms{19} {1888} {91}.}
  \ref{BandB}{Borwein,J.M. and Borwein,P.B. {\it Pi and the AGM}, Wiley,
  New York, 1998.}
  \ref{Resnikoff}{Resnikoff,H.L. \tams{124}{1966}{334}.}
  \ref{vandp}{Van der Pol, B. {\it Indag.Math.} {\bf18} (1951) 261,272.}
  \ref{Rankin}{Rankin,R.A. {\it Modular forms} C.U.P. Cambridge}
  \ref{Rankin2}{Rankin,R.A. {\it Proc. Roy.Soc. Edin.} {\bf76 A} (1976) 107.}
  \ref{Skoruppa}{Skoruppa,N-P. {\it J.of Number Th.} {\bf43} (1993) 68 .}
  \ref{Down}{Dowker.J.S. {\it Nucl.Phys.}B (Proc.Suppl) ({\bf 104})(2002)153;
  also Dowker,J.S. hep-th/ 0007129.}
  \ref{Eichler}{Eichler,M. \mz {67}{1957}{267}.}
  \ref{Zagier}{Zagier,D. \invm{104}{1991}{449}.}
  \ref{Lang}{Lang,S. {\it Modular Forms}, Springer, Berlin, 1976.}
  \ref{Kosh}{Koshliakov,N.S. {\it Mess.of Math.} {\bf 58} (1928) 1.}
  \ref{BandH}{Bodendiek, R. and Halbritter,U. \amsh{38}{1972}{147}.}
  \ref{Smart}{Smart,L.R., \pgma{14}{1973}{1}.}
  \ref{Grosswald}{Grosswald,E. {\it Acta. Arith.} {\bf 21} (1972) 25.}
  \ref{Kata}{Katayama,K. {\it Acta Arith.} {\bf 22} (1973) 149.}
  \ref{Ogg}{Ogg,A. {\it Modular forms and Dirichlet series} (Benjamin,
  New York,
   1969).}
  \ref{Bol}{Bol,G. \amsh{16}{1949}{1}.}
  \ref{Epstein}{Epstein,P. \ma{56}{1903}{615}.}
  \ref{Petersson}{Petersson.}
  \ref{Serre}{Serre,J-P. {\it A Course in Arithmetic}, Springer,
  New York, 1973.}
  \ref{Schoenberg}{Schoenberg,B., {\it Elliptic Modular Functions},
  Springer, Berlin, 1974.}
  \ref{Apostol}{Apostol,T.M. \dmj {17}{1950}{147}.}
  \ref{Ogg2}{Ogg,A. {\it Lecture Notes in Math.} {\bf 320} (1973) 1.}
  \ref{Knopp}{Knopp,M.I. \dmj {45}{1978}{47}.}
  \ref{Knopp2}{Knopp,M.I. \invm {}{1994}{361}.}
  \ref{LandZ}{Lewis,J. and Zagier,D. \aom{153}{2001}{191}.}
  \ref{DandK1}{Dowker,J.S. and Kirsten,K. {\it Elliptic functions and
  temperature inversion symmetry on spheres} hep-th/.}
  \ref{HandK}{Husseini and Knopp.}
  \ref{Kober}{Kober,H. \mz{39}{1934-5}{609}.}
  \ref{HandL}{Hardy,G.H. and Littlewood, \am{41}{1917}{119}.}
  \ref{Watson}{Watson,G.N. \qjm{2}{1931}{300}.}
  \ref{SandC2}{Chowla,S. and Selberg,A. {\it Proc.Nat.Acad.} {\bf 35}
  (1949) 371.}
  \ref{Landau}{Landau, E. {\it Lehre von der Verteilung der Primzahlen},
  (Teubner, Leipzig, 1909).}
  \ref{Berndt4}{Berndt,B.C. \tams {146}{1969}{323}.}
  \ref{Berndt3}{Berndt,B.C. \tams {}{}{}.}
  \ref{Bochner}{Bochner,S. \aom{53}{1951}{332}.}
  \ref{Weil2}{Weil,A.\ma{168}{1967}{}.}
  \ref{CandN}{Chandrasekharan,K. and Narasimhan,R. \aom{74}{1961}{1}.}
  \ref{Rankin3}{Rankin,R.A. {} {} ().}
  \ref{Berndt6}{Berndt,B.C. {\it Trans.Edin.Math.Soc}.}
  \ref{Elizalde}{Elizalde,E. {\it Ten Physical Applications of Spectral
  Zeta Function Theory}, \break (Springer, Berlin, 1995).}
  \ref{Allen}{Allen,B., Folacci,A. and Gibbons,G.W. \pl{189}{1987}{304}.}
  \ref{Krazer}{Krazer}
  \ref{Elizalde3}{Elizalde,E. {\it J.Comp.and Appl. Math.} {\bf 118}
  (2000) 125.}
  \ref{Elizalde2}{Elizalde,E., Odintsov.S.D, Romeo, A. and Bytsenko,
  A.A and
  Zerbini,S.
  {\it Zeta function regularisation}, (World Scientific, Singapore,
  1994).}
  \ref{Eisenstein}{Eisenstein}
  \ref{Hecke}{Hecke,E. \ma{112}{1936}{664}.}
  \ref{Hecke2}{Hecke,E. \ma{112}{1918}{398}.}
  \ref{Terras}{Terras,A. {\it Harmonic analysis on Symmetric Spaces} (Springer,
  New York, 1985).}
  \ref{BandG}{Bateman,P.T. and Grosswald,E. {\it Acta Arith.} {\bf 9}
  (1964) 365.}
  \ref{Deuring}{Deuring,M. \aom{38}{1937}{585}.}
  \ref{Mordell}{Mordell,J. \prs{}{}{}.}
  \ref{GandZ}{Glasser,M.L. and Zucker, {}.}
  \ref{Landau2}{Landau,E. \jram{}{1903}{64}.}
  \ref{Kirsten1}{Kirsten,K. \jmp{35}{1994}{459}.}
  \ref{Sommer}{Sommer,J. {\it Vorlesungen \"uber Zahlentheorie}
  (1907,Teubner,Leipzig).
  French edition 1913 .}
  \ref{Reid}{Reid,L.W. {\it Theory of Algebraic Numbers},
  (1910,MacMillan,New York).}
  \ref{Milnor}{Milnor, J. {\it Is the Universe simply--connected?},
  IAS, Princeton, 1978.}
  \ref{Milnor2}{Milnor, J. \ajm{79}{1957}{623}.}
  \ref{Opechowski}{Opechowski,W. {\it Physica} {\bf 7} (1940) 552.}
  \ref{Bethe}{Bethe, H.A. \zfp{3}{1929}{133}.}
  \ref{LandL}{Landau, L.D. and Lishitz, E.M. {\it Quantum
  Mechanics} (Pergamon Press, London, 1958).}
  \ref{GPR}{Gibbons, G.W., Pope, C. and R\"omer, H., \np{157}{1979}{377}.}
  \ref{Jadhav}{Jadhav,S.P. PhD Thesis, University of Manchester 1990.}
  \ref{DandJ}{Dowker,J.S. and Jadhav, S. \prD{39}{1989}{1196}.}
  \ref{CandM}{Coxeter, H.S.M. and Moser, W.O.J. {\it Generators and
  relations of finite groups} (Springer. Berlin. 1957).}
  \ref{Coxeter2}{Coxeter, H.S.M. {\it Regular Complex Polytopes},
   (Cambridge University Press, \break Cambridge, 1975).}
  \ref{Coxeter}{Coxeter, H.S.M. {\it Regular Polytopes}.}
  \ref{Stiefel}{Stiefel, E., J.Research NBS {\bf 48} (1952) 424.}
  \ref{BandS}{Brink, D.M. and Satchler, G.R. {\it Angular momentum theory}.
  (Clarendon Press, Oxford. 1962.).}
  \ref{Rose}{Rose}
  \ref{Schwinger}{Schwinger, J. {\it On Angular Momentum}
  in {\it Quantum Theory of Angular Momentum} edited by
  Biedenharn,L.C. and van Dam, H. (Academic Press, N.Y. 1965).}
  
  \ref{Ray}{Ray,D.B. \aim{4}{1970}{109}.}
  \ref{Ikeda}{Ikeda,A. {\it Kodai Math.J.} {\bf 18} (1995) 57.}
  \ref{Kennedy}{Kennedy,G. \prD{23}{1981}{2884}.}
  \ref{Ellis}{Ellis,G.F.R. {\it General Relativity} {\bf2} (1971) 7.}
  \ref{Dow8}{Dowker,J.S. \cqg{20}{2003}{L105}.}
  \ref{IandY}{Ikeda, A and Yamamoto, Y. \ojm {16}{1979}{447}.}
  \ref{BandI}{Bander,M. and Itzykson,C. \rmp{18}{1966}{2}.}
  \ref{Schulman}{Schulman, L.S. \pr{176}{1968}{1558}.}
  \ref{Bar1}{B\"ar,C. {\it Arch.d.Math.}{\bf 59} (1992) 65.}
  \ref{Bar2}{B\"ar,C. {\it Geom. and Func. Anal.} {\bf 6} (1996) 899.}
  \ref{Vilenkin}{Vilenkin, N.J. {\it Special functions},
  (Am.Math.Soc., Providence, 1968).}
  \ref{Talman}{Talman, J.D. {\it Special functions} (Benjamin,N.Y.,1968).}
  \ref{Miller}{Miller, W. {\it Symmetry groups and their applications}
  (Wiley, N.Y., 1972).}
  \ref{Dow3}{Dowker,J.S. \cmp{162}{1994}{633}.}
  \ref{Cheeger}{Cheeger, J. \jdg {18}{1983}{575}.}
  \ref{Cheeger2}{Cheeger, J. \aom {109}{1979}{259}.}
  \ref{Dow6}{Dowker,J.S. \jmp{30}{1989}{770}.}
  \ref{Dow20}{Dowker,J.S. \jmp{35}{1994}{6076}.}
  \ref{Dowjmp}{Dowker,J.S. \jmp{35}{1994}{4989}.}
  \ref{Dow21}{Dowker,J.S. {\it Heat kernels and polytopes} in {\it
   Heat Kernel Techniques and Quantum Gravity}, ed. by S.A.Fulling,
   Discourses in Mathematics and its Applications, No.4, Dept.
   Maths., Texas A\&M University, College Station, Texas, 1995.}
  \ref{Dow9}{Dowker,J.S. \jmp{42}{2001}{1501}.}
  \ref{Dow7}{Dowker,J.S. \jpa{25}{1992}{2641}.}
  \ref{Warner}{Warner.N.P. \prs{383}{1982}{379}.}
  \ref{Wolf}{Wolf, J.A. {\it Spaces of constant curvature},
  (McGraw--Hill,N.Y., 1967).}
  \ref{Meyer}{Meyer,B. \cjm{6}{1954}{135}.}
  \ref{BandB}{B\'erard,P. and Besson,G. {\it Ann. Inst. Four.} {\bf 30}
  (1980) 237.}
  \ref{PandM}{P\'{o}lya,G. and Meyer,B. \cras{228}{1948}{28}.}
  \ref{Springer}{Springer, T.A. Lecture Notes in Math. vol 585 (Springer,
  Berlin,1977).}
  \ref{SeandT}{Threlfall, H. and Seifert, W. \ma{104}{1930}{1}.}
  \ref{Hopf}{Hopf,H. \ma{95}{1925}{313}. }
  \ref{Dow}{Dowker,J.S. \jpa{5}{1972}{936}.}
  \ref{LLL}{Lehoucq,R., Lachi\'eze-Rey,M. and Luminet, J.--P. {\it
  Astron.Astrophys.} {\bf 313} (1996) 339.}
  \ref{LaandL}{Lachi\'eze-Rey,M. and Luminet, J.--P.
  \prp{254}{1995}{135}.}
  \ref{Schwarzschild}{Schwarzschild, K., {\it Vierteljahrschrift der
  Ast.Ges.} {\bf 35} (1900) 337.}
  \ref{Starkman}{Starkman,G.D. \cqg{15}{1998}{2529}.}
  \ref{LWUGL}{Lehoucq,R., Weeks,J.R., Uzan,J.P., Gausman, E. and
  Luminet, J.--P. \cqg{19}{2002}{4683}.}
  \ref{Dow10}{Dowker,J.S. \prD{28}{1983}{3013}.}
  \ref{BandD}{Banach, R. and Dowker, J.S. \jpa{12}{1979}{2527}.}
  \ref{Jadhav2}{Jadhav,S. \prD{43}{1991}{2656}.}
  \ref{Gilkey}{Gilkey,P.B. {\it Invariance theory,the heat equation and
  the Atiyah--Singer Index theorem} (CRC Press, Boca Raton, 1994).}
  \ref{BandY}{Berndt,B.C. and Yeap,B.P. {\it Adv. Appl. Math.}
  {\bf29} (2002) 358.}
  \ref{HandR}{Hanson,A.J. and R\"omer,H. \pl{80B}{1978}{58}.}
  \ref{Hill}{Hill,M.J.M. {\it Trans.Camb.Phil.Soc.} {\bf 13} (1883) 36.}
  \ref{Cayley}{Cayley,A. {\it Quart.Math.J.} {\bf 7} (1866) 304.}
  \ref{Seade}{Seade,J.A. {\it Anal.Inst.Mat.Univ.Nac.Aut\'on
  M\'exico} {\bf 21} (1981) 129.}
  \ref{CM}{Cisneros--Molina,J.L. {\it Geom.Dedicata} {\bf84} (2001)
  \ref{Goette1}{Goette,S. \jram {526} {2000} 181.}
  207.}
  \ref{NandO}{Nash,C. and O'Connor,D--J, \jmp {36}{1995}{1462}.}
  \ref{Dows}{Dowker,J.S. \aop{71}{1972}{577}; Dowker,J.S. and Pettengill,D.F.
  \jpa{7}{1974}{1527}; J.S.Dowker in {\it Quantum Gravity}, edited by
  S. C. Christensen (Hilger,Bristol,1984)}
  \ref{Jadhav2}{Jadhav,S.P. \prD{43}{1991}{2656}.}
  \ref{Dow11}{Dowker,J.S. \cqg{21}{2004}4247.}
  \ref{Dow12}{Dowker,J.S. \cqg{21}{2004}4977.}
  \ref{Dow13}{Dowker,J.S. \jpa{38}{2005}1049.}
  \ref{Zagier}{Zagier,D. \ma{202}{1973}{149}}
  \ref{RandG}{Rademacher, H. and Grosswald,E. {\it Dedekind Sums},
  (Carus, MAA, 1972).}
  \ref{Berndt7}{Berndt,B, \aim{23}{1977}{285}.}
  \ref{HKMM}{Harvey,J.A., Kutasov,D., Martinec,E.J. and Moore,G.
  {\it Localised Tachyons and RG Flows}, hep-th/0111154.}
  \ref{Beck}{Beck,M., {\it Dedekind Cotangent Sums}, {\it Acta Arithmetica}
  {\bf 109} (2003) 109-139 ; math.NT/0112077.}
  \ref{McInnes}{McInnes,B. {\it APS instability and the topology of the brane
  world}, hep-th/0401035.}
  \ref{BHS}{Brevik,I, Herikstad,R. and Skriudalen,S. {\it Entropy Bound for the
  TM Electromagnetic Field in the Half Einstein Universe}; hep-th/0508123.}
  \ref{BandO}{Brevik,I. and Owe,C.  \prD{55}{4689}{1997}.}
  \ref{Kenn}{Kennedy,G. Thesis. University of Manchester 1978.}
  \ref{KandU}{Kennedy,G. and Unwin S. \jpa{12}{L253}{1980}.}
  \ref{BandO1}{Bayin,S.S.and Ozcan,M.
  \prD{48}{2806}{1993}; \prD{49}{5313}{1994}.}
  \ref{Chang}{Chang, P., {\it Quantum Field Theory on Regular Polytopes}.
   Thesis. University of Manchester, 1993.}
  \ref{Barnesa}{Barnes,E.W. {\it Trans. Camb. Phil. Soc.} {\bf 19} (1903) 374.}
  \ref{Barnesb}{Barnes,E.W. {\it Trans. Camb. Phil. Soc.}
  {\bf 19} (1903) 426.}
  \ref{Stanley1}{Stanley,R.P. \joa {49Hilf}{1977}{134}.}
  \ref{Stanley2}{Stanley,R.P. {\it Enumerative Combinatorics} vols.1,2
  (C.U.P., Cambridge, 1997, 1999.}
  \ref{Stanley}{Stanley,R.P. \bams {1}{1979}{475}.}
  \ref{Hurley}{Hurley,A.C. \pcps {47}{1951}{51}.}
  \ref{IandK}{Iwasaki,I. and Katase,K. {\it Proc.Japan Acad. Ser} {\bf A55}
  (1979) 141.}
  \ref{IandT}{Ikeda,A. and Taniguchi,Y. {\it Osaka J. Math.} {\bf 15} (1978)
  515.}
  \ref{GandM}{Gallot,S. and Meyer,D. \jmpa{54}{1975}{259}.}
  \ref{Flatto}{Flatto,L. {\it Enseign. Math.} {\bf 24} (1978) 237.}
  \ref{OandT}{Orlik,P and Terao,H. {\it Arrangements of Hyperplanes},
  Grundlehren der Math. Wiss. {\bf 300}, (Springer--Verlag, 1992).}
  \ref{Shepler}{Shepler,A.V. \joa{220}{1999}{314}.}
  \ref{SandT}{Solomon,L. and Terao,H. \cmh {73}{1998}{237}.}
  \ref{Vass}{Vassilevich, D.V. \plb {348}{1995}39.}
  \ref{Vass2}{Vassilevich, D.V. \jmp {36}{1995}3174.}
  \ref{CandH}{Camporesi,R. and Higuchi,A. {\it J.Geom. and Physics}
  {\bf 15} (1994) 57.}
  \ref{Solomon2}{Solomon,L. \tams{113}{1964}{274}.}
  \ref{Solomon}{Solomon,L. {\it Nagoya Math. J.} {\bf 22} (1963) 57.}
  \ref{Obukhov}{Obukhov,Yu.N. \pl{109B}{1982}{195}.}
  \ref{BGH}{Bernasconi,F., Graf,G.M. and Hasler,D. {\it The heat kernel
  expansion for the electromagnetic field in a cavity}; math-ph/0302035.}
  \ref{Baltes}{Baltes,H.P. \prA {6}{1972}{2252}.}
  \ref{BaandH}{Baltes.H.P and Hilf,E.R. {\it Spectra of Finite Systems}
  (Bibliographisches Institut, Mannheim, 1976).}
  \ref{Ray}{Ray,D.B. \aim{4}{1970}{109}.}
  \ref{Hirzebruch}{Hirzebruch,F. {\it Topological methods in algebraic
  geometry} (Springer-- Verlag,\break  Berlin, 1978). }
  \ref{BBG}{Bla\v{z}i\'c,N., Bokan,N. and Gilkey, P.B. {\it Ind.J.Pure and
  Appl.Math.} {\bf 23} (1992) 103.}
  \ref{WandWi}{Weck,N. and Witsch,K.J. {\it Math.Meth.Appl.Sci.} {\bf 17}
  (1994) 1017.}
  \ref{Norlund}{N\"orlund,N.E. \am{43}{1922}{121}.}
   \ref{Norlund1}{N\"orlund,N.E. {\it Differenzenrechnung} (Springer--Verlag, 1924, Berlin.)}
  \ref{Duff}{Duff,G.F.D. \aom{56}{1952}{115}.}
  \ref{DandS}{Duff,G.F.D. and Spencer,D.C. \aom{45}{1951}{128}.}
  \ref{BGM}{Berger, M., Gauduchon, P. and Mazet, E. {\it Lect.Notes.Math.}
  {\bf 194} (1971) 1. }
  \ref{Patodi}{Patodi,V.K. \jdg{5}{1971}{233}.}
  \ref{GandS}{G\"unther,P. and Schimming,R. \jdg{12}{1977}{599}.}
  \ref{MandS}{McKean,H.P. and Singer,I.M. \jdg{1}{1967}{43}.}
  \ref{Conner}{Conner,P.E. {\it Mem.Am.Math.Soc.} {\bf 20} (1956).}
  \ref{Gilkey2}{Gilkey,P.B. \aim {15}{1975}{334}.}
  \ref{MandP}{Moss,I.G. and Poletti,S.J. \plb{333}{1994}{326}.}
  \ref{BKD}{Bordag,M., Kirsten,K. and Dowker,J.S. \cmp{182}{1996}{371}.}
  \ref{RandO}{Rubin,M.A. and Ordonez,C. \jmp{25}{1984}{2888}.}
  \ref{BaandD}{Balian,R. and Duplantier,B. \aop {112}{1978}{165}.}
  \ref{Kennedy2}{Kennedy,G. \aop{138}{1982}{353}.}
  \ref{DandKi2}{Dowker,J.S. and Kirsten, K. {\it Analysis and Appl.}
 {\bf 3} (2005) 45.}
  \ref{Dow40}{Dowker,J.S. \cqg{23}{2006}{1}.}
  \ref{BandHe}{Br\"uning,J. and Heintze,E. {\it Duke Math.J.} {\bf 51} (1984)
   959.}
  \ref{Dowl}{Dowker,J.S. {\it Functional determinants on M\"obius corners};
    Proceedings, `Quantum field theory under
    the influence of external conditions', 111-121,Leipzig 1995.}
  \ref{Dowqg}{Dowker,J.S. in {\it Quantum Gravity}, edited by
  S. C. Christensen (Hilger, Bristol, 1984).}
  \ref{Dowit}{Dowker,J.S. \jpa{11}{1978}{347}.}
  \ref{Kane}{Kane,R. {\it Reflection Groups and Invariant Theory} (Springer,
  New York, 2001).}
  \ref{Sturmfels}{Sturmfels,B. {\it Algorithms in Invariant Theory}
  (Springer, Vienna, 1993).}
  \ref{Bourbaki}{Bourbaki,N. {\it Groupes et Alg\`ebres de Lie}  Chap.III, IV
  (Hermann, Paris, 1968).}
  \ref{SandTy}{Schwarz,A.S. and Tyupkin, Yu.S. \np{242}{1984}{436}.}
  \ref{Reuter}{Reuter,M. \prD{37}{1988}{1456}.}
  \ref{EGH}{Eguchi,T. Gilkey,P.B. and Hanson,A.J. \prp{66}{1980}{213}.}
  \ref{DandCh}{Dowker,J.S. and Chang,Peter, \prD{46}{1992}{3458}.}
  \ref{APS}{Atiyah M., Patodi and Singer,I.\mpcps{77}{1975}{43}.}
  \ref{Donnelly}{Donnelly.H. {\it Indiana U. Math.J.} {\bf 27} (1978) 889.}
  \ref{Katase}{Katase,K. {\it Proc.Jap.Acad.} {\bf 57} (1981) 233.}
  \ref{Gilkey3}{Gilkey,P.B.\invm{76}{1984}{309}.}
  \ref{Degeratu}{Degeratu.A. {\it Eta--Invariants and Molien Series for
  Unimodular Groups}, Thesis MIT, 2001.}
  \ref{Seeley}{Seeley,R. \ijmp {A\bf18}{2003}{2197}.}
  \ref{Seeley2}{Seeley,R. .}
  \ref{melrose}{Melrose}
  \ref{DandW}{Douglas,R.G. and Wojciekowski,K.P. \cmp{142}{1991}{139}.}
  \ref{Dai}{Dai,X. \tams{354}{2001}{107}.}
\end{putreferences}

\bye